\documentclass[11pt]{amsart}
\usepackage{amsfonts,amscd,amsmath,amssymb,amsthm,mathrsfs}
\usepackage{array,latexsym,color,float,enumerate}
\usepackage{graphicx,epsfig}
\usepackage[top=4.05cm, bottom=4.05cm, left=3.5cm, right=3.5cm]{geometry}
\usepackage[all]{xy}
\usepackage{lmodern} %pdf
\numberwithin{equation}{section}

\theoremstyle{plain}
\newtheorem{theorem}{Theorem}[section]
\newtheorem{lemma}[theorem]{Lemma}
\newtheorem{proposition}[theorem]{Proposition}
\newtheorem{corrolary}[theorem]{Corollary}

\theoremstyle{definition}
\newtheorem{definition}[theorem]{Definition}
\newtheorem{example}[theorem]{Example}
\newtheorem{noname}[theorem]{}
\newtheorem{remark}[theorem]{Remark}
\newtheorem{construction}[theorem]{Construction}
\newtheorem{notation}[theorem]{Notation}

\theoremstyle{remark}
\newtheorem*{smallremark}{Remark}

\newtheorem{case}{Case} \makeatletter \@addtoreset{case}{theorem}\makeatother
\newtheorem{claim}{Claim} \makeatletter \@addtoreset{claim}{theorem}\makeatother

\newcommand{\bthm}{\begin{theorem}}                     \newcommand{\ethm}{\end{theorem}}
\newcommand{\bprop}{\begin{proposition}}                \newcommand{\eprop}{\end{proposition}}
\newcommand{\blem}{\begin{lemma}}                       \newcommand{\elem}{\end{lemma}}
\newcommand{\bcor}{\begin{corrolary}}                   \newcommand{\ecor}{\end{corrolary}}
\newcommand{\brem}{\begin{remark}}                      \newcommand{\erem}{\end{remark}}
\newcommand{\bdfn}{\begin{definition}}                  \newcommand{\edfn}{\end{definition}}
\newcommand{\bitem}{\begin{itemize}}                    \newcommand{\eitem}{\end{itemize}}
\newcommand{\bex}{\begin{example}}                      \newcommand{\eex}{\end{example}}
\newcommand{\bno}{\begin{noname}}                       \newcommand{\eno}{\end{noname}}
\newcommand{\bsrem}{\begin{smallremark}}                \newcommand{\esrem}{\end{smallremark}}
\newcommand{\bnot}{\begin{notation}}                    \newcommand{\enot}{\end{notation}}
\newcommand{\bcon}{\begin{construction}}                \newcommand{\econ}{\end{construction}}
\newcommand{\bca}{\begin{case}}                         \newcommand{\eca}{\end{case}}
\newcommand{\bcl}{\begin{claim}}                        \newcommand{\ecl}{\end{claim}}
\newcommand{\ba}{\begin{array}}                         \newcommand{\ea}{\end{array}}
\newcommand{\benum}{\begin{enumerate}}                  \newcommand{\eenum}{\end{enumerate}}

\newcommand{\un}{\underline}            \newcommand{\ov}{\overline}
              \newcommand{\wt}{\widetilde}
\newcommand{\cal}[1]{\mathcal{#1}}

\def\8{\infty}
\def\PP{\mathbb{P}}

\def\C{\mathbb{C}}
\def\Z{\mathbb{Z}}
\def\N{\mathbb{N}}
\def\Q{\mathbb{Q}}
\def\E{\widehat{E}}
\def\ovk{\overline\kappa}
\def\kp{\kappa}

\def\xra{\xrightarrow}
\def\embed{\hookrightarrow}

\def\:{\colon}

\def\med{\medskip}
\def\ssk{\smallskip}
\def\bsk{\bigskip}
\newcommand{\noin}{\noindent}

\def\Bk{\operatorname{Bk}}
\def\Spec{\operatorname{Spec}}

\def\Sing{\operatorname{Sing}}

\def\Supp{\operatorname{Supp}}
\def\Pic{\operatorname{Pic}}

\def\Exc{\operatorname{Exc}}

\usepackage[linkbordercolor={0 0 1}]{hyperref}

\begin{document}
\title[Singular $\Q$-homology planes]{Singular $\Q$-homology planes of negative Kodaira dimension have smooth locus of non-general type}

\author[Karol Palka]{Karol Palka}\address{Karol Palka: Institute of Mathematics, University of Warsaw, ul. Banacha 2, 02-097 Warsaw, Poland}\email{palka@impan.pl}

\author[Mariusz Koras]{Mariusz Koras}\address{Mariusz Koras: Institute of Mathematics, University of Warsaw, ul. Banacha 2, 02-097 Warsaw, Poland}\email{koras@mimuw.edu.pl}

\thanks{Both authors were supported by Polish Grant NCN}


\begin{abstract} \noindent We show that if a normal $\Q$-acyclic complex surface has negative Kodaira dimension then its smooth locus is not of general type. This generalizes an earlier result of Koras-Russell for contractible surfaces. \end{abstract}

\maketitle

\section{Main result}

We work in the category of complex algebraic varieties. In this paper we continue the program of classification of singular $\Q$-homology planes. A normal surface $S'$ is called a \emph{$\Q$-homology plane} if its rational cohomology is the same as that of the affine plane $\C^2$, i.e. $H^*(S',\Q)\cong \Q$. Properties of these surfaces have been analyzed for a long time, motivations come from studies on the Cancellation Conjecture of Zariski, on the two-dimensional Jacobian Conjecture, on quotients of actions of reductive groups on affine spaces or on exotic $\C^n$'s. For a review in the smooth case see \cite[\S 3.4]{Miyan-OpenSurf} and in the singular case \cite{Palka-review}. The complete counterparts of $\Q$-homology planes are complex surfaces with rational cohomology of $\PP^2$, so-called fake projective planes (they are algebraic by \cite[V.1.1]{BHPV}). These are well understood, for example it has been shown recently that there are exactly $100$ of them  up to biholomorphism, hence up to algebraic isomorphism, cf. \cite{fake_proj_spaces}. While the world of $\Q$-homology planes is also relatively well understood, in some cases still much needs to be done.

The basic invariants of $S'$ are the logarithmic Kodaira dimension $\ovk(S')$ and the \mbox{logarithmic} Kodaira dimension of the smooth locus, $\ovk(S'\setminus\Sing S')$. They take \mbox{values} in $\{-\8,0,1,2\}$ and satisfy the inequality $\ovk(S'\setminus\Sing S')\geq \ovk(S')$ (see \cite{Iitaka} for the definition and properties of the logarithmic Kodaira dimension $\ovk$).

In \cite{Palka-classification} the first author obtained a classification of singular $\Q$-homology planes with smooth locus of non-general type, i.e. with $\ovk(S'\setminus\Sing S')\leq 1$. Earlier a structure theorem for singular $\Q$-homology planes with singularities of at most quotient type, which have Kodaira dimension of the smooth locus equal to $1$ or $-\8$, has been given in \cite{MiSu-hPlanes}. We therefore concentrate on the remaining cases of the classification, i.e. when the smooth locus is a surface of general type. These are more difficult to tackle because of a lack of structure theorems for surfaces of (log-)general type. The only obviously available tool here is the logarithmic Bogomolov-Miyaoka-Yau inequality proved by Kobayashi \cite{Kob}. An important consequence of the latter (see for example \cite[3.3]{Palka-classification}) is that $S'$ has only one singular point and this point is of analytical type $\C^2/G$ for some finite subgroup $G<GL(2,\C)$. By a theorem of Pradeep-Shastri \cite{PS-rationality} $S'$ is rational. Many such surfaces have been constructed (see for example \cite[Theorem 1]{MiSu-Qhp_with_C**}). While a priori there is no bound on the Kodaira dimension of $S'$, we obtain in this paper the following result.

\bthm\label{thm:main result} Singular $\Q$-homology planes of negative Kodaira dimension have smooth locus of non-general type. \ethm

The theorem is a generalization of an analogous result of Koras-Russell \cite{KR-ContrSurf} on contractible surfaces and their earlier result on contractible surfaces with hyperbolic action of $\C^*$, which was a crucial step in the proof of linearizability of $\C^*$-actions (and hence actions of connected reductive groups) on $\C^3$, cf. \cite{KR-C*onC3}. Although we do not assume here that $S'$ is contractible, but only that it is $\Q$-acyclic, a significant part of methods from loc. cit. can be adapted to our situation. The result for contractible surfaces is recovered as a special case. We assume that $\ovk(S')=-\8$ and $\ovk(S'\setminus\Sing S')=2$ and analyze the consequences. The final contradiction is obtained in a long series of steps restricting more and more the possible geometry and derived numerical properties of the boundary and of the exceptional divisor of the resolution.

We now give a more detailed overview. In section \ref{sec:notation} we introduce basic notions and notation recalling in particular properties of the Fujita-Zariski decomposition and the notion of Hamburger-Noether pairs. In section \ref{sec:topology and inequalities} we describe homological and geometric properties of a $\Q$-homology plane $S'$, of its minimal resolution $S$ and its smooth locus $S_0$. Basic properties of the snc-minimal boundary $D$, the exceptional divisor $\E$ of the minimal resolution and of the logarithmic canonical divisor $K+D+\E$, where $K$ is a canonical divisor on a smooth minimal completion $(\ov S,D+\E)$ of $S_0$, are derived. In particular, $\E$ and $D$ are connected trees and $\E$ has at most one branching component. In the whole paper the fact that $S'$ does not contain curves which are topologically contractible is essential. We decompose the divisor $\E$ as $\E=E+\Delta$, where $\Delta$ consists of 'external' $(-2)$-curves. We define $\epsilon$ by the equality $(K+D+\E)^2=-1-\epsilon$. By an inequality of Miyaoka $\epsilon\geq 0$. In section \ref{sec:bounding exc divisor} we show that except one case the inequality $K\cdot E+2\epsilon\leq 5$ holds (cf. \ref{1prop:KE and eps}), which gives crucial bounds on $K\cdot E$ and $\epsilon$ and allows us to list possible dual graphs of $\E$ (cf. \ref{1prop:possible $E$}). We also show that one can find an affine ruling of $S$ for which $\Delta$ is contained in fibers. It is proved in section \ref{sec:pre-minimal rulings} that if $\E$ has only one component which is not a $(-2)$-curve then the process of resolving the base point of this ruling on $\ov S$ can be well controlled. In section \ref{sec:D is a fork} we show that the boundary $D$ has only one branching component $B$ and this leads to a precise description of the Fujita-Zariski decomposition of $K+D+\E$. With these tools in hand we show in section \ref{sec:surface W} that the surface $\ov S-(\E+D-B)$, which contains $S_0$, is of general type. This takes considerable amount of work, but bounds possible shapes of $\E$ to four cases (cf. \ref{3cor:final bounds for E^}) by the logarithmic Bogomolov-Miyaoka-Yau inequality. These are finally excluded in section \ref{sec:special cases} by analyzing properties of the affine ruling of $S\setminus\Delta$. In sections \ref{sec:surface W} and \ref{sec:special cases} we need to support our analysis by referring to results of computer programs.

\bsk\textsl{\textsf{Acknowledgements.}} The theorem was obtained in the Ph.D. thesis of the first author during his graduate studies at the University of Warsaw under the supervision of the second author. We present here a modified version of the proof. Both authors would like to thank the Institute of Mathematics of the University of Warsaw for very good working conditions. 

\tableofcontents

\section{\label{sec:notation}Notation and preliminaries}

We use standard notions and notation of the theory of open algebraic surfaces, we recall some of them. The reader is referred to \cite{Miyan-OpenSurf} for a detailed treatment as well as for basic theorems of the theory. We denote the linear and numerical equivalences of divisors by $\sim$ and $\equiv$ respectively.

Let $T$ be a divisor with simple normal crossings on a smooth complete surface. We write $\un T$ for the reduced divisor with the same support. If $U$ is a component of $T$ then $\beta_T(U)=U\cdot (\un T-U)$ is called the \emph{branching number of $U$ in $T$} and any $U$ with $\beta_T\geq 3$ is called a \emph{branching component} of $T$. If $T$ is reduced and its dual graph contains no loops then we say that $T$ is \emph{a forest}, it is \emph{a tree} if it is connected. A component with $\beta_T\leq 1$ is called a \emph{tip of $T$}. The dual graph of $T$ is weighted, the weights of vertices are the self-intersections of the corresponding components of $T$. We define the discriminant $d(T)$ as equal to $1$ if $T=\emptyset$ and as the determinant of the minus intersection matrix of $T$ otherwise. By elementary expansion properties of determinants we have:

\blem\label{1lem:d() expansion} Let $C$ be a component of a rational tree $R$, let $R_1,\ldots,R_k$ be the connected components of $R-C$. Let $C_i$ be the irreducible component of $R_i$ meeting $C$. Then
$$d(R)=-C^2\prod_i d(R_i)-\sum_i d(R_i-C_i)\prod_{j\neq i} d(R_j).$$  \elem

Suppose $T$ is a (reduced) rational chain, i.e. it can be written as $T=T_1+\ldots+T_n$, where $T_i\cong \PP^1$, $\beta_T(T_i)\leq 2$ and $T_i\cdot T_{i+1}=1$ for $i=1,\ldots,n-1$. There are at most two choices of the first component of a chain, each defines an orientation on it. We write $T=[-T_1^2,\ldots,-T_n^2]$ and by $T^t$ we mean the same chain considered with an opposite orientation (there is only one orientation if $n=1$). We define $d'(T)=d(T-T_1)$ and we put $d'(\emptyset)=0$. In case $T_1^2=\ldots=T_n^2=-2$ we write $T=[(n)]$. We call $T$ \emph{admissible} if $T_i^2\leq -2$ for each $i$. If $d(T)\neq 0$ we define $$\delta(T)=\frac{1}{d(T)}, e(T)=\frac{d'(T)}{d(T)} \text{\ and\ } \wt e(T)=e(T^t).$$

Suppose $T$ is a tree with exactly one branching component $T_0$. Then $T$ is called a \emph{wide fork} and is called \emph{a fork} if $\beta_T(T_0)=3$. The fork $T$ is \emph{admissible} if it is rational, the three connected components of $T-T_0$ are admissible chains and the intersection matrix of $T$ is negative definite. Admissible chains and forks are exactly the exceptional snc-divisors of minimal resolutions of quotient singular points. A singular point on a surface is of quotient type if and only if locally analytically it is isomorphic to the singular point of $\C^2/G$ for some finite subgroup $G<GL(2,\C)$.

A \emph{normal pair} $(X,D)$ consists of a complete normal surface $X$ and a reduced simple normal crossing divisor $D$, whose support is contained in the smooth locus of $X$. If $X$ is smooth then $(X,D)$ is a \emph{smooth pair}. An \emph{$n$-curve} is a smooth rational curve with self-intersection $n$. If $D$ contains no non-branching $(-1)$-curves then the pair $(X,D)$ is snc-minimal. If $X_0$ is a normal (smooth) surface then any normal pair $(X,D)$, such that $X\setminus D=X_0$ is called a \emph{normal (smooth) completion} of $X_0$. If $(X,D)$ is a normal pair then a blow-up of $X$ with center $c\in D$ is called sprouting (subdivisional) for $D$ if $c$ belongs to exactly one (two) irreducible component of $D$.

Let $(X,D)$ be a smooth pair. Denote the canonical divisor on $X$ by $K_X$. If $\sigma\:Y\to X$ is a blow-up we denote its exceptional divisor by $\Exc \alpha$, the total transform, the reduced total transform and the proper transform of $D$ by $\sigma^*D$, $\sigma^{-1}D$, $\sigma'D$ respectively. We need the following easy observations.

\blem\label{1lem:blow-up and transforms} Let $(X,D)$ be a smooth pair and let $\sigma\:Y\to X$ be a blow-up. \benum[(i)]

\item If $A,B$ are divisors on $X$ then $A\cdot B=\sigma'A\cdot \sigma^*B=\sigma^*A\cdot \sigma^*B.$

\item If $\sigma$ is sprouting with respect to $D$ or if $D=0$ then $\sigma^*(K_X+D)=K_Y+\sigma^{-1}D-\Exc \sigma$ and $K_X\cdot (K_X+D)=K_Y\cdot (K_Y+\sigma^{-1}D)+1$.

\item If $\sigma$ is subdivisional with respect to $D$ then $\sigma^*(K_X+D)=K_Y+\sigma^{-1}D$ and $K_X\cdot (K_X+D)=K_Y\cdot (K_Y+\sigma^{-1}D)$. \eenum\elem

To compute the negative part of the Zariski decomposition of the logarithmic canonical divisor $K_X+D$ it is useful to compute the \emph{bark of $D$} ($\Bk D$). Barks are defined independently for all connected components of $D$, so in what follows we will assume that $D$ is connected. If $D$ is an admissible chain or an admissible fork we define $\Bk D$ as a unique $\Q$-divisor with support in $\Supp D$ satisfying $$(K_X+D-\Bk D)\cdot D_i=0$$ for each component $D_i$ of $D$. If $D=T=T_1+\ldots+T_n$ is an admissible chain then it is also convenient to define a 'one-sided bark' $Bk(T,T_1)$ with support contained in $\Supp T$ by $$T_i\cdot \Bk (T,T_1)=-\delta_{i,1}$$ (Kronecker's delta). If in the last case the choice of $T_1$ is clear from the context we write $\Bk' T$ for $\Bk(T,T_1)$. Note that $\Bk T= \Bk(T,T_1)+\Bk(T,T_n)$.

To define the bark in general we need some additional notions. Suppose $D$ is not a chain. A chain $T\subseteq D$ is a \emph{twig} of $D$ if $\beta_D\leq 2$ for all components of $T$ and $\beta_D=1$ for some (unique in fact) component of $T$. If $T$ is a twig of $D$ then by a \emph{default orientation} of $T$ we mean the one in which the tip of $D$ contained in $T$ is the first component ($T_1$) of $T$. Analogously, if $D$ is not an admissible chain (it may or may not be a chain) we define admissible twigs and maximal admissible twigs of $D$.

Suppose now $D$ is neither an admissible chain nor an admissible fork. Let $R_1,\ldots, R_s$ be all the maximal admissible twigs of $D$. We define  $$\Bk D=\Bk' R_1+\ldots+\Bk 'R_s.$$ We put $D^{\#}=D-\Bk D$,  $$\delta(D)=\sum_{i=1}^s \delta(R_i),\ e(D)=\sum_{i=1}^s e(R_i)\text{\ and\ } \wt e(D)=\sum_{i=1}^s \wt e(R_i).$$ We will need the following properties of barks, most of which follow by a straightforward calculation (cf. \cite[\S2.3]{Miyan-OpenSurf}).

\blem\label{1lem:ChainBk} Let $T=T_1+\ldots+T_n$ be an admissible chain, write $\Bk' T=\sum_{i=1}^n m_i' T_i$ and $\Bk T=\sum_{i=1}^n m_iT_i$, then:\benum[(i)]

\item $d'(T)\leq d(T)-1$, $e(T)=(-T_1^2-e(T-T_1))^{-1}$, $\delta(T)\leq e(T)\leq 1-\delta(T)$,

\item $m_i'=d(T_{i+1}+\ldots+T_n)/d(T)$,

\item $0<m_i'<1$ and $0<m_i\leq1$ (in particular $\Supp \Bk' T=\Supp \Bk T=\Supp T$). Moreover, if $m_i=1$ for some $i$ then $T=[2,2,\ldots,2]$ and $m_i=1$ for each $i$,

\item $\Bk'^2 T=-e(T)$, $\Bk^2 T=-e(T)-\wt e(T)-2\delta(T)=-(d'(T)+d'(T^t)+2)/d(T)\geq -2$. \eenum\elem

\bsrem The formula $e(T)=(-T_1^2-e(T-T_1))^{-1}$ shows that knowing $e(T)$ one can recover $T$ in terms of continued fractions. \esrem

\blem\label{1lem:ForkBk} Let $F=B+R_1+R_2+R_3$ be an admissible fork with maximal twigs $R_i$. Write $\Bk F=\sum_{i=1}^n m_i F_i$, where $F_i$ are the irreducible distinct components of $F$. Then:\benum[(i)]

\item $0<m_i\leq1$ (in particular $\Supp \Bk F=\Supp F$). Moreover, if $m_i=1$ for some $i$ then $F$ consists of $(-2)$-curves and $m_i=1$ for each $i$,

\item $(d(R_1),d(R_2),d(R_3))$ is one of the \emph{Platonic triples}: $(2,3,3)$, $(2,3,4)$, $(2,3,5)$ or $(2,2,k)$ for some $k\geq2$,

\item $1<\delta(F)\leq \wt e(F)<2\leq-B^2$,

\item $d(F)=d(R_1)d(R_2)d(R_3)(-B^2-\wt e(F))$,

\item $\Bk^2 F=-(\delta(F)-1)^2(-B^2-\wt e(F))^{-1}-e(F)<-e(F)<-1$. \eenum\elem

\brem\label{1rem:-2 chains and forks} Note that since $\wt e(T)+\delta(T)\leq 1$ (and $e(T)+\delta(T)\leq 1$ too) for an admissible chain $T$, we have $\Bk^2 T=-2$ if and only if $T$ consists of $(-2)$-curves. Then for an admissible fork $F$ we get by \ref{1lem:ForkBk}(iii) that $\delta(F)+\wt e(F)\leq 3\leq 1-B^2$, so $-\Bk^2 F\leq \delta(F)-1+e(F)\leq 2$ and again the equality occurs if and only if $F$ consists of $(-2)$-curves (is a $(-2)$-fork).\erem

\blem\label{1lem:SmallChains} For every $d>2$ there exist at least two admissible chains with discriminant $d$: $[d]$ and $[(d-1)]$. Here is a full list of all other admissible chains for $d\leq11$:\benum[]

\item $d=5$ \emph{:} $[3,2]$,

\item $d=7$ \emph{:} $[4,2]$, $[3,(2)]$,

\item $d=8$ \emph{:} $[3,3]$, $[2,3,2]$,

\item $d=9$ \emph{:} $[5,2]$, $[3,(3)]$,

\item $d=10$ \emph{:} $[4,(2)]$,

\item $d=11$ \emph{:} $[6,2]$, $[4,3]$, $[3,(4)]$, $[2,3,(2)]$.  \eenum\elem

A $\PP^1$-ruling of a complete normal surface is a surjective morphism of the surface onto a smooth curve, for which general fibers are isomorphic to $\PP^1$. Let $(X,D)$ be a smooth pair and let $p\:X\to \PP^1$ be a $\PP^1$-ruling. The multiplicity of an irreducible component $L$ of a fiber will be denoted by $\mu(L)$ (see \cite[2.9]{Palka-classification} for summary of basic properties of singular fibers which we use here). The horizontal part $D_h$ of $D$ is defined as an effective divisor with support in $\Supp D$, such that $D-D_h$ is effective and intersects trivially with fibers. A horizontal irreducible curve $C$ is called an $n$-section of $p$ (or simply 'section' if $n=1$) if $C\cdot F=n$ for any fiber $F$ of $p$. The components of any fiber $F$ are either $D$-components (the ones contained in $D$) or $(X-D)$ -components. We denote the number of $(X-D)$ -components of $F$ by $\sigma(F)$, by $\nu$ the number of fibers with $\sigma=0$ (which are contained in $D$) and by $\Sigma_{X-D}$ the sum of numbers $(\sigma(F)-1)$ taken over the set of fibers not contained in $D$. Of course, for a general fiber $\sigma=1$. Put $h=\#D_h$. The basic observation is that if one contracts a vertical $(-1)$-curve and simultaneously changes $(X,D)$ for its image then the numbers $b_2(X)-b_2(D)-\Sigma+\nu$ and $h$ do not change. So since for a $\PP^1$-bundle over a smooth complete curve $b_2(D)=h+\nu$, $b_2(X)=2$ and $\Sigma=0$, we get the following relation (cf. \cite[4.16]{Fujita}).

\bprop\label{1prop:Fujita equation} If $(X,D)$ is a smooth pair then for any $\PP^1$-ruling of $X$ $$\Sigma_{X-D}=h+\nu-2+b_2(X)-b_2(D).$$ \eprop

Any $0$-curve on a smooth surface induces a $\PP^1$-ruling with this curve as one of the fibers (cf. \cite[2.8]{Palka-classification}). The structure of singular fibers of such rulings is well known (cf. 2.9 loc. cit.).

\bdfn\label{1def:ruling completion} A \emph{rational ruling} of a surface is a surjective morphism of the surface onto a smooth curve, for which general fibers are rational curves.  If $p_0\: X_0\to B_0$ is a rational ruling of a normal surface then by a \emph{completion of $p_0$} we mean a triple $(X,D,p)$, where $(X,D)$ is a normal completion of $X_0$ and $p\:X\to B$ is an extension of $p_0$ to a $\PP^1$-ruling with $B$ being a smooth completion of $B_0$. We say that $p$ is a \emph{minimal completion of $p_0$} if $p$ does not dominate any other completion of $p_0$. \edfn

Note that if $p$ is a minimal completion of $p_0$ then every vertical $(-1)$-curve contained in $D$ intersects at least three other components of $D$.

\ssk We recall the notion of Hamburger-Noether pairs. For details see \cite{Russell2} and \cite[Appendix]{KR-C*onC3}.

\bdfn\label{1def:char pairs} Suppose we are given an irreducible germ of a singular analytic curve $(\chi_1,q_1)$ on a smooth algebraic surface and a curve $C_1$ passing through $q_1$, smooth at $q_1$. Put $c_1=(C_1\cdot \chi_1)_{q_1}$ and choose a local coordinate $y_1$ at $q_1$ in such a way that $Y_1=\{y_1=0\}$ is transversal to $C_1$ at $q_1$ and $p_1=(Y_1\cdot\chi_1)_{q_1}$ is not bigger than $c_1$. Blow up over $q_1$ until the proper transform $\chi_2$ of $\chi_1$ meets the reduced total inverse image $F_1$ of $C_1$ in a point $q_2$, which does not belong to components of $F_1$ other than the unique exceptional component $C_2$ of $\un{F}_1-C_1$. We then say that \emph{$C_2$ (and $F_1$) is produced from $C_1$ by the pair $\binom{c_1}{p_1}$}. Put $c_2=(C_2\cdot \chi_2)_{q_2}$. We repeat this procedure and define successively $(\chi_i,q_i)$ and $C_i$ until $\chi_{h+1}$ is smooth for some $h\geq 1$. Then we refer to the sequence $\binom{c_1}{p_1}, \binom{c_2}{p_2}, \ldots, \binom{c_h}{p_h}$ as the sequence of \emph{Hamburger-Noether pairs} (or \emph{characteristic pairs} for short) \emph{of the resolution of $(\chi_1,q_1)$} or the sequence of \emph{characteristic pairs of $F$}, where $F$ is the (reduced) total transform of $C_1$. It is convenient to extend the definition to the case when $(\chi_1,q_1)$ is smooth by defining it sequence of characteristic pairs to be $\binom{1}{0}$. \edfn

The convention that $c_i\geq p_i$ seems artificial, but will be useful in our situation. Note also that the definitions make sense for $(\chi_1,q_1)$ reducible, as long as each blow-up (except possibly the last one) leaves irreducible branches of $\chi_1$ unsplitted, so that the center of the succeeding blow-up is uniquely determined.

\blem\label{1lem:sum of mu and squares} Assume that the sequence of blow-ups $(\sigma_j)_{j\in I_i}$, leading from $(\chi_i,q_i)$ to $(\chi_{i+1},q_{i+1})$ is described as above by the characteristic pair $\binom{c_i}{p_i}$. Let $\mu_j$ be the multiplicity of the center of $\sigma_j$. Then we have:\benum[(i)]

\item $c_{i+1}=gcd(c_i,p_i)$,

\item  $\sum_{I_i} \mu_j=c_i+p_i-gcd(c_i,p_i)$,

\item $\sum_{I_i}\mu_j^2=c_ip_i$. \eenum \elem

\begin{proof}The formulas hold in case $c_i=p_i$. If $c_i>p_i$ then perform the first blow-up and note that the remaining part of the sequence $(\sigma_j)_{j\in I_i}$ is described by $\binom{c_i-p_i}{p_i}$ in case $c_i-p_i\geq p_i$ or by $\binom{p_i}{c_i-p_i}$ otherwise. The multiplicity of the first center is $p$. Now the result follows by induction on $max(c_i,p_i)$.\end{proof}

Consider a fiber $F$ of a $\PP^1$-ruling of some smooth complete surface, such that $F$ contains at most one $(-1)$-curve. Suppose $U$ is a component of $F$ with $\mu_F(U)=1$. There is a uniquely determined sequence of contractions of $(-1)$-curves in $F$ and its subsequent images which makes $F$ a smooth $0$-curve and does not contract $U$. The reverting sequence of blow-ups orders naturally the set of components of $F$ in order they are produced. Let $B_1, \ldots, B_k$ be the branching components of $F$ ordered as described. We call the chain consisting of $U$, the components produced before $B_1$ and of $B_1$ \emph{the first branch of $F$}, the chain consisting of components produced after $B_1$ but before $B_2$ and of $B_2$ \emph{the second branch of $F$}, etc. The \emph{(k+1)-st branch} is a chain of components produced after $B_k$.

\bdfn\label{1def:char pairs for fibers} Let $F$ and $U$ be as above. Denote the birational transform of $U$ after contractions (the image of $F$) by the same letter. If $F$ is singular let $L$ be the $(-1)$-curve of $F$. For some $q\in L$ let $(\chi,q)$ be an irreducible germ of a smooth analytic curve intersecting $L$ transversally at $q$. Denote its image after contractions by $(\chi_1,q_1)$. Then the sequence of characteristic pairs of the resolution of $(\chi_1,q_1)$ produces $L$ (and $F$) from $U$ (cf. \ref{1def:char pairs}). If the choice of $U$ is clear from the context we refer to this sequence as \emph{the sequence of characteristic pairs of $F$}.\edfn

Note that by definition if $\binom{c_i}{p_i}$, $i=1,\ldots,h$ is the sequence of characteristic pairs of $F$ then $gcd(c_h,p_h)=1$ and the last curve produced by the sequence  (the unique $(-1)$-curve in case $F$ is singular) has multiplicity $c_1$. As in \ref{1def:char pairs} the sequence of characteristic pairs of a smooth fiber is $\binom{c_1}{p_1}=\binom{1}{0}$.

\bex Consider a $\PP^1$-ruling of some complete surface. Let the notation be as above. Let $F=A_n+\ldots+A_1+L+B_1+\ldots+B_m$ be a non-branched singular fiber with a unique $(-1)$-curve $L$. Only the tips of $F$, $A_n$ and $B_m$, have multiplicity one. $F$ is produced from $A_n$ by one characteristic pair, call it $\binom{c}{p}$ (we have $gcd(c,p)=(\chi\cdot L)_q=1$). The algorithm to recover $F$ when $\binom{c}{p}$ is known reduces to some simple observations. Let $C_1$ be the birational transform of $A_n$ after contraction of remaining components of the fiber. We have $c=(C_1\cdot \chi_1)_{q_1}$ and $p=(Y_1\cdot \chi_1)_{q_1}$. Consider a blow-up at $q_1$, let $E$ be the exceptional curve and let $(\chi',q')$, $q'\in E$ be the proper transform of $(\chi_1,q_1)$. If $c=p$ then $q'$ does not belong to $C_1+Y_1$ and we are done. If $c>p$ then $q'\in C_1$, $(C_1\cdot \chi')_{q'}=c-p$ and $(E\cdot \chi')_{q'}=p$. In case $c-p\geq p$ we continue with the pair $\binom{c-p}{p}$ and with $(C_1,E,\chi')$ replacing $(C_1,Y_1,\chi_1)$. In case $c-p<p$ we continue with the pair $\binom{p}{c-p}$ and with $(E,C_1,\chi')$ replacing $(C_1,Y_1,\chi_1)$. Put $A=A_n+\ldots+A_1$. One proves that $$c=d(A)\text{\ and\ } p=d'(A).$$ Here are some examples. If $F=[k,1,(k-1)]$ then $\binom{c}{p}=\binom{k}{1}$. If $F=[(k-1),1,k]$ then $\binom{c}{p}=\binom{k}{k-1}$. If $F=[5,3,1,2,3,(3)]$ then $\binom{c}{p}=\binom{14}{3}$. \eex

\blem\label{1lem:Eff-NegDef is Eff} Let $A$ and $B$ be $\Q$-divisors on a smooth complete surface, such that $Q(B)$ is negative definite and $A\cdot B_i\leq 0$ for each irreducible component $B_i$ of $B$. Denote the integral part of a $\Q$-divisor by $[\ ]$.\benum[(i)]

\item If $A+B$ is effective then $A$ is effective.

\item If $n\in \N$ and $n(A+B)$ is a $\Z$-divisor then $h^0(n(A+B))=h^0([nA])$.\eenum\elem

\begin{proof} See 2.2 loc. cit. \end{proof}

For any divisor $D$ on a smooth complete surface $X$ we define the arithmetic genus of $D$ by Put $p_a(D)=\frac{1}{2}D\cdot(K_X+D)+1$. We have $p_a(D_1+D_2)=p_a(D_1)+p_a(D_2)+D_1\cdot D_2-1$. One shows by induction that if $D$ is a rational reduced snc-tree then $p_a(D)=0$. For the notion and properties of the Kodaira dimension of a divisor see \cite{Iitaka}.

\blem\label{1lem:linear systems} Let $D$ be an effective divisor on a complete smooth rational surface $X$.  \benum[(i)]

\item We have $h^0(K_X+D)+h^0(-D)\geq p_a(D)$. If $|K_X+D|=\emptyset$ then $\underline{D}$ is a rational snc-forest and if moreover $D=D_1+D_2$ with $p_a(D_1)=p_a(D_2)=0$ then $D_1\cdot D_2\leq 1$.

\item If $D$ has smooth rational components and $X$ in neither a Hirzebruch surface nor $\PP^2$ then $D\sim \sum C_i$, where $C_i\cong \PP^1$ and $C_i^2\leq -1$.

\item If $\kappa(K_X+D)=-\8$ then for any divisor $F$ one has $\kappa(F+m(K_X+D))=-\8$ for $m\gg 0$.

\eenum\elem

\begin{proof} (i) The Riemann-Roch theorem on a rational surface gives $h^0(K_X+D)+h^0(-D)\geq p_a(D)$ and the other properties follow by applying it in various ways (cf. \cite[2.1, 2.2]{Russell2}), for (ii) see \cite[4.1]{KR-C*onC3}, for (iii) see \cite[2.5]{Fuj-Zar}. \end{proof}

One of the fundamental facts used in this paper is the inequality of Bogomolov-Miyaoka-Yau type proved by Kobayashi (\cite{Kob}). It is most convenient for us to refer to the following corollary from a generalization proved by Langer (see \cite[5.2]{Langer} for the generalization and \cite[2.5]{Palka-k(S_0)=0} for the proof of the proposition).

\bprop\label{1lem:KobIneq} Let $(X,D)$ be a smooth pair with $\kp(K_X+D)\geq 0$. \benum[(i)]

\item The following inequality holds: $$3\chi(X-D)+\frac{1}{4}((K_X+D)^-)^2\geq (K_X+D)^2.$$

\item For each connected component of $D$, which is a connected component of $\Bk D$ (hence contractible to a quotient singularity) denote by $G_P$ the local fundamental group of the respective singular point $P$, put $D^\#=D-\Bk D$. Then $$\chi(X-D)+\sum_P\frac{1}{|G_P|}\geq \frac{1}{3} (K_X+D^\#)^2.$$ \eenum \eprop

\med\section{\label{sec:topology and inequalities}Basic properties and some inequalities}

Let $S'$ be a complex $\Q$-homology plane, i.e. a normal complex algebraic surface, such that $H^*(S',\Q)\cong \Q$. We assume that $S'$ is singular. We denote by $\rho\:S\to S'$ be the snc-minimal resolution of singularities and by $\E$ be the reduced exceptional divisor of $\rho$. In the whole paper we assume for a contradiction that $\ovk(S')=-\8$ and $\ovk(S_0)=2$ and we derive consequences. Since $\ovk(S_0)=2$, $S_0$ is neither affine- nor $\C^*$-ruled, so it admits a unique snc-minimal completion $(\ov S,D+\E)$ (cf. \cite[6.11]{Palka-classification}).

Following \cite[3.3]{Palka-k(S_0)=0} we call a curve $C$ on $(\ov S,D+\E)$ \emph{simple} if and only if $C\cong \PP^1$ and $C$ has at most one common point with each connected component of $D+\E$. Once we know that $S'$ is affine we get that $C$ on $(\ov S,D+\E)$ is simple if and only if $\rho(C\cap S)$ is topologically contractible. Decompose $\E$ as $\E=E+\Delta$, where $\Delta$ is the divisor of external $(-2)$-curves in $\E$, i.e. $\Delta$ is a reduced divisor with the smallest support, such that $E$ does not contain a $(-2)$-tip.

Let us first collect some basic results, mainly following from \cite{Palka-classification}. We note that besides affiness all properties of $S'$ stated in (i) and (ii) are based on the Kobayashi inequality \ref{1lem:KobIneq}(ii). For open surfaces and for smooth pairs we have a notion of minimality called \emph{almost minimality}, which generalizes the notion of minimality for complete smooth surfaces, we refer to \cite[2.3.11]{Miyan-OpenSurf} for the definition. We use the fact that for almost minimal pairs the Zariski decomposition of the logarithmic canonical divisor can be computed in terms of barks. We denote the canonical divisor of $\ov S$ by $K$.

\bprop\label{1prop:basic S' properties} With the notation as above one has: \benum[(i)]

\item $S'$ is affine, rational and its singular locus consists of one singular point of quotient type,

\item there is no simple curve on $(\ov S,D+\E)$, in particular the pair $(\ov S,D+\E)$ is almost minimal and $(K+D+\E)^-=\Bk D+\Bk \E$,

\item not every component of $\E$ is a $(-2)$-curve, i.e. $\E\neq\Delta$,

\item $d(D)=-d(\E)\cdot |H_1(S',\Z)|^2$, $\pi_1(S')=\pi_1(S)$ and $H_i(S',\Z)=0$ for $i>1$,

\item $D$ is a rational tree and if it has a component with non-negative self-intersection then this component is branching and $D$ is not a fork,

\item the inclusion $D\cup \E\to \ov S$ induces an isomorphism on $H_2(-,\Q)$,

\item $\Sigma_{S_0}=h+\nu-2$ and $\nu\leq 1$,

\item $\Pic S_0\cong H_1(S_0,\Z)$ is of order $d(\E)\cdot |H_1(S',\Z)|$. \eenum\eprop

\begin{proof} (i) $S'$ is affine and logarithmic by \cite[3.2, 3.3]{Palka-classification}, so it is rational by \cite{PS-rationality}. (ii) The non-existence of simple curves is proved for example in \cite[3.4]{Palka-k(S_0)=0} (or one can refer to the nonexistence of contractible curves on $S'$, see \cite{GM-Affine-lines}). Then $(\ov S,D+\E)$ is almost minimal and $(K+D+\E)^-=\Bk D+\Bk \E$ by \cite[2.3.15]{Miyan-OpenSurf} and by the uniqueness of the Zariski decomposition. (iii) If $\E=\Delta$ then $(K+D)\cdot\E=0$, so since $\ovk(S_0)\geq 0$ and since $\E$ has negative definite intersection matrix, $\kappa(K+D)\geq 0$ by \ref{1lem:Eff-NegDef is Eff}, a contradiction. For (iv), (vi)-(viii) see \cite[3.1, 3.2]{Palka-classification}.

(v) Since $S'$ is affine, $D$ is connected, so it is a rational tree by 3.4 loc. cit. Let $B$ be a component of $D$ with $B^2\geq 0$. We blow up over $B$ until $B^2=0$. Let $(\wt S,\wt D)\to (\ov S, D)$ be the resulting birational morphism. We can choose the centers of subsequent blow-ups so that $\wt D$ contains at most one non-branching $(-1)$-curve and, unless $D=B$, so that the blow-ups are subdivisional for $D$ and its total transforms. In any case it follows that $B$ has to be a branching component $(\beta_D(B)\geq 3)$, otherwise we get a $\PP^1$-, a $\C^1$- or a $\C^*$- ruling of $S_0$, hence $\ovk(S_0)\leq 1$ by Iitaka's addition theorem (cf. \cite[10.4]{Iitaka}), which is a contradiction. Suppose now that $D$ is a fork and $B$ is its unique branching component. Then $B$ gives a $\PP^1$-ruling of $\wt S$ for which $\wt D_h$ consists of three sections. By \ref{1prop:basic S' properties}(vii) we have $\Sigma_{S_0}=2$, because $\E$ is vertical. Note that every vertical $(-1)$-curve is an $S_0$-component. Suppose there is a singular fiber $F$ containing a unique $(-1)$-curve $L$. We have $\mu(L)>1$, so $\wt D_h$ does not intersect $L$. However, $\un F-L$ has at most two connected components, so $\wt D$ contains a loop, a contradiction. Thus every singular fiber has at least two $(-1)$-curves. Denote the fiber containing $\E$ by $F_0$. Let $D_0$ be the divisor of $\wt D$-components of $F_0$ and let $L_1$, $L_2$ be some $(-1)$-curves in $F_0$. We have $D_0\neq 0$, otherwise one of the $S_0$-components of $F_0$ would be simple. Any $(-1)$-curve in $F_0$ intersecting $\E$ is a tip of $F_0$, otherwise it would have $\mu>1$ and so it could not intersect $\wt D_h$, hence would be simple. We have $\sigma(F_0)\leq 3$, so since $F_0$ is connected, there is an $S_0$-component $M\subseteq F_0$ intersecting $\E$ and $D_0$ which is not exceptional (not a $(-1)$-curve). It follows that $\sigma(F_0)=3$, so $F_0$ is the only singular fiber.

Suppose $F_0$ is branched. Let $T$ be a maximal twig containing $L_1$ and let $R$ be the component of $\un{F_0}-T$ meeting $T$. Since $L_1,L_2$ are the only $(-1)$-curves of $F_0$, renaming $L_1$ and $L_2$ if necessary by a sequence of contractions of $(-1)$-curves different than $L_2$ we can contract the whole $T$. We have $\mu(R)>1$, otherwise this contraction would make $R$ into a non-tip component of a fiber with a unique $(-1)$-curve, which is impossible for $\mu(R)=1$ (cf. \cite[2.9(v)]{Palka-classification}). It follows that all components of $T$ have multiplicity bigger than $1$, so $\wt D_h\cdot T=0$. But $\wt D$ is connected, so this gives $\wt D\cdot L_1\leq 1$, a contradiction with (ii).

Since $F_0$ is a chain, $M$ is not branching, so (ii) implies that it intersects $\wt D_h$, hence $\wt D_h\cdot (L_1+L_2+D_0)\leq 2$. Since $\wt D_h\cdot D_0>0$, this gives, say, $\wt D_h\cdot L_1=0$. As $L_1$ is not simple, $L_1$ intersects two different connected components of $D_0$, which gives $\wt D_h\cdot D_0=2$ and $\wt D_h\cdot L_2=0$. Thus $L_2$ is simple, a contradiction. \end{proof}

The unique singular point of $S'$ is analytically of type $\C^2/G$ for some $G<GL(2,\C)$. We can and will assume that $G$ is small, i.e. it does not contain pseudo-reflections. Then $G$ is isomorphic to the local fundamental group of the singular point (cf. \cite{Brieskorn}, \cite[1.5.3.5]{Miyan-OpenSurf}). The divisor $\E$ is an admissible chain if $G$ is cyclic and an admissible fork otherwise. The discriminant is given by $d(\E)=|G/[G,G]|$ (see \cite{Mumford}). From (v) we see that the maximal twigs of $D$ are admissible, so since $d(D)<0$ by (iv), $D$ is not a chain. Moreover, (v) implies that $(\ov S,D+\E)$ is the unique snc-minimal completion of $S_0$ (cf. \cite[6.11(1)]{Palka-classification}). Let $T_i$ for $i=1,\ldots,s$ be all the maximal twigs of $D$, put $T=T_1+\ldots+T_s$. We put $$d_i=d(T_i), \delta_i=\delta(T_i), e_i=e(T_i), \wt e_i=e(T_i^t)$$ and $$\delta=\delta(D), e=e(D), \wt e=\wt e(D).$$ We write $\cal P$ for $(K+D+\E)^+$ and $\cal N$ for $(K+D+\E)^-$.

\blem\label{1lem:eps properties} The integer $\epsilon$ defined by the equality $(K+D+\E)^2=-1-\epsilon$ depends only on the isomorphism type of $S'$ and has the following properties (cf. \cite[5.3]{KR-ContrSurf}): \benum[(i)]

\item $\epsilon\geq 0$,

\item $K\cdot (K+D)=3-\epsilon-K\cdot E\leq 0$,

\item $\#\E+\#D=7+\epsilon+K\cdot D+K\cdot E$,

\item $\delta\leq e=-\Bk^2 D\leq 1+\epsilon+\Bk^2\E+\frac{3}{|G|}$. \eenum\elem

\begin{proof} Since the snc-minimal completion of $S_0$ is unique, $\epsilon$ is determined by the isomorphism type of $S'$. (i) Since $\cal N\neq 0$, by \ref{1lem:KobIneq}(i) we get $-1-\epsilon=(K+D+\E)^2<3\chi(S_0)=3(\chi(S')-1)=0$. (iii) Since $D$ and $\E$ are connected rational trees, their arithmetic genera vanish and we get  $K\cdot(K+D+\E)=3-\epsilon$, so $K^2=3-\epsilon-K\cdot D-K\cdot E$ and the formula follows from the Noether formula $K^2+\chi(\ov S)=12$. (ii) Suppose $K\cdot E+\epsilon\leq 2$. By the Riemann-Roch theorem $$h^0(-K-D)+h^0(2K+D)\geq K\cdot(K+D)+p_a(D)=3-\epsilon-K\cdot E>0,$$ so $-K-D\geq 0$, otherwise we would have $\kp(K+D)\geq 0$. We have $K\cdot\E>0$ and $K\cdot E_i\geq 0$ for every component $E_i$ of $\E$, hence $\E$ is in the fixed part of $-K-D$, so $-K-D-\E\geq 0$, which contradicts $\kp(K+D+\E)=2$. (iv) We have $\Bk^2D=-e$ by \ref{1lem:ChainBk}(iv) and $\cal N=\Bk D+\Bk \E$ by \ref{1prop:basic S' properties}(ii), so  $$-1-\epsilon=(K+D+\E)^2=\cal P^2+\Bk^2 D+\Bk^2 \E$$ and then (iv) is a consequence of \ref{1lem:KobIneq}(ii) applied to $(\ov S,D+\E)$. \end{proof}

\blem\label{1lem:if eps<2} Suppose $\epsilon<2$. Then:\benum[(i)]

\item $|2K+D+E|\neq\emptyset$,

\item $s-2-\frac{6}{|G|}\leq \delta$,

\item $s-3\leq \epsilon+\Bk^2 \E+\frac{9}{|G|}$, and if the equality holds then all twigs of $D$ are tips,

\item if $\Delta=\emptyset$ then $e+\delta\geq s+\epsilon+\frac{K\cdot E}{4}-\frac{5}{2}$. \eenum\elem

\begin{proof} (i) Riemann-Roch's theorem gives $h^0(-K-D-E)+h^0(2K+D+E)\geq 2-\epsilon$. If $-K-D-E\geq 0$ then $-K-D-\E\geq 0$, which contradicts $\kp(K+D+\E)=2$. Thus $2K+D+E\geq 0$. (ii) Let $R=D-T$. Each component of $\E+T$ is in the support of $\cal N$, hence intersects trivially with $\cal P$. By (i) and \ref{1lem:KobIneq}(ii) we have $$0\leq \cal P\cdot(2K+D+\E)=2\cal P\cdot(K+D+\E)-\cal P\cdot(D+\E)=2\cal P^2-\cal P\cdot R\leq \frac{6}{|G|}-\cal P\cdot R.$$ As $R$ is a rational tree, its arithmetic genus vanishes, so $$\cal P\cdot R=(K+D-\Bk D)\cdot R=-2+(T-\Bk D)\cdot R=-2+s-\delta$$ by \ref{1lem:ChainBk}(ii). (iii) is a consequence of \ref{1lem:eps properties}(iv), (ii) and the fact that if the inequality can become equality only if $e=\delta$.

(iv) Let $m$ be the biggest natural number, for which $|E+m(K+D)|\neq \emptyset$, $m\geq 2$ by (i). Write $E+m(K+D)\sim\sum a_iC_i$, where $a_i$ are positive integers and $C_i$ are distinct irreducible curves. We have $|K+D+\sum a_iC_i|=\emptyset$, so by \ref{1lem:linear systems}(i) $C_i$ are smooth rational curves, such that $C_i\cdot D\leq 1$. By \ref{1lem:linear systems}(ii) we can assume that they have negative self-intersections. Since $E+m(K+D)$ is effective, $E+m(K+D^\#)$ is effective by \ref{1lem:Eff-NegDef is Eff}, so we can write it as $$E+m(K+D^\#)\equiv \sum c_iC_i,$$ where $c_i>0$ and $C_i$ are as above. Note that $K\cdot E\geq 2$, otherwise $E=\E=[3]$ and $E\cdot (2K+D+E)=-1<0$, which would lead to $\ovk(K+D)\geq 0$ by (i). Suppose $(E+2K)\cdot C_i<0$ for some $i$, say $i=1$. If $C_1\nsubseteq E$ then, since $C_1\cdot D\leq 1$ and since $\Delta=\emptyset$, we have $C_1\cdot E\geq 2$ by \ref{1prop:basic S' properties}(ii), so $K\cdot C_1<-\frac{1}{2}C_1\cdot E\leq -1$, which contradicts $C_1^2<0$, as $C_1\cong \PP^1$. Thus $C_1\subseteq E$. But then $K\cdot C_1\geq 0$ and $$0>(E+2K)\cdot C_1=K\cdot C_1+\beta_E(C_1)-2,$$ so since $\Delta=\emptyset$, we get $E=C_1$ and $K\cdot E\leq 1$, a contradiction. We infer that $0\leq (E+2K)\cdot (E+m(K+D^\#))$. We have $$(E+2K)\cdot(K+D)=2K\cdot(K+D+E)-K\cdot E=6-2\epsilon-K\cdot E$$ and $$\Bk D\cdot K=\Bk D\cdot(K+D^\#)+\Bk^2D-\Bk D\cdot (D-T)-\Bk D\cdot T=0-e-\delta+s,$$ so from the above inequality we get $$s-\delta-e\leq \frac{1}{2m}(K\cdot E-2)+3-\epsilon-\frac{1}{2}K\cdot E\leq \frac{1}{4}(K\cdot E-2)+3-\epsilon-\frac{1}{2}K\cdot E,$$ which gives (iv). \end{proof}

\med\section{Bounding the shape of the exceptional divisor}\label{sec:bounding exc divisor}

\bprop\label{1prop:Zhp has a smooth covering} Let $X$ be $\Z$-homology plane with a unique singular point, which is of analytical type $\C^2/\Z_a$. Then there exists a smooth affine surface $Y$ with an action of $\Z_a$ on it, which has a unique fixed point, is free on its complement and for which $X\cong Y/\Z_a$. \eprop

\begin{proof} We modify a bit the arguments of \cite[2.2]{Koras-A2}. Let $q\in X$ be the singular point. Then there is a (contractible) neighborhood $N\subseteq X$ of $q$, which is analytically isomorphic to $\C^2/\Z_a$. Let $p\:(\C^2,0)\to (N,q)$ be the quotient map and let $j$ be the embedding of $N-q$ into $X-q$. Let $G$ be the commutator of $\pi_1(X-q)$ and let $Y_0\to X-q$ be the covering corresponding to the inclusion $G\embed \pi_1(X-q)$. We show that $Y=Y_0\cup \{0\}$ is smooth. Since $\C^2-0$ is simply connected, $p_{|\C^2-0}$ has a lifting $\wt p:\C^2-0\to Y_0$. The embedding $(N,N-q)\embed (X,X-q)$ induces a morphism of long homology exact sequences of respective pairs. The reduced homology groups of $N$ and $X$ vanish, so in both sequences the boundary homomorphisms are isomorphisms. By the excision theorem $H_2(N,N-q,\Z)\cong H_2(X,X-q,\Z)$, hence $H_1(N-q,\Z)\to H_1(X-q,\Z)$ is an isomorphism. Since $\pi_1(N-q)$ is abelian, it follows that the composition $$\pi_1(N-q)\to \pi_1(X-q)\to H_1(X-q,\Z)$$ is an isomorphism. Let $y_1,y_2\in \C^2-0$ be two points lying over the same point in $N-q$, such that $\wt p(y_1)=\wt p(y_2)$. The path joining $y_1$ and $y_2$ in $\C^2-0$ maps by $\wt p$ to a loop $Y_0$. Let $\alpha\in \pi_1(N-q)$ be a loop which is an image in $N-q$ of the same path. Then $\pi_1(j)(\alpha)\in \pi_1(X-q)$ belongs to $G$, hence $\alpha$ is in the kernel of the composition $$\pi_1(N-q)\to \pi_1(X-q)\to H_1(X-q,\Z),$$ which is trivial. We get that $y_1=y_2$, so $\wt p$ is a monomorphism and we see that the local fundamental group of $Y$ at $0$ is trivial. By \cite{Mumford} (the proof is topological and works for non-algebraic surfaces) we see that $Y$ is smooth.

As a finite unbranched cover of an algebraic variety $Y_0$ is algebraic and the map $Y_0\to X-q$ is finite, so $\C[Y_0]$ is an integral extension of $\C[X-q]\cong \C[X]$, hence it is a finitely generated and integrally closed $\C$-algebra. The homomorphism $\C[X]\to \C[Y_0]$ induces a morphism $r\:\Spec \C[Y_0]\to X$. The natural embedding $\psi\:Y_0\to \Spec \C[Y_0]$ is an isomorphism onto $r^{-1}(X-q)$ and extends to a morphism by smoothness of $Y$. The inverse extends to a morphism from $\Spec \C[Y_0]$ to $X$ by normality of $\Spec \C[Y_0]$. \end{proof}

The following theorem is a key step in the proof of the main result of the paper. It is based on the method of finding well-behaved exceptional curves on open surfaces of negative Kodaira dimension introduced in \cite[4.2, 4.3]{KR-C*onC3} and which has its origin in \ref{1lem:linear systems}(iii).

\bprop\label{1prop:KE and eps} Either $K\cdot E+2\epsilon\leq 5$ or $\epsilon=2$ and then $\E=[4]$ and $D$ consists of $(-2)$-curves.\eprop

\begin{proof} Note that $(2K+E)\cdot(K+D)=6-2\epsilon-K\cdot E$, so $K\cdot E+2\epsilon\leq 5$ is equivalent to $(2K+E)\cdot(K+D)>0$. Under two additional assumptions, that there exists a $(-1)$-curve $A\subseteq \ov S$, such that $A\cdot\E\leq 1$ and that $S'$ is contractible, it is proved in  \cite[5.10, 5.11]{KR-ContrSurf} that the inequality $(2K+E)\cdot(K+D)\leq 0$ implies that there exists an exceptional simple curve on $(\ov S,D+\Delta)$, which intersects $\Delta$. Of course, it also intersects $D$, as $S'$ is affine. Moreover, under the above assumptions the process of contracting and finding such $(-1)$-curves can be iterated to infinity. By the definition of simplicity this is impossible, because the number of connected components of $\Delta$ is finite. The proof of 5.10 loc. cit. does not require the contractibility, but only the $\Q$-acyclicity of $S'$, so it can be simply repeated in our situation. However, the existence of the 'initial' curve $A$, which is assured by lemma 5.7 loc. cit. in case $S'$ is contractible, has to be reconsidered in our situation.

Suppose $K\cdot E+2\epsilon> 5$. From the above remarks  it follows that we can assume that there is no $(-1)$-curve $A\subseteq \ov S$ with $A\cdot \E\leq 1$. We can repeat the proof by contradiction in 5.7 loc. cit. up to 5.7.4(i). In 5.7.4(ii) an argument referring to \cite{Koras-A2} (and hence to contractibility) is used and it needs to be modified in our situation. We are therefore in a situation where $K+\E^\#\equiv 0$, $\Bk^2\E$ is an integer and $D$ consists of $(-2)$-curves. As $\E$ does not consist of $(-2)$-curves, by \ref{1rem:-2 chains and forks} and \ref{1lem:ForkBk}(v) $\Bk^2\E=-1$ and $\E$ is a chain. We have now $$-1-\epsilon=(K+D+\E)^2=(D+\Bk \E)^2=D^2-1,$$ hence $\epsilon=-D^2=2+K\cdot D=2$. By Riemann-Roch's theorem $$h^0(\E+2K)+h^0(-K-\E)\geq K\cdot(K+\E)=3-\epsilon-K\cdot D=1.$$ If $-K-\E\sim U$ for an effective divisor $U$ then $K+\E^\#\equiv 0$ implies $U+\Bk\E\equiv 0$, hence $\Bk \E=0$, which is impossible by \ref{1lem:ChainBk}(iii). Recall that for a $\Q$-divisor $T$ we denote the integral and fractional parts of $T$ by $[T]$ and $\{T\}$ respectively. We get $2(K+\E)\geq 0$, which by \ref{1lem:Eff-NegDef is Eff}(ii) implies that $[2(K+\E^\#)]\sim U$ for some effective divisor $U$. Then $$0\equiv 2(K+\E^\#)\equiv [2(K+\E^\#)]+\{2(K+\E^\#)\}\equiv U+\{-2\Bk\E\},$$ so since $\{-2\Bk\E\}$ is effective, $\{-2\Bk\E\}=U=0$. Thus $2\Bk\E$ is a $\Z$-divisor. Since $\E$ is not a $(-2)$-chain, $\E\neq \Bk \E$ and we get $2\Bk \E=\E$ and $2K+\E=2K+2\E^\#\sim U=0$. It follows that $\Delta=0$ and $K\cdot E=2$. Moreover, as $E_i\cdot(2K+\E)=0$ for each component $E_i$ of $\E$, we get that either $\E=[4]$ or $\E=[3,(k),3]$ for some $k\geq 0$ (recall that $[(k)]$ is a chain of $(-2)$-curves of length $k$). To finish the proof we need to exclude cases other than $\E=[4]$.

Suppose $\E=[3,(k),3]$ for some $k\geq 0$. We have $\#D=9-k$ by \ref{1lem:eps properties}(iii), so there are only finitely many possibilities for the weighted dual graph of $D$. The inequality \ref{1lem:eps properties}(iv) gives $$e(D)\leq 3+\Bk ^2\E+\frac{3}{|G|}=2+\frac{3}{d(E)}=2+\frac{3}{4(k+2)}.$$ $D$ consists of $(-2)$-curves, so $e(D)=s-\delta$. Taking a square of \ref{1prop:basic S' properties}(ii) we get $-3=\cal P^2-e(D)-1$, so $\cal P^2=s-2-\delta$. Since $\cal P^2>0$, we obtain: $$0<s-2-\delta \leq\frac{3}{4(k+2)} =\frac{3}{4(11-\#D)}.$$ In particular, $s-2\leq\delta+\frac{3}{8}\leq \frac{s}{2}+\frac{3}{8}$, so $s\leq 4$. Another condition is given by \ref{1prop:basic S' properties}(iv): $$\sqrt{-\frac{d(D)}{d(E)}}\in\N.$$ We check by a direct computation that there are only two pairs of weighted dual graphs of $D$ and $\E$ satisfying both conditions (one checks first that the first condition implies that $k\leq 1$ for $s=3$ and $k\leq 2$ for $s=4$):\benum[(1)]

\item $s=3$, $T_1=[2,2]$, $T_2=[2,2,2]$, $T_3=[2,2,2]$, $\E=[3,3]$,

\item $s=4$, $T_1=[2]$, $T_2=[2]$, $T_3=[2]$, $T_4=[2,2,2]$, $\E=[3,3]$. \eenum

\noin Note that in case (2) $D-T_1-T_2-T_3-T_4$ has three components. In both cases $-d(D)=d(\E)=8$, so $H_1(S',\Z)=0$ by \ref{1prop:basic S' properties}(iv). By \ref{1prop:Zhp has a smooth covering} $S'$ can be identified with the image of a quotient morphism $p\:Y\to Y/\Z_8$ of some smooth affine surface $Y$. Let $(x,y)$ be local parameters which are semi-invariant with respect to the action of $\Z_8$ (recall that $t\in \C(Y)$ is semi-invariant with respect to the action of $G$ on $Y$ if there exists a character $\chi\:G\to \C^*$, such that $g^*t=\chi(g)t$). As in the case of $\C^2\to \C^2/\Z_8$, if $C$ is the proper transform on $S$ of $p(\{x=0\})$ then $C\cdot \E=1$ and $C$ meets $\E$ is a tip (cf. \cite{Hirzebruch-cyclic_singularities}). Thus $K\cdot C=-\frac{1}{2}\E\cdot C=-\frac{1}{2}$, a contradiction. \end{proof}

\bcor\label{1cor:possible KE and eps} If $\epsilon=0$ then $K\cdot E\in\{3,4,5\}$. If $\epsilon=1$ then $K\cdot E\in\{2,3\}$. If $\epsilon=2$ then either $K\cdot E=1$ or $\E=[4]$. \ecor

\begin{proof} We have $K\cdot E+\epsilon\geq 3$ and $\epsilon\geq 0$ by \ref{1lem:eps properties}(i),(ii). By \ref{1prop:KE and eps} we have $K\cdot E+2\epsilon\leq 5$ for $(\E,\epsilon)\neq ([4],2)$, so the corollary follows. \end{proof}

\bprop\label{1prop:if eps=0 or E^ fork}\ \benum[(i)]

\item If $\epsilon=0$ then $\#\E=1$ and $D$ is a fork,

\item If $\E$ is a fork then $\epsilon=2$,

\item $\Delta$ does not contain a fork. \eenum\eprop

\begin{proof}
\ssk (i) Since $D$ is not a chain we have $s\geq 3$. For $\epsilon=0$ lemma \ref{1lem:if eps<2}(iii) gives $0\leq s-3\leq \Bk^2 \E+\frac{9}{|G|}$. If $\E$ is a fork then $\Bk^2 \E<-1$ by \ref{1lem:ForkBk}(v), so $|G|\leq 8$. Since $G$ is small and non-abelian, it is the quaternion group, for which the resolution consist of $(-2)$-curves (the abelianization of the group is $\Z_2\times\Z_2$, row $2$ is the table \cite[Satz 2.11]{Brieskorn}), a contradiction with \ref{1prop:basic S' properties}(iii). Thus $\E$ is a chain, so $d(\E)=|G|$ and we get $d'(\E)+d'(\E^t)\leq7$ by \ref{1lem:ChainBk}(iv). Suppose $\E$ has more than one component. Taking into account \ref{1cor:possible KE and eps} there are two possibilities for $\E$: $[3,4]$ and $[2,5]$. In both cases we obtain $\Bk^2\E+\frac{9}{|G|}=0$, so $s=3$ and the inequalities \ref{1lem:eps properties}(iv) and \ref{1lem:if eps<2}(ii) become equalities. We get $e=\delta<1$, which is possible only if maximal twigs of $D$ are tips. Denoting the branching component of $D$ by $B$ we have $d(D)=d_1 d_2 d_3 (-B^2-\delta)$, so since $d(D)<0$, we get $-B^2<\delta<1$, a contradiction with \ref{1prop:basic S' properties}(v). Therefore $\#\E=1$. If $s\neq3$ then \ref{1lem:if eps<2}(iii) and \ref{1cor:possible KE and eps} give subsequently $(s-3)d(\E)\leq5$, $s=4$ and $\E=[5]$. Then $e=\delta=\frac{4}{5}$, so the inequality \ref{1lem:if eps<2}(iv) fails, a contradiction.

\ssk (ii) Let $\E$ be a fork. By (i) $\epsilon\neq 0$. Suppose $\epsilon=1$. Then $\Bk^2\E+\frac{9}{|G|}+1\geq 0$, so since $\Bk^2\E<-e(\E)$, we get $|G|(e(\E)-1)\leq 9$. One checks using \cite[Satz 2.11]{Brieskorn} that the last inequality is satisfied only for the fork $\E$ which has $[2], [2], [3]$ as maximal twigs and $[2]$ as a branching curve. In this case $\Bk^2\E=-\frac{3}{2}$ and $|G|=24$, so the initial inequality fails.

\ssk (iii) Suppose $\Delta$ contains a fork. Then $\epsilon=2$ by (ii), so $\#E=1$ by \ref{1cor:possible KE and eps}. By \ref{1lem:Eff-NegDef is Eff} we have $$\ovk(S\setminus\Delta)=\kappa(K_{\ov S}+D+\Delta)=\kappa(K_{\ov S}+D)=\ovk(S)=-\8.$$ Suppose $S\setminus\Delta$ is affine-ruled. Consider a minimal completion $(\wt S,\wt D+\Delta)\to B$ of this ruling (cf. \ref{1def:ruling completion}). Since $S'$ is affine, the horizontal component is contained in $\wt D$. If $E$ is vertical then $S_0$ is affine-ruled, which contradicts $\ovk(S_0)=2$. Thus there are two horizontal components in $\wt D+E$. Since $E\cap \wt D=\emptyset$, we have $\nu=0$, so $\Sigma_{S_0}=0$ by \ref{1prop:basic S' properties}(vii), hence each singular fiber has a unique $(-1)$-curve. Then each connected component of $\Delta$ is a chain, a contradiction. By \cite{MiTs-PlatFibr} $S\setminus\Delta$ contains an open subset $U$, which is Platonically $\C^*$-fibred. In particular $S\setminus \Delta$ is $\C^*$-ruled (we have shown that it is not affine-ruled). The component $E$ cannot be vertical for this ruling, otherwise $S_0$ is $\C^*$-ruled, which contradicts $\ovk(S_0)=2$. Consider a minimal completion of this ruling. We have $\nu=0$, so $\Sigma_{S_0}=1$. By the description of the Platonic fibration in loc.cit. the branching component of the fork contained in $\Delta$ is horizontal. Let $F_0$ be the fiber containing two $S_0$-components, call them $L_1$ and $L_2$. By minimality only these curves can be $(-1)$-curves of $F_0$. Decompose $\Delta$ into $\Delta_1+\Delta_2$, where $\Delta_1$ is a fork and $\Delta_2$ is a chain (may be empty). Since $\wt D\cap F_0$ is connected and since $S'$ is affine, we have $L_1\cdot \wt D=L_2\cdot \wt D=1$. This gives $(L_1+L_2)\cdot \Delta_1=1$ because $F_0$ and $\Delta_1$ are trees. Say $L_1\cdot \Delta_1=1$ and $L_2\cdot \Delta_1=0$. If only one of the $L_i$'s is a $(-1)$-curve then it follows from the structure of a singular fiber with a unique $(-1)$-curve that it has to be $L_2$, as $\Delta_1$ intersects a component of $F_0$ of multiplicity one. In any case we get that $L_2^2=-1$, $L_2\cdot \Delta_1=0$ and by the negative definiteness of the intersection matrix of the fiber $L_2+\Delta_2$ is a chain. Analyzing the contraction of this chain as in \cite[6.1]{KR-ContrSurf} one shows that the fact that $K\cdot E=1$ leads to $L_2\cdot \E=1$, i.e. $L_2$ is simple on $(\wt S,\wt D+\E)$, hence on $(\ov S,D+\E)$, which contradicts \ref{1prop:basic S' properties}(ii). \end{proof}

\bcor\label{1cor:S-Delta affine-ruled} $S\setminus\Delta$ is affine-ruled. \ecor

\begin{proof} The logarithmic Kodaira dimension of  $S\setminus\Delta$ is negative, so by the structure theorems mentioned above $S\setminus\Delta$ is affine-ruled or it contains a Platonic fibration as an open subset. The last case is possible only if $\Delta$ contains a fork, which is excluded by \ref{1prop:if eps=0 or E^ fork}(iii). \end{proof}

Recall that $[(k)]$ denotes a chain of $(-2)$-curves of length $k$ and that the default orientation of a twig is the one in which the first component is a tip of the divisor and the last component intersects some component of the divisor not contained in the twig.

\bprop\label{1prop:possible $E$} $\E$ is of one of the following types:\benum

\item[(a)] $[5]$, $[6]$, $[7]$

\item[(b1)] fork: $$ \xymatrix{{A}\ar@{-}[r] &{-2}\ar@{-}[r]\ar@{-}[d]& {B}\\
{} & {-2} & {} } $$ with $(A,B)$ equal to: $([3],[2,2])$, $([3],[2,2,2])$, $([3],[2,2,2,2])$, $([2,3],[2,2])$ or $([(n),3],[2])$, where $n\geq 0$,

\item[(b2)] fork: $$ \xymatrix{{A}\ar@{-}[r] &{-3}\ar@{-}[r]\ar@{-}[d]& {B}\\
{} & {-2} & {} } $$ with $(A,B)$ equal to one of: $([2,2],[2,2])$, $([2,2],[2,2,2])$, $([2,2],[2,2,2,2])$ or $([2],[(n)])$, where $n\geq 0$,

\item[(b3)] $[(r),3,(x)]$ for $r,x\geq0$,

\item[(c1)] $[(r),4]$ or $[(r),5]$ for $r\geq 0$,

\item[(c2)] $[(x),3,(y),3]$ or $[(x),3,(y),4]$ or $[(x),4,(y),3]$ for $x,y\geq 0$,

\item[(c3)] $[(r),3,(x),3,(y),3]$ for $r,x,y\geq 0$,

\item[(c4)] $[2,4,2]$, $[2,5,2]$, $[2,3,3,2]$, $[2,3,4,2]$, $[2,4,2,2]$, $[2,5,2,2]$.\eenum\eprop

\begin{proof} If $\E$ is a fork then $\epsilon=2$ by \ref{1prop:if eps=0 or E^ fork}(ii), so $E=[3]$ by \ref{1cor:possible KE and eps}. We know that $\Delta$ does not contain a fork, so all possible $\E$'s satisfying \ref{1lem:ForkBk}(ii)-(iii) are listed in (b1) and (b2). Chains for $\epsilon=2$ other than $[4]$ are in (b3) and $\E$'s for $\epsilon=0$ are in (a) (cf. \ref{1cor:possible KE and eps} and \ref{1prop:if eps=0 or E^ fork}(i)). Now we can assume that $\E$ is a chain and $\epsilon=1$, so $K\cdot E\in\{2,3\}$ by \ref{1cor:possible KE and eps}. The possibilities with $E\cdot \Delta\leq 1$ are listed in (c1), (c2) and (c3), so we can now assume $E\cdot \Delta=2$. If $T$ is an oriented chain with the first component $T_1$ then we write $d''(T)$ for $d'(T-T_1)$. From \ref{1lem:if eps<2}(iii) we get $d'(\E)+ d'(\E^t)\leq d(\E)+7$ and since $$d(\E)=2d'(\E)-d''(\E)=2d'(\E^t)-d''(\E^t),$$ we have $$\frac{1}{2}(d(\E)+d''(\E))+\frac{1}{2}(d(\E)+d''(\E^t))\leq d(\E)+7,$$ so $d''(\E)+d''(\E^t)\leq 14$. This gives six possibilities for $\E$: $[2,4,2]$, $[2,5,2]$, $[2,3,3,2]$, $[2,3,4,2]$, $[2,4,2,2]$ and $[2,5,2,2]$, which are listed in (c4). \end{proof}

\med\section{Special affine rulings of the resolution}\label{sec:pre-minimal rulings}

In this section we assume that $\#E=1$, i.e. the exceptional divisor of the snc-minimal resolution $S\to S'$ has a unique component with self-intersection different than $(-2)$ (in terms of the list \ref{1prop:possible $E$} this holds in cases (a),(b),(c1) and part of (c4)). Under this assumption we will produce and analyze special affine rulings of $S\setminus\Delta$ (hence of $S$).

We keep the notation $(\ov S,D)$ for the unique snc-minimal smooth completion of $S$. Consider an affine ruling of $S\setminus\Delta$ (it exists by \ref{1cor:S-Delta affine-ruled}). There exists a modification $(\ov S^\dagger,D^\dagger)\to (\ov S,D)$ and a $\PP^1$-ruling $f\:(\ov S^\dagger,D^\dagger+\Delta)\to\PP^1$, which is a minimal completion of the affine ruling. Clearly, $E$ is horizonal, otherwise $S_0$ is affine-ruled, which contradicts $\ovk(S_0)=2$. It follows that $\nu=0$ and since $\#E=1$, we have $h=2$ and hence $\Sigma_{S_0}=0$ by \ref{1prop:basic S' properties}(vii). Thus every fiber of $f$ contains a unique $S_0$-component and since $f$ is minimal, it is the unique $(-1)$-curve of the fiber in case the fiber is singular. As we have seen in \ref{1def:char pairs for fibers}, once we fix a component of $F$ of multiplicity one, $F$ can be uniquely described by a sequence of characteristic pairs recovering $F$ from (the birational transform of) the component. In our situation the default choice is the component of $F$ intersecting the horizontal component of $D^\dagger$.

\begin{notation}\label{2no:pairs notation} Let $f$ be a completion of an affine ruling of $S\setminus\Delta$ as above. Let $F$ be some fiber of $f$ and let $H$ be the section contained in $D^\dagger$. Put $\gamma=-E^2$, $n=-H^2$ and $d=E\cdot F$. Let $h$ be the number of characteristic pairs of $F$. We write $\Delta\cap F=\Delta_1+\ldots+\Delta_k$, $k\geq 0$ where $\Delta_i$ are irreducible and $\Delta_k$ is a tip of $F$. If the fiber is singular then it follows that the last pair of $F$ is $\binom{c_h}{p_h}=\binom{k+1}{1}$. If $\Delta\neq \emptyset$ then $E\cdot \Delta_{i_0}=1$ for a unique $1\leq i_0\leq k$, because $\E$ is a tree. In case $\Delta\cap F=\emptyset$ put $i_0=0$. Define $F'$ as the image of $F$ after contraction of curves produced by $\binom{c_h}{p_h}$ and let the sequence of characteristic pairs for $F'$ be $\binom{\un c_i}{\un p_i}$ with $i=1,\ldots,h-1$ (if $h=1$ then $\binom{\un c_1}{\un p_1}=\binom{1}{0}$). Put $c_h'=c_h-i_0$ and $\mu=\mu_F(C)$, where $C$ is the unique $(-1)$-curve of $F$. We define $$\kappa=c_h C\cdot E+c_h'\text{\ and \ } \rho=\kappa C\cdot E+c_{h}'C\cdot E+c_{h}'.$$ If $f$ has exactly two singular fibers, we write the analogous quantities for the second fiber with $\wt{(\ )}$: $\wt \kappa$, $\wt C$, $\wt{\un p}_i$, $\wt c_h'$ etc. If $f$ has more singular fibers then instead of $\kappa$, $C$, $\un p_i$, $c_h'$, etc. we write $\kappa(F)$, $C_F$, $\un p_i(F)$, $c_h'(F)$, etc. \end{notation}

It follows from the definition that $\un c_i=c_i/c_h$ and $\un p_i=p_i/c_h$, so  $gcd(\un c_i,\un p_i)=\un c_{i+1}$ for $i=1,\ldots,h-1$ and $gcd(\un c_{h-1},\un p_{h-1})=1$ if $h>1$. The multiplicities of $C$ and $\Delta_{i_0}$ in $F$ are $\mu=\un c_1 c_h$ and $\un c_1 c_h'$, so $$d=E\cdot F=c_1E\cdot C+\un c_1 c_h'E\cdot \Delta_{i_0}=\un c_1 \kappa.$$ Note that $c_h'=0$ if and only if $\Delta\cap F=\emptyset$ if and only if $c_h=1$.

We denote the least common multiple of a set $M$ of natural numbers by $lcm(M)$.

\bprop\label{2prop:pairs} With the notation as in \ref{2no:pairs notation} the following equations hold (cf. \cite[6.10, 6.11]{KR-ContrSurf}):

\begin{align}d(n+2)+\gamma-2&=\sum_F \kappa(F) (\un c_1(F)+\sum_{i=1}^{h(F)-1}\un p_i(F)) , \label{eq:pairs-general I}\\ nd^2+\gamma&=\sum_F(\kappa^2(F)\sum_{i=1}^{h(F)-1}\un c_i(F)\un p_i(F)+\rho(F)) \label{eq:pairs-general II}, \\ d\cdot |H_1(S',\Z)|&=\prod_F\un c_1(F), \label{eq:pairs-general Pic}\\ d&=\underset{F}{lcm}\{\un c_1(F)\},\label{eq:pairs-general-d-is-lcm}\end{align} where $F$ runs over all singular fibers of $f$. \eprop

\begin{proof} First we derive the equations (1) and (2). For simplicity we assume that there is a unique singular fiber, the general case follows. We have $\Sigma_{S_0}=0$. Consider the sequence of blow-downs $$\ov S=S^{(m)}\xra{\sigma_m} S^{(m-1)}\xra{\sigma_{m-1}}\ldots\xra{\sigma_1} S^{(0)},$$ $S^{(0)}$ a Hirzebruch surface, which contracts $F$ to a smooth $0$-curve without touching $H$. Denote by $K^{(j)}$ and $E^{(j)}$ the canonical divisor of $S^{(j)}$ and the birational transform of $E$ on $S^{(j)}$ respectively. Denoting the multiplicity of the center of $\sigma_{j}$ on $E^{(j-1)}$ by $\mu_{j}$ we have $$K^{(j)}\cdot E^{(j)}-K^{(j-1)}\cdot E^{(j-1)}=\mu_j\text{\ \ and\ \ } (E^{(j-1)})^2-(E^{(j)})^2=\mu_j^2,$$ $j=1,\ldots,m$. We have $E^{(0)}\equiv d(nF^{(0)}+H)$, where $F^{(0)}$ is some fiber of the induced $\PP^1$-ruling of $S^{(0)}$ and $d=E^{(0)}\cdot F^{(0)}=E\cdot F$. We compute $$K^{(m)}\cdot E^{(m)}-K^{(0)}\cdot E^{(0)}= K\cdot E+d(n+2)= \gamma-2+d(n+2)$$ and $$(E^{(0)})^2-(E^{(m)})^2=nd^2+\gamma,$$ which gives left sides of the above equations. We thus need to compute $\sum \mu_j$ and $\sum \mu_j^2$. Let $F',\un c_i,\un p_i,\kappa$ be as defined above. Let us first consider the case $\Delta\cap F=\emptyset $. We then have $\kappa=C\cdot E$ and the sequence of characteristic pairs for $F$ is $\binom{\un c_1}{\un p_1},\ldots,\binom{\un c_{h-1}}{\un p_{h-1}},\binom{1}{1}$. The sequence of blow-downs $\sigma_j$ is divided into groups described by these pairs. The set of indices $j$, for which the blow-up $\sigma_j$ is a part of the group of blow-downs determined by the characteristic pair $\binom{c_i}{p_i}$ will be denoted by $I_i$. In case $\kappa=C\cdot E=1$ we get by \ref{1lem:sum of mu and squares} $$\sum_{j\in I_i}\mu_j=c_i+p_i-gcd(c_i,p_i)\text{\ \ and\ }\sum_{j\in I_i}\mu_j^2=c_ip_i.$$ Now for $C\cdot E=\kappa\geq 1$ the multiplicity of each center is $\kappa$ times bigger, hence in general we get $$\sum_{j\in I_i}\mu_j=\kappa(c_i+p_i-gcd(c_i,p_i))\text{\ \ and\ }\sum_{j\in I_i}\mu_j^2=\kappa^2c_ip_i.$$ We have $c_h'=0$ and $c_h=1$, so this gives $$\sum \mu_j=\kappa\sum_{i=1}^h(\un c_i+\un p_i-gcd(\un c_i,\un p_i))= \kappa(\un c_1+\sum_{i=1}^{h}\un p_i-1)=\kappa (\un c_1+\sum_{i=1}^{h-1}\un p_i)$$ and $$\sum \mu_j^2=\kappa^2\sum_{i=1}^h \un c_i\un p_i=\kappa^2(\sum_{i=1}^{h-1}\un c_i\un p_i+1),$$ as required.

We now consider the case $\Delta\cap F\neq \emptyset$. Let $E'$ be the image of $E$ after contracting $F$ to $F'$. It follows from the above arguments that $$K'\cdot E'-K^{(0)}\cdot E^{(0)}=\kappa (\un c_1+\sum_{i=1}^{h-1}\un p_i-1)\text{\ and\ }(E^{(0)})^2-(E')^2=\kappa^2\sum_{i=1}^{h-1}\un c_i\un p_i,$$ so we need to compute $K\cdot E-K'\cdot E'$ and $E'^2-E^2$. We are now left with the last pair $\binom{c_h}{p_h}$, which groups $c_h=c_h'+i_0$ blow-ups. The proper transform of $E'$ after making first $c_h'$ blow-ups is $E^{(m-i_0)}$. The multiplicity of the center of each of these blow-ups is $C\cdot \E=C\cdot E+1$, so $$K^{(m-i_0)}\cdot E^{(m-i_0)}-K'\cdot E'=c_h'(C\cdot E+1)\text{\ and\ }E'^2-(E^{(m-i_0)})^2=c_h'(C\cdot E+1)^2.$$ Now $E^{(m-i_0)}$ may intersect the fiber in more than one point. The multiplicity of the center of each of the remaining $i_0$ blow-ups is $C\cdot E$, hence $$K\cdot E-K^{(m-i_0)}\cdot E^{(m-i_0)}=i_0 C\cdot E \text{\ and\ } (E^{(m-i_0)})^2-E^2=i_0 (C\cdot E)^2.$$ This gives \eqref{eq:pairs-general I} and \eqref{eq:pairs-general II}.

We now derive (3). Put $Q(F)=\sum_{i=1}^{h(F)-1}\un c_i(F)\un p_i(F)$ and $e(F)=d(F\cap \Delta-\Delta_{i_0(F)})/c_h(F)=c_h'(F)(c_h(F)-c_h'(F))/c_h(F)$. Then, as in \cite[3.4.6]{KR-C*onC3} $\rho(F)=\kappa(F)^2/c_h(F)+e(F)$, so we can rewrite \eqref{eq:pairs-general II} as: $$nd^2+\gamma-\sum_F e(F)=\sum_F\kappa^2(F)(Q(F)+1/c_h(F)),$$ which by 3.5.5 loc. cit. gives

\begin{align}nd^2+d(\E)/\prod_Fc_h(F)=\sum_F\kappa^2(F)(Q(F)+1/c_h(F)).\label{eq:pairs-general IIprim}\end{align}

$\Pic \ov S$ is a free abelian group with generators $f$ (general fiber), $H$ and vertical components not intersecting $H$. Let $G(F)$ be the component of $F$ intersecting $H$. Then $\Pic S_0$ is a generated by $f$ and $S_0$-components $C_F$ with defining relations coming from $E\sim 0$ and $G(F)\sim 0$ for any singular fiber $F$. The latter gives $f\sim \mu(C_F) C_F$. Expand $E$ in terms of the above generators, let $-k_F$ be the coefficient of $C_F$ and let $a,b$ be the coefficients of $f$ and $H$. Intersecting with $f$ and then with $H$ we get $b=d=E\cdot f$ and $a=bn=dn$, hence the relation coming from $E\sim 0$ is $\sum_F k_F C_F\sim dn f$. In the proof of 3.6 loc. cit. it is shown that $k_F=\kappa(F)(c_h(F)Q(F)+1)$, so taking the determinant of the defining relations we obtain $$\pm |\Pic S_0|/\prod_F\mu(C_F)=-nd+\sum_F\kappa(F)/\mu(C_F)(c_h(F)Q(F)+1).$$ Multiplying both sides by $d$ we have $$nd^2 \pm d|\Pic S_0|/\prod_F\mu(C_F)=\sum_Fd\kappa(F)c_h(F)/\mu(C_F)(Q(F)+1/c_h(F)).$$ Since $dc_h(F)/\mu(C_F)=\un c_1(F)\kappa(F)c_h(F)/(\un c_1(F)c_h(F))=\kappa(F)$, left sides of the above equation and of \eqref{eq:pairs-general IIprim} are the same, which gives $$d\cdot |\Pic S_0|=d(\E)\cdot \prod _F\un c_1(F).$$ Now (3) follows from by \ref{1prop:basic S' properties}(viii).

We have $\pi_1(S')=\pi_1(S)$ by \ref{1prop:basic S' properties}(iv). Note that the greatest common divisor of $S$-components of a fiber equals $\un c_1(F)$. Then by \cite[4.19, 5.9]{Fujita} $\pi_1(S)$ is generated by $\sigma_F$, where $F$ runs over singular fibers of $F$, and the defining relations are $(\sigma_F)^{\un c_1(F)}=1$ and $\prod \sigma_F=1$. Hence $H_1(S,\Z)$, which is the abelianization of $\pi_1(S)$, is the quotient of $\bigoplus_F \Z_{\un c_1(F)}$ by the subgroup generated by $(1,\ldots,1)$. We obtain $|H_1(S',\Z)|=(\prod_F \un c_1(F))/m$, where $m=lcm_F\{\un c_1(F)\}$, i.e. $m$ is the least common multiple of all $\un c_1(F)$'s. Plugging into (3) gives (4). \end{proof}

\bdfn\label{2def:pre-min ruling} Let $\pi:X\to C$ be a dominating morphism of a normal surface to a complete curve $C$. We say that $\pi$ is \emph{pre-minimal} if for some normal completion $(\ov X,\ov X\setminus X)$ it has an extension $\ov \pi:\ov X\to C$, such that the boundary divisor $\ov X\setminus X$ can be made snc-minimal using only subdivisional blow-downs. Then we will say also that $\ov \pi:(\ov X,\ov X\setminus X)\to C$ is pre-minimal.\edfn

\bcor\label{2cor:two fibers equations} Let $\#E=1$ and let $f$ be a minimal completion of an affine ruling of $S\setminus\Delta$. Then $f$ has at least two singular fibers and if it has two then in the notation of \ref{2no:pairs notation} one has: \benum[(i)]

\item $\un c_1=\wt \kappa\cdot |H_1(S',\Z)|$ and $\wt {\un c}_1=\kappa\cdot |H_1(S',\Z)|,$

\item $h,\wt h\geq 2$,

\item $d(D)=-d(\E)\cdot gcd(\un c_1,\wt {\un c}_1)^2$.

\item if $f$ is pre-minimal then $h+\wt h=n+1+\epsilon+K\cdot E$.
\eenum\ecor

\begin{proof}  Note that by \ref{1prop:basic S' properties}(ii) $\kappa(F)\geq 2$ for every fiber $F$. If $f$ has only one singular fiber then \eqref{eq:pairs-general Pic} gives $\un c_1=d\cdot |H_1(S',\Z)|= \un c_1 \kappa \cdot |H_1(S',\Z)|$, so $\kappa=1$, a contradiction. Assume $f$ has two singular fibers. (i) By \eqref{eq:pairs-general Pic} we have $$\un c_1\wt {\un c}_1=d\cdot |H_1(S',\Z)|=\wt {\un c}_1\wt \kappa \cdot |H_1(S',\Z)|,$$ so $\un c_1=\wt \kappa \cdot |H_1(S',\Z)|$ and analogously $\wt {\un c}_1=\kappa \cdot |H_1(S',\Z)|$. (ii) If, say, $\wt h=1$ then by definition $\wt {\un c}_1=1$, so again $\kappa=1$, a contradiction.
(iii) By \eqref{eq:pairs-general Pic} and \eqref{eq:pairs-general-d-is-lcm} $|H_1(S',\Z)|=\un c_1\wt {\un c}_1/lcm(\un c_1,\wt {\un c}_1)=gcd(\un c_1,\wt {\un c}_1)$, so (iii) follows from \ref{1prop:basic S' properties}(iv).

(iv) Since $f$ is pre-minimal, contractions in $\varphi_f$ are subdivisional with respect to $D^\dagger$, hence $$K_{\ov S^\dagger}\cdot(K_{\ov S^\dagger}+D^\dagger) =K\cdot(K+D)=3-\epsilon-K\cdot E.$$ Contract singular fibers to smooth fibers without touching $H$, denote the image of $D$ by $\wt D$ and the resulting Hirzebruch surface by $\wt S$. We have $$K_{\wt S}\cdot(K_{\wt S}+\wt D)=K_{\wt S}^2+K_{\wt S}\cdot H+2K_{\wt S}\cdot F=8+n-2-4=n+2.$$ A blow-down which is sprouting for a divisor $T$ increases $K\cdot(K+T)$ by one, so $$K^\dagger\cdot (K^\dagger+D^\dagger+C+\wt C+\Delta)+h+\wt h=K_{\wt S}\cdot(K_{\wt S}+\wt D)$$ and we get (iv).  \end{proof}

We will see that in case $\#E=1$ one can always find a pre-minimal affine ruling of $S\setminus\Delta$, often having additional good properties. We follow the original notation of \cite[5.3]{KR-C*onC3}.

\begin{figure}[h]\centering\includegraphics[scale=0.5]{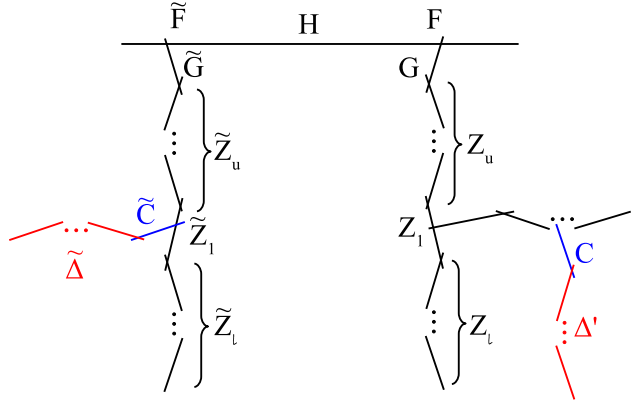}\caption{Notation for affine rulings of $S\setminus\Delta$.}  \label{fig:pre-minimal ruling}\end{figure}

\begin{notation}\label{2no:premin ruling notation} Assume $\#E=1$. Let $f:(\ov S^\dagger,D^\dagger+\Delta)\to\PP^1$ be a minimal completion of an affine ruling of $S\setminus\Delta$. We have $\Sigma_{S_0}=h+\nu-2=\nu=0$ by \ref{1prop:basic S' properties}(vii), because $E$ is irreducible and horizontal. Let $H^2=-n$, where $H$ is the horizontal component of $D^\dagger$. If $\beta_{D^\dagger}(H)>2$ then $(\ov S^\dagger,D^\dagger)=(\ov S,D)$ and the ruling is pre-minimal. Assume $\beta_{D^\dagger}(H)\leq 2$. If $n=1$ then $D^\dagger$ is not snc-minimal. In any case by successive contractions of exceptional curves in $D^\dagger$ (and its images) we obtain a morphism $\varphi_f:\ov S^\dagger\to \ov S$. Let $F$ be a singular fiber of $f$, such that $F\cap D^\dagger$ is branched. Denote the component of $F$ meeting $H$ by $G$. Let $Z$ be the chain consisting of curves produced by the first characteristic pair of $F$ and let $Z_1$ be the curve of highest multiplicity in $Z$. Let $Z_u$ and $Z_l$ (upper, lower) be the connected components of $Z-Z_1$ with $Z_u$ meeting $G$ (see Fig.~\ref{fig:pre-minimal ruling}). Let $Z_{lu}$ be the component of $Z_l$ meeting $Z_1$ and $C$ the unique $(-1)$-curve of $F$. Let $h$ be the number of characteristic pairs of $F$ and $\mu$ the multiplicity of $C$. If there is another singular fiber denote it by $\wt F$. Analogously for $\wt F$ define $\wt G, \wt Z_1, \wt h$, etc. Put $H^\dagger=Z_u+G+ H+\wt G+\wt Z_u$. Define $\Delta'=\Delta\cap F$ and $\wt \Delta=\Delta\cap \wt F$. \end{notation}

\bdfn\label{2def:almost minimal} In the situation as above $f$ is \emph{almost minimal} if $\varphi_f$ does not touch vertical $S_0$-components. \edfn

\bsrem By \ref{2cor:two fibers equations} $f$ has at least two singular fibers. If it has more than two then $\beta_{D^\dagger}(H)>2$ because each singular fiber contains a $D^\dagger$-component, hence $D^\dagger=D$ is snc-minimal, so $\varphi_f=id$ and $f$ is almost (and pre-) minimal. If $f$ is almost minimal with two singular fibers two then $h,\wt h\geq 2$ by \ref{2cor:two fibers equations} and the contractions in $\varphi_f$ take place within $H^\dagger$. It follows that an almost minimal ruling is pre-minimal.\esrem

\bprop[Koras-Russell, {\cite[5.3]{KR-C*onC3}}] \label{2prop:pre-minimal ruling exists} Let $C$ be a $(-1)$-curve in $\ov S$, such that $\kp(K_{\ov S}+D+\Delta+C)=-\8$. Then there exists a pre-minimal affine ruling of $S\setminus\Delta$ with $C$ in a fiber, such that either\benum[(i)]
\item $f$ is almost minimal or
\item $f$ has exactly two singular fibers, $\wt \Delta=0$ and $\varphi_f$ contracts precisely $H^\dagger+\wt Z_1$. If $Z_1$ is touched $x$ times in this process then $x\geq 4$ and $\wt V^2=2-x$, where $\wt V\subseteq D$ is the birational transform of $\wt Z_{lu}$. \eenum \eprop

Having the results established above the proof of the proposition and of all preliminary results (except 5.3.3(i) loc. cit, which is not necessary) goes without modifications as in loc. cit. The proposition implies that we have a good control over curves that are contracted when minimalizing the boundary. Note that in case (ii) $\wt Z_{lu}^2=1-x$ (as $\wt Z_{lu}$ is touched once in the contraction process), $\wt F$ has two characteristic pairs and the second is $\binom{1}{1}$.

\bcor\label{2cor:S-Delta pre-minimally ruled} If $\#E=1$ then there exists a pre-minimal affine ruling of $S\setminus\Delta$ with properties as in \ref{2prop:pre-minimal ruling exists}. \ecor

\begin{proof} Consider a minimal completion of some affine ruling of $S\setminus\Delta$. Since at least one of the branching components of $D^\dagger$ remains branching in $D$, there exists a singular fiber $F$, such that its $S_0$-component $C$ is not touched by the minimalization of $D^\dagger$ to $D$. By \ref{1lem:Eff-NegDef is Eff} we have $\kappa(K_{\ov S}+D+C+\Delta)=\kappa(K_{\ov S}+D+C+\Delta\cap F)$, because $\Delta-\Delta\cap F$ has a negative definite intersection matrix and its components intersect $K_{\ov S}+D+C+\Delta\cap F$ trivially. The snc-minimalization of a divisor or adding to a divisor a $(-1)$-curve intersecting it transversally in one point do not change the Kodaira dimension of the divisor, hence $\kappa(K_{\ov S}+D+C+\Delta\cap F)=\kappa(K_{\ov S}+D)=-\8$. Thus we can apply \ref{2prop:pre-minimal ruling exists}. \end{proof}

\med\section{The boundary is a fork}\label{sec:D is a fork}

\blem\label{2lem:KE=1 for eps=2} If $\epsilon=2$ then $K\cdot E=1$.\elem

\begin{proof} Suppose $\epsilon=2$ and $K\cdot E\neq 1$, then $\E=[4]$ by \ref{1cor:possible KE and eps}. Let $f:(\ov S^\dagger,D^\dagger)\to \PP^1$ be a pre-minimal affine ruling of $S\setminus\Delta$ (we use the notation of \ref{2no:premin ruling notation}). Let $F_1,\ldots F_N$ be the singular fibers. Put $U=H+\un F_1+\ldots+\un F_N$. We have $\Sigma_{S_0}=0$ and by \ref{2cor:two fibers equations} $N\geq 2$. Let $h_i=h(F_i)$ be the number of characteristic pairs of $F$. By \ref{1prop:KE and eps} $D$ consists of $(-2)$-curves and $\Delta=\emptyset$. In particular, $h_i\geq 2$. Suppose $N>2$. Then $D^\dagger=D$.   If we contract all $F_i$'s to smooth fibers without touching $H$ we make $h_1+h_2+\ldots+h_N$ sprouting blow-downs inside $U$. Let $\wt D$ and $\wt K$ be the image of $D$ and the canonical divisor of the resulting Hirzebruch surface. We have $K\cdot(K+U)=K\cdot(K+D)-N=-1-N$ and $\wt K\cdot (\wt K+\wt D)=8+\wt K\cdot H-2N=8-2N$. We obtain $-1-N+h_1+\ldots+h_N=8-2N$. Therefore $N=3$ and $h_1=h_2=h_3=2$, hence $D$ has three maximal twigs and, since $D$ consists of $(-2)$-curves, they are all equal to $[2,2,2]$. By \eqref{eq:pairs-general-d-is-lcm} $\kappa(F_1)\cdot \un c_1(F_1)=d=lcm(2,2,2)=2$, so $\kappa(F_1)=1$, a contradiction with \ref{1prop:basic S' properties}(ii). Thus $N=2$.

Suppose $f$ is not almost minimal. Then $n=1$ and $\wt h=2$, so $h=4$. By \ref{2prop:pre-minimal ruling exists} $\varphi_f:\ov S^\dagger\to \ov S$ contracts precisely $H^\dagger+\wt Z_1$ and $Z_1$ is touched exactly $2-\wt V^2=4$ times, hence $Z_1^2=-6$. $D$ consists of $(-2)$-curves, so the second branch of $F$ (see the definitions after \ref{1lem:sum of mu and squares}) is now necessarily $[(5)]$ and the third $[1,2]$ (the first component, $[1]$, is a tip of $F$). We have also $Z_l=[(k)]$ and $\wt Z_l=[(m),3]$ for some  non-negative integers $k,m$, hence $G=[k+1]$ and $\wt G=[m+2]$. If $k\neq 1$ then $\wt G$ is contracted before $G$, so $m=0$ and we see that $Z_1$ is touched at most once, a contradiction. Therefore $k=1$ and then $m=1$. Then $D$ has two branching components meeting each other, $B_1$ and $B_2$, such that $D-B_1-B_2=T_1+T_2+T_3+T_4$, with $T_1\cdot B_1=T_2\cdot B_1=1$, $T_1=[2,2]$, $T_2=[2]$, $T_3=[2]$ and $T_4=[2,2,2,2]$. We compute $d(D)=-25$, which contradicts \ref{2cor:two fibers equations}(iii). Thus $f$ is almost minimal with two singular fibers.

We have now $Z_l=[(k)]$ and $\wt Z_l=[(m)]$ for some positive integers $k,m$, so $Z_u=\wt Z_u=0$, $\wt G=[m+1]$ and $G=[k+1]$. Suppose $n=1$. Then $(\wt h,h)=(2,4)$ or $(\wt h,h)=(3,3)$. Consider the case $(\wt h,h)=(2,4)$. Note that $\wt Z_1^2=-2$, so $\wt G$ is not contracted by $\varphi_f$, hence $m>1$. If $k\neq 1$ then $\varphi_f$ contracts only $H$, so $m=k=2$ and the second branch of $F$ is $[1,2,2]$. In this case $d(D)=-9$, a contradiction with \ref{1prop:basic S' properties}(iv). Therefore $k=1$. We get $m=3$ and $Z_1^2=-3$ and we infer that the second branch of $F$ is $[2,2]$ and the third is $[1,2]$. Thus $D$ has two branching components, $B_1$ and $B_2$, and $D-B_1-B_2=T_1+T_2+T_3+T_4$ with $T_1=[(5)]$, $T_2=[2]$, $T_3=[2]$ and $T_4=[2]$. We get $d(D)=-16$ and $gcd(\wt {\un c}_1,\un c_1)=4$, a contradiction with \ref{2cor:two fibers equations}(iii).
Consider the case $(\wt h,h)=(3,3)$. We can assume $k\geq m$. If $m=1$ and $k=2$ then the second branch of $\wt F$ is $[2,2,2]$ and the second branch of $F$ is $[2,2]$, $gcd(\wt {\un c}_1,\un c_1)=6$ and $d(D)=-36$, a contradiction with \ref{2cor:two fibers equations}(iii). If $m=1$ and $k=3$ then the second branch of $\wt F$ is $[2,2]$ and the second branch of $F$ is $[1,2]$, $gcd(\wt {\un c}_1,\un c_1)=4$ and $d(D)=-16$, a contradiction with \ref{2cor:two fibers equations}(iii). It follows that $m=k=2$. Then second branches of $\wt F$ and $F$ are both $[1,2]$, so $d(D)=-9$, again a contradiction with \ref{2cor:two fibers equations}(iii).

We have now $n=2$, so $(\wt h,h)=(2,5)$ or $(\wt h,h)=(3,4)$. Now $Z_l$, $\wt Z_l$, $G$ and $\wt G$ are irreducible $(-2)$-curves. If $(\wt h,h)=(2,5)$ then $gcd(\wt {\un c}_1,\un c_1)=2$ and the second branch of $F$ is $[1,2,2,2]$, hence $d(D)=-4$. If $(\wt h,h)=(3,4)$ then $gcd(\wt {\un c}_1,\un c_1)=2$, the second branch of $\wt F$ is $[1,2]$ and the second branch of $F$ is $[1,2,2]$, so $d(D)=-4$. In both cases we get a contradiction with \ref{2cor:two fibers equations}(iii). \end{proof}

To prove that $D$ is a fork we need the following lemma. Recall that $s$ is the number of maximal twigs of $D$.

\blem\label{2lem:premin and short D elimination} Assume $\#E=1$. \benum[(i)]

\item If no twig of $D$ of length $\geq 2$ contains a $(-2)$-tip then there exists an affine ruling of $S\setminus \Delta$ with no base points on $\ov S$.

\item If $s=4$ and $\Delta$ is connected then $D$ has a twig of length $\geq 2$. \eenum \elem

\begin{proof} (i) Let $f:(\ov S^\dagger,D^\dagger+\Delta)\to \PP^1$ be a minimal completion of a pre-minimal affine ruling of $S\setminus\Delta$. Suppose $D^\dagger\neq D$. Then $f$ has two singular fibers, $F$ and $\wt F$, and $n=1$ (cf. \ref{2no:premin ruling notation}). By \ref{2prop:pre-minimal ruling exists}(ii) we can assume that the components of $Z_l$ are not contracted by $\varphi_f$. Since $h\geq 2$, by our assumption about maximal twigs of $D$ either $Z_l=[2]$ or $Z_l$ has a $\leq(-3)$-tip, in any case $G=[2]$. Analogous argument holds for $\wt F$, hence $H$ meets two $(-2)$-curves in $D^\dagger$. Therefore $D$ contains a non-branching component with non-negative self-intersection, a contradiction with \ref{1prop:basic S' properties}(v).

(ii) Suppose that $s=4$ and all maximal twigs of $D$ are tips. Then $D^\dagger=D$ by the first part of the lemma. From the geometry of the ruling we see that $H$ does not intersect a branching component of $D$, so it cannot be a maximal twig of $D$. If $H$ is non-branching in $D$ then $D$ has at least two branching components, which being contained in fibers, cannot be $(-1)$-curves, a contradiction with \cite[4.2]{Palka-k(S_0)=0}. Thus $H$ is branching in $D$, so there are at least three singular fibers. Two of them (at least) do not contain a branching component of $D$, hence contain unique $D$-components by our assumption. Then they both contain a component of $\Delta$, so $\Delta$ is not connected. \end{proof}

\bprop\label{2prop:D is a fork} $D$ is a fork.\eprop

\begin{proof} Suppose $D$ is not a fork. We first show that $\E=[5]$, $\epsilon=1$ and $s=4$ and then we eliminate this case in several steps. We prove successive statements.

\medskip (1) $\#E=1$ and $\epsilon=1 $ or $2$.

\begin{proof} We have $\epsilon\neq 0$ by \ref{1prop:if eps=0 or E^ fork}(i). To prove $\#E=1$ we can assume $\epsilon=1$ by \ref{1cor:possible KE and eps}. Thus $\E$ is a chain by \ref{1prop:if eps=0 or E^ fork}(ii) and it satisfies $(s-4)d(\E)+d'(\E)+d'(\E^t)\leq 7$ by \ref{1lem:if eps<2}(iii). Using $2\leq K\cdot E\leq 3$ this gives only two cases for which $\#E\neq 1$: $s=4$ and $\E=[3,3]$ or $s=4$ and $\E=[3,4]$. By \ref{1lem:eps properties}(iv) in both cases $e+\delta<3$, which is impossible by \ref{1lem:if eps<2}(iv).\end{proof}

(2) If $K\cdot(K+D)\neq 0$ then $\E=[5]$, $\epsilon=1$ and $s=4$.

\begin{proof} Assume $K\cdot(K+D)\neq 0$. For $\epsilon=2$ we have $K\cdot(K+D)=3-\epsilon-K\cdot E=0$ by \ref{2lem:KE=1 for eps=2}, so $\epsilon=1$ by (1). Again by \ref{1lem:if eps<2}(iii) $(s-4)d(\E)+d'(\E)+d'(\E^t)\leq 7$, so since $K\cdot E=3$ and $\#E=1$, we obtain $s=4$ and $\E=[2,5]$ or $s\leq 5$ and $\E=[5]$. In the first case we have $e=\delta=\frac{4}{3}$ by \ref{1lem:eps properties}(iv) and \ref{1lem:if eps<2}(ii), so maximal twigs of $D$ are tips, a contradiction with \ref{2lem:premin and short D elimination}. Suppose $s=5$ in the second case. Then similarly $e=\delta=\frac{9}{5}$, which is impossible by \ref{1lem:if eps<2}(iv).\end{proof}

We choose a minimal completion $f:(\ov S^\dagger,D^\dagger)\to \PP^1$ of a pre-minimal affine ruling of $S\setminus\Delta$. Subdivisional modifications of $D$ do not change $K\cdot(K+D)$, so $K^\dagger\cdot(K^\dagger+D^\dagger)=K\cdot(K+D)$, where $K^\dagger=K_{\ov S^\dagger}$. According to \ref{2cor:two fibers equations} $f$ has at least two singular fibers.

\medskip (3) If $D^\dagger\cap F$ is not a chain for some fiber $F$ of $f$ then $K\cdot(K+D)\neq 0$.

\begin{proof} Suppose $F\cap D^\dagger$ is branched and $K\cdot(K+D)=0$. Write $F$ as $F=F\cap D^\dagger+C+\Delta'$, where $C$ is a $(-1)$-curve, and $\Delta'\subset \Delta$. We contract the chain $C+\Delta'$ and successive $(-1)$-curves in $F$ as long as they are subdivisional for $D^\dagger$. Denote the images of $D^\dagger$, $E$ and $F$ by $D^{(1)}$, $E^{(1)}$ and $F^{(1)}$. Let $K^{(1)}$ be the canonical divisor of the image of $\ov S$. In general, if after some sequence of contractions we define $D^{(i)}$ then we denote the respective images of $E$, $F$, etc. by $E^{(i)}$, $F^{(i)}$ etc. and the canonical divisor on the respective image of $\ov S$ by $K^{(i)}$. The contraction of $C+\Delta'$ and contractions subdivisional with respect to the image of $D^\dagger$ do not change $K^\dagger\cdot(K^\dagger+D^\dagger)$ and $E\cdot(K^\dagger+D^\dagger)$ (cf. \ref{1lem:blow-up and transforms}), i.e. $$K^{(1)}\cdot(K^{(1)}+D^{(1)})=K\cdot(K+D)=0 \text {\ and\ } E^{(1)}\cdot(K^{(1)}+D^{(1)})=E\cdot(K+D)=K\cdot E.$$ Moreover, $F^{(1)}\cap D^{(1)}$ is branched.

Let $D^{(1)}_\alpha$ be the $(-1)$-tip of $D^{(1)}$, and let $D^{(2)}$ be the image of $D^{(1)}$ after the contraction of $D^{(1)}_\alpha$. Let $D^{(1)}_\beta$ be the unique $D^{(1)}$-component intersecting $D^{(1)}_\alpha$. Note that $\kp(K^{(2)}+D^{(2)})=\ovk(S\setminus (C\cup \Delta))=\ovk(S)=-\8$, so since by the Riemann-Roch theorem $$h^0(-K^{(2)}-D^{(2)})+h^0(2K^{(2)}+D^{(2)})\geq K^{(2)}\cdot(K^{(2)}+D^{(2)})=1,$$ we get $-K^{(2)}-D^{(2)}\geq 0$. For every component $V$ of $D^{(2)}$ we have $V\cdot (-K^{(2)}-D^{(2)})=2-\beta_{D^{(2)}}(V)$. Since $s\geq 4$, $D^{(2)}$ is branched and every branching curve of $D^{(2)}$, and hence every component of $D^{(2)}$ which is not a tip, is in the fixed part of $-K^{(2)}-D^{(2)}$. Suppose $D^{(2)}_\beta$ is not a tip of $D^{(2)}$, then $-K^{(2)}-D^{(2)}-D^{(2)}_\beta\geq 0$, so $-K^{(1)}-D^{(1)}-D^{(1)}_\beta\geq 0$. Clearly, $E^{(1)}$ is in the fixed part of the latter divisor, so $-K^{(1)}-D^{(1)}-E^{(1)}\geq 0$. It follows that $-(K^\dagger+D^\dagger+E)\geq 0$, a contradiction with $\kp(K^\dagger+D^\dagger+E)=2$. Thus $D^{(2)}_\beta$ is a tip of $D^{(2)}$.

Let $D^{(3)}$ be the image of $D^{(2)}$ after the contraction of $D^{(2)}_\beta$. Since $D^{(2)}_\beta$ is a tip, $D^{(2)}$ has the same number of branching components as $D^{(1)}$ (greater than one by our assumptions about $D$), hence $D^{(3)}$ is not a chain. Moreover, $F^{(3)}$ is not a $0$-curve, as no branching component of $D^\dagger\cap F$ has been contracted. We made two sprouting blow-downs, so $$K^{(3)}\cdot(K^{(3)}+D^{(3)})= K^{(1)}\cdot(K^{(1)}+D^{(1)})+2=K\cdot(K+D)+2=2.$$ Riemann-Roch's theorem gives $h^0(-K^{(3)}-D^{(3)})\geq 2$. Since $f$ has at least two singular fibers, $H$ is not a tip of $D^{(3)}$. Since $D^{(3)}$ is not a chain, $H$ is in the fixed part of $-K^{(3)}-D^{(3)}$. Let's write $-K^{(3)}-D^{(3)}=H+R+M$, where $M$ is effective, $h^0(M)\geq 2$ and the linear system of $M$ has no fixed component.  Intersecting with a general fiber $F'$ we have $1=1+F'\cdot R+F'\cdot M$, so $F' \cdot M=F'\cdot R=0$ and $R$ and $M$ are vertical, hence $M\sim tF'$ for some $t>0$. We get that $K^{(3)}+D^{(3)}+H+t F'+R\sim0$. Intersecting with $E^{(3)}$ gives $$0\geq E^{(3)}\cdot(K^{(3)}+D^{(3)}+F')= E^{(2)}\cdot(K^{(2)}+D^{(2)}-D^{(2)}_\beta+F')=$$ $$=E^{(1)}\cdot(K^{(1)}+D^{(1)})+E^{(1)}\cdot(F'-2D_\alpha^{(1)}-D_\beta^{(1)})= K\cdot E+E^{(1)}\cdot(F_0^{(1)}-2D_\alpha^{(1)}- D_\beta^{(1)}),$$ which implies $E^{(1)}\cdot(F^{(1)}-2D_\alpha^{(1)} -D_\beta^{(1)})<0$. This is a contradiction, because $F^{(1)}$ is branched, so the multiplicities of $D_\alpha^{(1)}$ and $D_\beta^{(1)}$ in it are greater than one.\end{proof}

(4) $\E=[5]$, $\epsilon=1$ and $s=4$.

\begin{proof} Suppose (4) does not hold. Then by (2) and (3) $H$ is the only branching curve in $D^\dagger$, so $D^\dagger=D$, every singular fiber $F$ of $f$ has at most one branching component and $F\cap D$ is a chain. In particular, there are exactly $s$ singular fibers. Let $c$ be the number of singular fibers which are chains. If $F$ is such a fiber then $F\cap\Delta\neq\emptyset$ and $F\cap D$ is a tip, so $\wt e(F\cap D)\leq \frac{1}{2}$. Since $s\geq 4$ and since $\Delta$ has at most three connected components, we see that $c<s$, so we have an inequality $\wt e(D)<(s-c)+\frac{c}{2}=s-\frac{c}{2}$. Let's contract all singular fibers to smooth $0$-curves without touching $H$. The contraction of chain fibers does not affect $K\cdot(K+D)$ and the contraction of any other singular fiber increases $K\cdot(K+D)$ by one, so if $\wt D$ and $\wt S$ are the images of $D^\dagger$ and $\ov S^\dagger$ after contractions then $\wt D\equiv H+sF'$ for a general fiber $F'$ and $$K_{\wt S}\cdot(K_{\wt S}+\wt D)=K\cdot (K+D)+s-c=s-c.$$ We get $s-c= K_{\wt S}\cdot(K_{\wt S}+ \wt D)=8-H^2-2-2s$, so $n=-H^2=3s-c-6$. By the Laplace expansion we have (cf. \cite[2.1.1]{KR-ContrSurf}) $d(D)=d_1\cdot\ldots\cdot d_s(n-\wt e(D))$, where $d_i$ are discriminants of maximal twigs, so by \ref{1prop:basic S' properties}(iv) $\wt e(D)>n$. Thus $s-\frac{c}{2}> \wt e(D)>3s-c-6$, so $12>4s-c>3s$ and then $s\leq 3$, a contradiction. \end{proof}

\noin Recall that $T$ is the sum of maximal twigs of $D$.

\medskip (5) \begin{minipage}[t]{0.93\textwidth} If $R\subseteq D$ is a $\leq(-4)$-tip of $D$ then for every irreducible component $V$ of $T$ we have $0\leq V\cdot (2K+R)\leq 1$ and for at most one $V\cdot (2K+R)\neq 0$. \end{minipage}

\begin{proof} Let $m$ be a maximal natural number, such that $E+m(K+D)\geq 0$. It exists by \ref{1lem:linear systems}(iii) and is greater than one by (4) and \ref{1lem:if eps<2}(i). By \ref{1lem:linear systems}(ii) we can write $E+m(K+D)=\sum C_i$, where $C_i\cong \PP^1$ and $C_i^2< 0$. Moreover, $C_i\neq E$, as $\kappa(K+D)=-\8$. Multiplying both sides by $E+2K+R$ we have $$K\cdot E-2+m(4-2\epsilon-K\cdot E+R(D-R))=\sum_iC_i\cdot (E+2K+R),$$ so $\sum_iC_i\cdot (E+2K+R)=1$ by (4). Suppose $C_{i_0}\cdot (E+2K+R)<0$ for some $i_0$. If $C_{i_0}\cdot K\geq 0$ then we get $C_{i_0}=R$ and $0>R\cdot (2K+R)=R\cdot K-2$, which is impossible by our assumption on $R$. Thus $C_{i_0}\cdot K<0$. Then $C_{i_0}^2=-1$ and $C_{i_0}\cdot (E+R)\leq 1$. Simultaneously $|K+D+C_{i_0}|=\emptyset$ by the definition of $m$, so by \ref{1lem:linear systems}(i) $D\cdot C_{i_0}\leq 1$. Thus either $C_{i_0}$ is simple or it is a non-branching $(-1)$-curve in $D$, a contradiction.  Therefore $C_i\cdot (E+2K+R)\geq 0$ for each $i$. If $V$ is a component of $T$ then $V\cdot (E+m(K+D))=m(\beta_D(V)-2)$, so tips of $D$, and hence all components of $T$, appear among $C_i$'s and we are done. \end{proof}

(6) There are no $\leq (-4)$-tips in $D$.

\begin{proof} Suppose $T_1$ contains a $\leq -4$-tip of $D$, denote it by $R$. By (5) $T-R$ consists of $(-2)$-curves and $-5\leq R^2\leq -4$. Maximal twigs of $D$ other than $T_1$ are tips, otherwise $e\geq \frac{1}{5}+\frac{1}{2}+ \frac{1}{2}+\frac{2}{3}>\frac{9}{5}$, a contradiction with \ref{1lem:eps properties}(iv). If $R^2=-5$ then $V\cdot (2K+R)=0$ for every component of $T-R$, so $R$ is a maximal twig, a contradiction with \ref{2lem:premin and short D elimination}. Thus $T_1=[4,(k-1)]$ for some positive integer $k$, hence by \ref{1lem:eps properties}(iv) $\frac{9}{5}\geq e=\frac{3}{2}+\frac{1}{3+1/k}$, so $k\leq 3$. By \ref{2lem:premin and short D elimination} there is an affine ruling of $S\setminus\Delta$ which extends to a $\PP^1$-ruling $f$ of $(\ov S,D)$. If $F$ is a singular fiber of $f$ then, since $\Delta=\emptyset$, $D\cap F$ contains at least four components, otherwise we would have $F\cap D=[2,2,2]$, which is impossible by the description of maximal twigs. Thus for every singular fiber $F$ the divisor $F\cap D$ is branched, so by \ref{2cor:two fibers equations} $f$ has two singular fibers, $h,\wt h\geq 3$ and $h+\wt h=n+5$. Since $Z_{l}$ and $\wt Z_l$ are equal to $[4,(k-1)]$ or $[2]$, $G=[2]$ and $\wt G=[2]$, so $n>1$ by \ref{1prop:basic S' properties}(v). This implies that one of $h$ or $\wt h$, say $h$, is at least $4$, so the second branch of the respective singular fiber contains at least two $D$-components, hence contains $T_1$. Let $C$ be the unique $S_0$-component of $F$. Now $T_1+C$ should contract to a smooth point. This is possible only for $k=4$, a contradiction. \end{proof}

(7) Maximal twigs of $D$ are $[2]$, $[2]$, $[3]$ and $[3,2]$.

\begin{proof} We assume that $d_1\leq d_2\leq d_3\leq d_4$. By \ref{1lem:eps properties}(iv) and \ref{1lem:if eps<2}(iv) we have $e\leq \frac{9}{5}$ and $\delta\geq \frac{13}{4}-e\geq \frac{13}{4}-\frac{9}{5}=\frac{29}{20}$, so $d_1=2$ and $2\leq d_2\leq 3$. If $d_2=3$ then the lower bound on $\delta$ gives $d_3=d_4=3$, and since by \ref{2lem:premin and short D elimination} not all maximal twigs are tips, $e\geq \frac{1}{2}+\frac{1}{3}+\frac{1}{3}+\frac{2}{3}>\frac{9}{5}$, a contradiction. Thus $d_2=2$ and we have $\frac{1}{d_3}+\frac{1}{d_4}\geq \frac{9}{20}$, so $d_3\leq 4$. Since there are no $(-4)$-tips in $D$ by (6), $e_4>\frac{1}{3}$, so for $d_3=4$ we get $e\geq 1+\frac{3}{4}+\frac{1}{3}>\frac{9}{5}$, which is impossible. Thus $d_3\leq 3$. In fact $T_3=[3]$, otherwise $e\geq \frac{3}{2}+\frac{1}{3}>\frac{9}{5}$. We get $d_4\leq 8$ and $e_4\leq \frac{9}{5}-1-\frac{1}{3}<\frac{1}{2}$, so $T_4$ contains a $(-3)$-tip, hence $T_4=[3,3]$ or $T_4=[3,(k)]$ for some $k\in\{0,1,2\}$. Only $T_4=[3]$ and $T_4=[3,2]$ satisfy \ref{1lem:if eps<2}(iv), so other cases are excluded. The case $T_4=[3]$ is excluded by \ref{2lem:premin and short D elimination}. \end{proof}

\noin Now we see by \ref{2lem:premin and short D elimination} that there is an affine ruling $f$ of $(\ov S,D)$. As in (6) we see that $f$ has two singular fibers and the second branch of one of them consists of an $S_0$-component $C$ and $T_4$. Now again $T_4+C$ should contract to a smooth point. But this is impossible for $T_4=[3,2]$, a contradiction. \end{proof}

\blem\label{2lem:Zar and others for D fork} Let $\cal P=(K+D+\E)^+$ and let $B$ be the branching component of $D$. Put $b=-B^2$. Then:\benum[(i)]

\item $b\in\{1,2\}$ and $b<\wt e$,

\item $\delta<1$,

\item $\cal P\equiv\frac{1-\delta}{\wt e -b}(B+\sum_{i=1}^3\Bk' T_i^t)$,

\item $\Bk^2 \E=-\frac{(1-\delta)^2}{\wt e-b}+e-1-\epsilon$.\eenum\elem

\begin{proof}

(i) $0>d(D)=d_1 d_2 d_3(b-\wt e)\geq b-\wt e$ by \ref{1lem:ForkBk}(iv) and \ref{1prop:basic S' properties}(iv). Now $\wt e_i<1$, so $b<\wt e<3$ and we get $b\in\{1,2\}$ by \ref{1prop:basic S' properties}(v).

(ii) $\cal P\cdot V=0$ for every component $V$ of $T+\E$, because $T+\E \subseteq (K+D+\E)^-$. Components of $D+\E$ generate $\Pic \ov S\otimes \Q$ by \ref{1prop:basic S' properties}(vi), so $\cal P\cdot B\neq0$, otherwise $\cal P\equiv0$, which contradicts $\ovk(S_0)=2$. We infer that $$0<B\cdot \cal P=B\cdot(K+D-\Bk D)=1-\delta.$$

(iii) Both $\cal P$ and $B+\sum_{i=1}^3\Bk' T_i^t$ intersect trivially with all components of $T+\E$, so they are linearly dependent in $\Pic \ov S\otimes\Q$. Moreover $\cal P\cdot B=1-\delta$ and $(B+\sum_{i=1}^3\Bk' T_i^t)\cdot B=\wt e-b$.

(iv) We compute $$\cal P^2=\frac{(1-\delta)^2}{(\wt e-b)^2}(B^2+\sum_{i=1}^3 \wt e_i)=\frac{(1-\delta)^2}{\wt e-b},$$ so since $\Bk^2 D=-e$, (iv) follows from \ref{1prop:basic S' properties}(ii). \end{proof}

\brem\label{2rem:finite number of cases} If $K\cdot T$ is bounded (for example this is the case when we can bound the determinants $d_1,d_2,d_3$) then there is only finitely many possibilities for the weighted dual graphs of $D$ and $\E$. Indeed, by \ref{1prop:KE and eps} and \ref{2lem:KE=1 for eps=2} $K\cdot E+\epsilon\leq 5$ and by \ref{2lem:Zar and others for D fork}(i) $b\in\{1,2\}$, so $K\cdot E+K\cdot D$ is bounded. It is therefore enough to bound $\#\E+\#D$. This is possible using Noether formula \ref{1lem:eps properties}(iii). \erem

\blem\label{2lem:three fibers only for b=2} If $b=\#E=1$ then any affine ruling of $S\setminus\Delta$ has two singular fibers.\elem

\begin{proof} Let $f:(\ov S^\dagger,D^\dagger+\Delta)\to \PP^1$ be a minimal completion of an affine ruling of $S\setminus\Delta$. We have $\Sigma_{S_0}=0$, because $\#E=1$. By \ref{2cor:two fibers equations} $f$ has more than one singular fiber. Suppose it has more than two singular fibers. Each singular fiber contains a $D$-component, so we infer that $D^\dagger=D$, $B$ is horizontal and $f$ has three singular fibers $F_1,F_2,F_3$. Let $C_i$ and $\Delta_i$ for $i=1,2,3$ be respectively the $S_0$-component and the connected component of $\Delta$ contained in $F_i$ (it is possible that $\Delta_i=0$). By \ref{1lem:linear systems}(iii) there exists a greatest integer $m$, such that $B+m(K+D)\geq0$. By \ref{1lem:linear systems}(i) $m\geq 1$, because $B\cdot D=3-b>1$. Write $B+m(K+D)\sim L$ with $L$ effective. Multiplying by a general fiber $F'$ we get $1-m=F'\cdot L\geq 0$, so $m=1$ and $L$ is vertical. Denote the $D$-component of $D$ intersecting $B$ by $D_i$. Denote the number of characteristic pairs of $F_i$ by $h_i$ and assume $h_1\leq h_3\leq h_3$. Note that for any component $D_0$ of $D$ we have $D_0\cdot(K+D)=-2+\beta_D(D_0)$, so all components of $D-B-D_1-D_2-D_3$ are contained in $L$. Now if $h_i\neq 1$ then $C_i+\Delta_i\subseteq L$. Indeed, if $h_i\neq 1$ then $C_i\cdot (K+D+B)=0$ and the $D$-component intersecting $C_i$ is contained in $L$, hence so is $C_i$ and then by induction all components of $\Delta_i$. By \ref{1prop:basic S' properties}(ii) $E\cdot (C_i+\Delta_i)\geq 2$ for each $i$, so $h_1=1$, otherwise $$K\cdot E=E\cdot (K+D+B)=E\cdot L\geq \sum_{i=1}^3 E\cdot (C_i+\Delta_i)\geq 6,$$ which contradicts \ref{1cor:possible KE and eps}. It follows that $\Delta\neq \emptyset$, hence $\epsilon\neq 0$ by \ref{1prop:if eps=0 or E^ fork}. Then $K\cdot E\leq 3$ by \ref{1cor:possible KE and eps}, so as above we infer that $h_2=1$. By \ref{2prop:pairs}(4) $d=\un c_1(F_3)$, so $\kappa_3=1$ and $C_3$ is simple on $(\ov S,D)$, a contradiction. \end{proof}

\bcor\label{2cor:Delta with three components} If $\Delta$ has three connected components then $b=\epsilon=2$.\ecor

\begin{proof} If $\Delta$ has three connected components then $\E$ is a fork, so $\epsilon=2$ by \ref{1prop:if eps=0 or E^ fork}(ii) and $\#E=1$ by \ref{2lem:KE=1 for eps=2}. Each connected component of $\Delta$ is contained in a different singular fiber of a minimal completion of an affine ruling of $S\setminus \Delta$. By \ref{2lem:three fibers only for b=2} and \ref{2lem:Zar and others for D fork}(i) $b=2$. \end{proof}

\med\section{Some intermediate surface containing the smooth locus}\label{sec:surface W}

Recall that $T=D-B$, where $B$ is the branching component of $D$. We define $W=\overline{S}-T-\E$. Clearly, $S_0=W\setminus B$ and hence $\chi(W)=\chi(S_0)+\chi(\C^{**})=-1$. Our goal is to prove that $\ovk(W)=2$. This takes a lot of work but allows later to strongly restrict possible shapes of $\E$ using the logarithmic Bogomolov-Miyaoka-Yau inequality. To achieve the goal we prove couple of technical lemmas and use results of computer programs.

\blem\label{3lem:eT+x/d=1 solutions}Let $R$ be an oriented admissible chain and let $\alpha$ be such that \begin{equation} e(R)+\frac{\alpha}{d(R)}=1. \tag{*}\end{equation} Then:\benum[(i)]

\item $R=[2,\ldots,2,2]$ or $R=0$ if and only if $\alpha=1$,

\item $R=[2,\ldots,2,3]$ if and only if $\alpha=2$,

\item $R=[2,\ldots,2,3,2]$ or $R=[2,\ldots,2,4]$ if and only if  $\alpha=3$.\eenum\elem

\begin{proof} Note that by \ref{1lem:d() expansion} we have a recurrence formula $d([a_1,a_2,\ldots,a_k])= a_1d([a_2,\ldots,a_k])-d([a_3,\ldots,a_k])$. Using it we see that $R=[2,a_1\ldots,a_k]$ satisfies (*) if and only if $[a_1,\ldots,a_k]$ does, so we may assume that $R=[a_1,\ldots,a_k]$ with $a_1\geq 3$. If the equation holds then $d'(R)+\alpha=d(R)=a_1 d'(R)-d''(R)$, so $2d'(R)\leq(a_1-1)d'(R)=d''(R)+\alpha< d'(R)+\alpha$, hence $d'(R)<\alpha\leq3$ and $k\leq2$. For $d'(R)=2$ we get $R=[3,2]$, for $d'(R)=1$ we get $R=[4]$ or $R=[3]$ and for $d'(R)=0$ we get $R=0$. \end{proof}

\blem\label{3lem:bounds for e(F)} If $R=[(k),c,a_1,\ldots,a_n]$ is admissible then $$\frac{k(c-1)+1}{k(c-1)+c}\leq e(R)<\frac{k(c-2)+1}{k(c-2)+c-1}.$$\elem

\begin{proof} For a chain $R=[u,\ldots]$ we have $d(R)=ud'(R)-d''(R)$ and hence $e(R)=\frac{1}{u-e'(R)}$. Since $0\leq e'(R)<1$, we get $\frac{1}{c}\leq e(R)<\frac{1}{c-1}$. The formula for $k\neq 0$ follows by induction. \end{proof}

\blem\label{3lem:$W$ properties}\ \benum[(i)]

\item $W$ is almost minimal and $K+T+\E\equiv\lambda \cal P+\Bk T+\Bk \E$, where $\lambda=1-\frac{\wt e -b}{1-\delta}$.

\item If $\ovk(W)\geq0$ then $\lambda \cal P\equiv(K+T+\E)^+$.

\item If $\ovk(W)\geq 0$ then $\wt e+\delta\leq b+1$, $\delta+\frac{1}{|G|}\geq 1$ and $\epsilon\neq 0$. The inequalities are strict if $\ovk(W)=2$.

\item If $\ovk(W)\neq 2$ then $\ovk(W)\leq 0$, $\wt e+\delta\geq 2$ and $b=1$. The inequality is strict if $\ovk(W)=-\8$.

\item If $K\cdot T_i=0$ for some $i$ then $h^0(2K+T+\E)\geq 3-b-\epsilon$. \eenum\elem

\begin{proof} (i) Recall that $\Bk T_i=\Bk' T_i+\Bk' T_i^t$. Using \ref{2lem:Zar and others for D fork}(iii) we have $$K+T+\E\equiv\cal P-B+\Bk D+\Bk \E=$$ $$=\cal P-B-\sum_{i=1}^3 \Bk' T_i^t+\sum_{i=1}^3 \Bk T_i+\Bk \E=(1-\frac{\wt e-b}{1-\delta})\cal P+\Bk T+\Bk \E.$$ Suppose $W$ is not almost minimal. Then by \cite[2.3.11]{Miyan-OpenSurf} there exists a $(-1)$-curve $C$, such that $C+\Bk \E+\Bk T$ has negative definite intersection matrix. Since the support of $\Bk \E+\Bk T$ is $\E\cup T$, $(K+T+\E)^-$ has at least $\#T+\#\E+1=b_2(\ov S)$ numerically independent components (cf. \ref{1prop:basic S' properties}(vi)), a contradiction with the Hodge index theorem.

(ii) From (i) and from the definition of $\Bk$ we see that $\cal P$ intersects trivially with every component of $T+\E$. If $\ovk(W)\geq 0$ then by the properties of Fujita-Zariski decomposition the same is true for $(K+T+\E)^+$. Since $\Pic \ov S\otimes \Q$ is generated by the components of $D+\E$, we get $(K+T+\E)^+\equiv \alpha \cal P$ for some $\alpha\in \Q$. We have $\cal P\cdot B=1-\delta$ and $$(K+T+\E)^+\cdot B=(K+T+\E)\cdot B-\Bk T\cdot B=b+1-\wt e-\delta,$$ hence intersecting with $B$ gives $\alpha=\lambda$.

(iii) We have $\chi(W)=-1$, so $\delta+\frac{1}{|G|}\geq 1+\frac{1}{3}\lambda^2\cal P^2$ by \ref{1lem:KobIneq}(ii). By (ii) and \cite[6.11]{Fujita} $\ovk(W)> 0$ ($\ovk(W)=0$) if and only if $\lambda>0$ (respectively $\lambda=0$), which is equivalent to $b+1> \wt e+\delta$ (respectively $b+1=\wt e+\delta$). Suppose $\epsilon=0$. Then $\E=[|G|]$ by \ref{1prop:if eps=0 or E^ fork}(i), so by \ref{1lem:eps properties}(iv) $\delta+\frac{1}{|G|}\leq e+\frac{1}{|G|}\leq 1$. Together with the inequality above this implies $e=\delta$, so maximal twigs of $D$ are tips, a contradiction with \ref{1lem:eps properties}(iii).

(iv) Suppose $\ovk(W)=1$. Then by (ii) $\lambda^2\cal P^2=0$, so $\lambda=0$ and hence $(K+T+\E)^+\equiv 0$ and $\ovk(W)=0$ by \cite[6.11]{Fujita}, a contradiction. Thus $\ovk(W)\leq 0$. Note that if $\ovk(W)=-\8$ then $\kappa(K+D+T)=-\8$ and by rationality of $W$ the divisor $K+T+\E$ cannot be numerically equivalent to an effective divisor, hence $\lambda<0$. Thus for $\ovk(W)\leq 0$ we have $b+1\leq\wt e+\delta$ and the inequality is strict for $\ovk(W)=-\8$. Suppose $b=2$. Since $\wt e_i+\frac{1}{d_i}\leq 1$, we get $\wt e_i+\frac{1}{d_i}=1$ for each $i$, so $D$ consist of $(-2)$-curves by \ref{3lem:eT+x/d=1 solutions}(i). By \ref{2lem:Zar and others for D fork}(iv) $0>\Bk^2 \E=1-\epsilon$, so $\epsilon=2$, $\E$ is a chain by \ref{1lem:ForkBk}(v) and $d'(\E)+d'(\E^t)+2=d(\E)$. By \ref{3lem:bounds for e(F)} if $\Delta$ is not connected then $e(\E),\wt e(\E)\geq \frac{1}{2}$, so $d'(\E)+d'(\E^t)\geq d(\E)$. Thus $\Delta$ is connected and by \ref{2lem:KE=1 for eps=2} $\E=[3,(k)]$ for some $k\geq 0$. Then $d'(\E)+d'(\E^t)+2-d(\E)=k+1$, a contradiction.

(v) Assume $K\cdot T_1=0$. Riemann-Roch's theorem gives $$h^0(-K-T_2-T_3-\E)+h^0(2K+T_2+T_3+\E)\geq$$ $$\geq\frac{1}{2}(K+T_2+T_3+\E)\cdot (2K+T_2+T_3+\E)+1=3-\epsilon-b.$$ If $-K-T_2-T_3-\E\geq 0$ then $B$, and hence $T_1$, is in the fixed part, so $-K-D-\E\geq 0$, which contradicts $\ovk(S_0)=2$. Thus $h^0(2K+T_2+T_3+\E)\geq 3-b-\epsilon$. \end{proof}

\bprop\label{3prop:no small 0-chains in D} If $D$ contains $[2,1,2]$ or $[3,1,2,2]$ then $\#E>1$ and $\ovk(W)=2$.\eprop

\begin{proof} Assume $D$ contains $F_\8=[2,1,2]$ or $F_\8=[3,1,2,2]$. Since $D$ is snc-minimal, the $(-1)$-curve of $F_\8$ is $B$, the branching component of $D$. The divisor $F_\8$ snc-minimalizes to a $0$-curve, hence gives a $\PP^1$-ruling $p:\ov S\to \PP^1$ with $F_\8$ as a fiber. $\E$ is vertical because $F_\8\cdot \E=0$, so $\Sigma_{S_0}=h+\nu-2=h-1\leq2$. Denote the fiber of $p$ containing $\E$ by $F_E$. We have $F_E\cdot D\leq 5$ because $\mu(B)\leq 3$. Note that for every $S_0$-component $L$ we have $L\cdot \E\leq 1$, because $F_E$ is a tree, so by \ref{1prop:basic S' properties}(ii) $\# L\cap D\geq 2$. There are no $(-1)$-curves in $D$ other than $B$, so all vertical $(-1)$-curves are $S_0$-components. We prove successive statements.

\medskip (1) If $\ovk(W)\neq 2$ then $E=[3]$.

\begin{proof} Suppose $\ovk(W)\neq 2$. By \ref{3lem:$W$ properties}(iv) $\ovk(W)\leq 0$, $\wt e+\delta\geq 0$ and $\lambda\leq 0$. We first show that all $S_0$-components are exceptional. For any $S_0$-component $L$ we have $L\cdot(K+T^\#+\E^\#)=\lambda L\cdot \cal P$. By \ref{2lem:Zar and others for D fork} $\Supp \cal P=D$, so $L\cdot \cal P>0$ because $L\cdot D>0$. Suppose $L^2\leq -2$. Then $L\cdot (T^\#+\E^\#)\leq\lambda L\cdot \cal P$, which, since $\lambda\leq 0$, is possible only if $\lambda=L\cdot T^\#=L\cdot \E^\#=0$. If $L$ intersects at least two twigs of $D$, say, $T_1$ and $T_2$ then $L\cdot T^\#=0$ implies that $T_1^\#=T_2^\#=0$, so $T_1$ and $T_2$ are $(-2)$-chains and then $\lambda=0$ gives $\wt e_3+\frac{1}{d_3}=0$, which is impossible. Thus $L\cdot T_1=L\cdot T_2=0$ and $\# L\cap T_3\geq 2$, which implies that $T_3$ contains the multiple section of $D$ and, as before, that it consists of $(-2)$-curves. We get $\wt e_3+\frac{1}{d_3}=1$ and now $\lambda=0$ gives $\wt e_1+\wt e_2<1$. However, by \ref{3lem:bounds for e(F)} in case $F_\8=[3,1,2,2]$ we have $\wt e_1+\wt e_2\geq \frac{1}{3}+\frac{2}{3}=1$ and in case $F_\8=[2,1,2]$ we have $\wt e_1+\wt e_2\geq \frac{1}{2}+\frac{1}{2}=1$, a contradiction.

Let $D_h$ and $D_v$ be respectively the divisor of horizontal components of $D$ and the divisor of $D$-components contained in $F_E$. Let $D_1$ be the multiple section contained in $D_h$. Denote the $S_0$-components of $F_E$ by $L_1,L_2,\ldots,L_{\sigma(F_E)}$. Clearly, $D_v$ has at most three connected components and they are chains. We prove that $D_h$ contains a section and $D_v\neq 0$. Suppose $D_h$ does not contain a section. In this case $D_h$ is irreducible, so $\Sigma_{S_0}=0$ and $\sigma(F_E)=1$. We have now $F_E\cdot D\leq 3$ and $\mu(L_1)\geq 2$, so since $\# L_1\cap D\geq 2$, $D_h$ intersects $L_1$ in exactly one point and $D_v\neq 0$. This gives $$\mu(L_1)+1\leq F_E\cdot D_h\leq 3,$$ so $\mu(L_1)=2$ and we get $\E=[2]$, a contradiction. Suppose $D_v=0$. Since $\#L_i\cap D\geq 2$ for each $i$, $\sigma(F_E)\leq 2$. As $D_h$ contains a section, the $S_0$-component intersecting it, say $L_1$, has multiplicity one, so $\sigma(F_E)=2$. Then $\mu(L_2)=1$, otherwise $L_2$ could intersect no other component of $D$ than $D_1$, which would imply $$F_E\cdot D_1\geq \mu(L_2)D_1\cdot L_2\geq 4.$$ This shows that $F_E=[1,(k),1]$ for some $k\geq 0$, which contradicts $K\cdot\E\neq 0$.

Let $\alpha\geq 1$ be the number of connected components of $D_v$. We can assume that $L_1$ intersects $\E$ and $D_v$, because $F_E$ is connected. In particular $\mu(L_1)\geq 2$. Note that every vertical $(-1)$-curve intersects at most two other vertical components, hence each $L_i$ meeting $\E$ intersects $D_h$, otherwise it would be simple. Moreover, if such $L_i$ does not intersect $D_v$, which happens for example if $\mu(L_i)=1$, then $\#L_i\cap D_h\geq 2$. We consider two cases.

Suppose $L_i\cdot \E=0$ for $i\neq 1$, i.e. $L_1$ is the only $S_0$-component intersecting $\E$. Consider the contraction of $(-1)$-curves in $F_E$ different than $L_1$ (if there are any) until $L_1$ is the unique exceptional component in the image $F_E'$ of the fiber. This contraction does not touch $\E+L_1$, so $\E$ is one of the connected components of $\un F_E'-L_1$. Since $L_1\cdot D_h>0$, we have $\mu(L_1)\leq 3$, otherwise $D_h$ would have to contain an $n$-section for some $n>3$. It follows that either $F_E'=[2,1,2]$ or $F_E'=[3,1,2,2]$, hence $\E=[3]$. We have also $\mu(L_1)=3$, so $D_h$ contains a 3-section, which implies $F_\8=[3,1,2,2]$.

Now suppose $\E$ intersects more than one $L_i$, say $L_2\cdot \E>0$. We have $$5\geq F_E\cdot D_h\geq (D_v+\mu(L_1)L_1+\mu(L_2)L_2)\cdot D_h$$ and $\mu(L_2)L_2\cdot D_h\geq 2$, so $\alpha+\mu(L_1)L_1\cdot D_h\leq 3$, hence $\alpha=1$ and $\mu(L_1)=2$. This gives $F_E\cdot D=5$, so $F_\8=[3,1,2,2]$ and $D$ contains three horizontal components. In particular, no maximal twig of $D$ is contained in $F_\8$. We have now $L_2\cdot D_v=0$ and $\#L_2\cap D\geq 2$, so $\mu(L_2)=1$. Moreover, there are no more $(-1)$-curves in $F_E$. Defining $F_E'$ as the fiber $F_E$ with $L_1$ contracted we find that $F_E'$ has at most two $(-1)$-curves and they are of multiplicity one. Hence all components of $F_E'$ have multiplicity one, so $F_E'=[1,(k),1]$ for some $k\geq 0$. It follows that $F_E=[1,(k-1),3,1,2]$, hence $E=[3]$ and we are done. \end{proof}

(2) \begin{minipage}[t]{0.93\textwidth} If $\#E=1$ then $(B,T_1,T_2,T_3,\E)=([1],[(5)],[3],[2,2,3],[3])$ and $\ovk(W)=-\8$. \end{minipage}

\begin{proof} Suppose $\#E=1$ (and $\ovk(W)$ any). By \ref{2cor:S-Delta pre-minimally ruled} there exists a pre-minimal affine ruling of $S\setminus\Delta$, let $f$ be its extension as in \ref{2no:premin ruling notation}. We use the notation of \ref{2no:premin ruling notation}. In general $f$ need not be defined on $\ov S$, but at least the components of $\un F-Z_1-Z_u$ are not touched by $\varphi_f$ ($F$ is the fiber of $f$, not of $p$). In particular, the divisor of $D$-components of the second branch of $F$ and $Z_l$ are maximal twigs of $D$, denote them by $T_1$ and $T_2$ respectively. The unique $(-1)$-curve $C$ contained in $F$ is not touched by $\varphi_f$, so it is exceptional on $\ov S$ and satisfies $C\cdot D=1$, $C\cdot B=0$ and, since it is not simple, $\# C\cap \E\geq 2$. Now let us look at how $C$ behaves with respect to $p$. Fibers of $p$ cannot contain loops, so since $\E$ is connected and vertical for $p$, $C$ is horizontal for $p$ and $F_\8\cdot C=F_E\cdot C\geq 2$. We have $C\cdot D=1$, so $C$ intersects $F_\8-B$ in a component $D_0\subseteq T_1$ of multiplicity greater than one, hence $F_\8=[3,1,2,2]$, $D_0\cdot B=1$ and $D_0^2=-2$. In particular, we may assume that $D$ does not contain $[2,1,2]$.

We now look back at the fiber $F$ of $f$ and we find that since $D_0^2=-2$, $\Delta'=0$ and $T_1$ consists of $(-2)$-curves. Note that if $f$ is almost minimal then applying the above argument to $\wt C$ instead of $C$ we get that $\wt C$ intersects $D_0$, which contradicts the fact that $C$ and $\wt C$ intersect different maximal twigs of $D$. Thus $f$ is not almost minimal. Contraction of $T_1+C$ touches $Z_1$ precisely $x=\#T_1$ times, so $Z_1^2=-x-1$, hence $\varphi_f$ touches $Z_1$ precisely $x$ times, because $b=1$. We have $\wt Z_{lu}^2=1-x$. The proper transform of $\wt Z_{lu}$ on $\ov S$ is not a $(-2)$-curve, otherwise $D$ would contain the chain $[2,1,2]$, which was already ruled out. Therefore by \ref{2prop:pre-minimal ruling exists}(ii) we get $x\geq 5$ and $\Delta=0$.

Note that at least one of $T_2$, $T_3$, contains a $(-2)$-tip, otherwise we get a contradiction as in \ref{2lem:premin and short D elimination}. We check now that this implies $\ovk(W)=-\8$ and $\E=[3]$. Indeed, if $\ovk(W)\geq 0$ then by \ref{3lem:$W$ properties} $\wt e +\delta\leq 2$ and $\delta+\frac{1}{d(\E)}\geq 1$, so if, say, $T_2$ contains a $(-2)$-tip then $d_2\geq 5$ and we get $\frac{1}{d_1}+\frac{1}{d(\E)}\geq 1-\frac{1}{6}-\frac{1}{5}=\frac{19}{30}$, hence $d_1=d(\E)=3$. But then $T_2=[2,3]$ and $T_3=[3]$, so $\wt e +\delta=1+\frac{3}{5}+\frac{2}{3}>2$, a contradiction. By \ref{3lem:$W$ properties}(v) we infer that $\epsilon=2$, hence $\E=[3]$.

Suppose $\wt e_2+\frac{1}{d_2}>\frac{1}{2}$ and write $T_2^t=[c]+R$. We have $c\geq 3$, because $D$ does not contain $[2,1,2]$. The inequality gives $c d(R)-d'(R)\leq 2d(R)+1$, hence $(c-2)d(R)\leq d'(R)+1\leq d(R)$. Thus $c=3$ and $e(R)+\frac{1}{d(R)}=1$, so $R=[(y)]$ for some $y\geq 0$ by \ref{3lem:eT+x/d=1 solutions}. We have now $Z_l=T_2=[(y),3]$, so $Z_u=[2]$, $G=[y+2]$ and, since $f$ is pre-minimal, $\wt G+\wt Z_u=[(y),4,(x-3)]$ and hence $\wt Z_l=[y+2,2,x-1]$. We get $T_3=[y+2,2,x-2]$, and the inequality $\wt e+\delta>2$ reduces now to $x(3+5y+2y^2)<9y^2+27y+20$. Since $x\geq 5$, we get $(x,y)\in\{(6,0),(5,3),(5,2),(5,1),(5,0)\}$. By \ref{2cor:two fibers equations}(iii) $-\frac{1}{3}d(D)$ should be a square, which happens only for $(x,y)=(5,0)$, i.e. in the case listed above.

Thus we can assume $\wt e_2+\frac{1}{d_2}<\frac{1}{2}$. Since $\ovk(W)=-\8$, by \ref{3lem:$W$ properties}(iv) we get $\wt e_3+\frac{1}{d_3}>\frac{1}{2}$. As before, this is possible only if $T_3=[(y),3]$ for some $y\geq 0$. It follows that $\wt Z_l=[(y),4]$, because $\wt Z_{lu}$ is touched once by $\varphi_f$. Then $\wt Z_u=[2,2]$ and $\wt G=[y+2]$, so since the ruling is pre-minimal,  $G+Z_u=[(y)]$ and hence $T_2=Z_l=[y+1]$. Now $Z_1=[x+1]$ and $Z_1$, which is a proper transform of $B$, is touched $5$ times by $\varphi_f$, so $x=5$. Now the inequality $\wt e+\delta>2$ yields $y\leq 3$. We check that $-\frac{1}{3}d(D)$ is a square only for $y=2$, which again gives the case listed above. \end{proof}

We are therefore left with the case $(B,T_1,T_2,T_3,\E) = ([1], [(5)], [3], [2,2,3], [3])$. To exclude it we look more closely at the ruling $p$ induced by $F_\8=[3,1,2,2]$ contained in $D$ (the case is quite difficult to rule out, as one can check that all the equalities and inequalities derived so far in this paper are satisfied). We use the notation from (1). In fact there are two different chains $[3,1,2,2]$ contained in $D$, we consider the one not containing $T_2$. We have therefore $F_E\cdot D=5$. By (1) we know that $F_E=[1,3,1,2]$ or $[3,1,2,2]$ ($F_E'=F_E$ because $D_v$ consists of $(-2)$-curves), but in the second case the $1$-section contained in $T_3$ would have to intersect $L_1$, which is impossible, as $\mu(L_1)=3$. Thus $F_E=[1,3,1,2]$ and, as above, we denote the $(-1)$-curve intersecting $D_v$ by $L_1$ and the second one by $L_2$. Let $D'$ denote the divisor of vertical components of $D$ not contained in $F_\8\cup F_E$. Clearly, $D'=[2,2]\subseteq T_1$. Let $F'$ be the singular fiber containing $D'$. Since $F'$, which satisfies $d(F')=0$, consists of $D'$ and some number of $(-1)$-curves, we necessarily have $F'=[1,2,2,1]$. Denote the $(-1)$-curves of $F'$ by $M_1$, $M_2$, where $M_1$ intersects $T_3$. A fiber of $p$ other than $F_\8$, $F_E$ and $F'$ consists only of $S_0$-components, hence is smooth, because $\Sigma_{S_0}=2$. Let $\zeta\:\ov S\to \wt S$ be the contraction of $$B+F_\8\cap T_1+M_2+F'\cap T_1+L_2+L_1+T_3\cap F_\8+T_3',$$ where $T_3'$ is the section contained in $T_3$. Since the contracted divisor consists of disjoint chains of type $[1,(t)]$, $\wt S$ is smooth, hence $\wt S=\PP^2$. As $\mu(L_1)=2$, we have $T_2\cdot L_1=1$, so $T_2\cdot L_2=1$. The contractions of $B+F_\8\cap T_1$, $L_2+L_1+T_3\cap F_\8+T_3'$ and $M_2+F'\cap T_1$ touch $T_2$ respectively $3$, $4$ and $3(T_2\cdot M_2)^2$ times. The curve $\zeta(T_2)$ has degree $3$, which yields $T_2^2+3+4+3(T_2\cdot M_2)^2=9$, so $3(T_2\cdot M_2)^2=5$, a contradiction. \end{proof}

\blem\label{3lem:if k(W)<=0 then T1=[2]} If $\ovk(W)\leq 0$ then $\epsilon=2$ and one of the maximal twigs of $D$ equals $[2]$. \elem

\begin{proof} By \ref{3lem:$W$ properties}(iv) $b=1$. By \ref{3prop:no small 0-chains in D} $D$ does not contain $[2,1,2]$ or $[3,1,2,2]$ and by \ref{3lem:$W$ properties}(iv) we have $\wt e+\delta \geq 2$. We explore intensively these facts. Note that $\wt e_i+\frac{1}{d_i}\leq 1$ for each $i$. Assume that $d_1\leq d_2\leq d_3$ and write $T_i=[...,t_i',t_i]$ with $t_i'=\emptyset$ if $\#T_i=1$. Recall that by our convention the last component of $T_i$, the one with self-intersection $t_i$, intersects $B$. We prove successive statements.

\medskip (1) $T_1=[3]$ or $t_1=2$.

\begin{proof} Suppose $t_1=3$. Then $(t_2',t_2), (t_3',t_3)\neq (2,2)$ by \ref{3prop:no small 0-chains in D} and if $t_2=2$ (or $t_3=2$) then $t_3\neq 2$ ($t_2\neq 2$), so using \ref{3lem:bounds for e(F)} we get $\wt e_1<\frac{1}{2}$, $\wt e_2+\wt e_3<\frac{2}{3}+\frac{1}{2}$, hence $\wt e<\frac{5}{3}$. We use continuously this type of argument below having in mind \ref{3prop:no small 0-chains in D} and the inequality $\wt e+\delta \geq 2$. Suppose $t_1\geq 4$. If some other $t_i$ equals $3$ then $\wt e<\frac{1}{3}+\frac{1}{2} +\frac{2}{3}=\frac{3}{2}$ and if not then $\wt e<\frac{1}{3}+\frac{1}{3}+1=\frac{5}{3}$. Thus in any case $t_1\neq 2$ implies $\frac{3}{d_1}\geq\delta\geq2- \wt e>2-\frac{5}{3}=\frac{1}{3}$, so $d_1\leq8$. By \ref{1lem:SmallChains} we have to consider the following possibilities for $T_1$: $[4]$, $[5]$, $[6]$, $[7]$, $[8]$, $[2,3]$, $[2,4]$, $[2,2,3]$, $[3,3]$.

\bca $T_1$ is one of $[2,4]$, $[5]$, $[6]$, $[7]$ or $[8]$. In each case $\wt e_1+\frac{1}{d_1}\leq\frac{3}{7}$. If $(t_3',t_3)=(2,2)$ (or similarly $(t_2',t_2)=(2,2)$) then $\wt e_2<\frac{1}{3}$ and we get $\frac{1}{d_2}>2-\frac{3}{7}-1-\frac{1}{3}$, so $d_2\leq 4$, a contradiction with $d_2\geq d_1$. In other case $\wt e+\frac{1}{d_1}<\frac{3}{7}+\frac{2}{3}+\frac{1}{2}$, so $\frac{2}{d_2} \geq \frac{1}{d_2}+\frac{1}{d_3} \geq 2-\wt e-\frac{1}{d_1}>\frac{17}{42}$ and again $d_2\leq 4$, a contradiction.\eca

\bca $T_1$ is $[2,2,3]$ or $[3,3]$. Then $\wt e_1+\frac{1}{d_1}\leq\frac{4}{7}$ and $\wt e_2+\wt e_3< \frac{1}{2}+ \frac{2}{3}$, so $\frac{2}{d_2}\geq2-\wt e-\frac{1}{d_1}>\frac{1}{4}$ and $d_2\leq7$. Since $d_1\leq d_2$ we get $T_1=[2,2,3]$ and $d_1=d_2=7$. By renaming $T_1$ and $T_2$ we can assume that $t_2\neq 2$. In fact we can assume that $T_2=[2,2,3]$ because other cases ($[7]$ and $[2,4]$) were excluded above. Thus $\wt e_3+\frac{1}{d_3}\geq \frac{6}{7}$. We have $\wt e_3<\frac{2}{3}$, because $(t_3',t_3)\neq (2,2)$, so $\frac{1}{d_3}>\frac{6}{7} -\frac{2}{3}$ and then $d_3\leq 5<d_1$, a contradiction.\eca

\bca $T_1=[4]$. We have $\wt e_1+\frac{1}{d_1}=\frac{1}{2}$, so $\frac{1}{d_2}+\frac{1}{d_3}\geq\frac{3}{2}-\wt e_2-\wt e_3$. We have $t_2+t_3\geq 5$. If $t_2\geq 4$ (or similarly $t_3\geq 4$) then $\frac{1}{d_2}\geq\frac{3}{2}-\wt e_2-1 >\frac{1}{6}$, so $d_2\leq 5$. If $t_2=3$ (or similarly $t_3=3$) then $\frac{2}{d_2}>\frac{3}{2}-\frac{2}{3}- \frac{1}{2}= \frac{1}{3}$, so again $d_2\leq 5$. Note that since $\wt e_3+\frac{1}{d_3}\leq 1$, $\wt e_2+\frac{1}{d_2}\geq \frac{1}{2}$, so $T_2\neq [5]$ (and similarly $T_3\neq [5]$). If $T_2$ is one of $[2,3]$, $[3,2]$ or $[2,2,2,2]$ then we have respectively $\wt e_2+\frac{1}{d_2}=\frac{3}{5}, \frac{4}{5}, 1$ and using \ref{3prop:no small 0-chains in D} we bound $\wt e_3$ from above respectively by $\frac{2}{3}, \frac{1}{2}$ and $\frac{1}{3}$, which gives $d_3=5$. However, we check easily that for $d_2=d_3=5$ the inequality $\frac{1}{d_2}+\wt e_2+\frac{1}{d_3}+\wt e_3\geq\frac{3}{2}$ cannot be satisfied. Thus $d_2=4$. By renaming $T_1$ and $T_2$ we can assume that $T_2\neq [2,2,2]$, so $T_2=[4]$. Then $\wt e_3+\frac{1}{d_3}\geq1$ so $T_3=[2,2,2]$ by \ref{3lem:eT+x/d=1 solutions} and after renaming $T_1$ and $T_3$ we are done.\eca

\bca $T_1=[2,3]$. We have $\wt e_2+\wt e_3+\frac{1}{d_2}+\frac{1}{d_3} \geq\frac{7}{5}$ and $\wt e_2+\wt e_3<\frac{2}{3}+ \frac{1}{2}$, so $d_2\leq8$. Suppose $d_2=5$. We can assume that $T_2=[2,3]$, because the case $T_1=[5]$, $T_2=[2,3]$ was considered above and in other cases $t_2=2$, so after renaming $T_1$ and $T_2$ we are done. If $d_3\neq 5$ then $\wt e_3\geq \frac{4}{5}- \frac{1}{d_3}>\frac{3}{5}$, hence $(t_3',t_3)=(2,2)$, a contradiction. Therefore $d_3=5$ and again we can assume that $T_3=[2,3]$, so $\wt e_2+\wt e_3+\frac{1}{d_2}+\frac{1}{d_3}=\frac{6}{5}$, a contradiction. Thus $6\leq d_2\leq 8$. If $T_2=[d_2]$ then $\frac{1}{d_3}+\wt e_3>\frac{7}{5}- \frac{2}{5}=1$, a contradiction. It follows that $T_2$ is one of $[2,2,3]$, $[2,4]$, $[3,3]$, $[4,2]$ or $[2,3,2]$ (in particular $d_2>6$). By \ref{3prop:no small 0-chains in D} $\wt e_3<\frac{2}{3}$ in first three cases and $\wt e_3<\frac{1}{2}$ in the latter two cases. In each case we obtain $\wt e_3+\wt e_2+\frac{1}{d_2}\leq \frac{5}{4}$, hence $d_3\leq 6<d_2$, a contradiction. \qedhere \eca  \end{proof}

(2) $T_1=[3]$ or $T_1=[2]$.

\begin{proof} Suppose $\#T_1\neq 1$. We have $\wt e_2+\wt e_3+\frac{1}{d_2}+\frac{1}{d_3}\geq 1$. By (1) $t_1=2$, so $t_2,t_3\neq 2$, hence $\wt e_2+\wt e_3<\frac{1}{2}+\frac{1}{2}=1$ and from the inequality $\wt e+\delta\geq 2$ we get $\wt e_1+\frac{3}{d_1}>1$. This gives $d'(T_1^t)=d(T_1^t)-1$ or $d'(T_1^t)=d(T_1^t)-2$, so $T_1=[(k)]$ or $[3,(k)]$ for some $k>0$ by \ref{3lem:eT+x/d=1 solutions}.

Suppose $k\geq 2$. In this case $t_2,t_3\geq 4$, so $\wt e_2,\wt e_3<\frac{1}{3}$. Then $\frac{1}{d_2}+\frac{1}{d_3}>\frac{1}{3}$ and we get $d_1\leq d_2\leq 5$, which is possible only if $T_2$ is a tip and $T_1=[(k)]$ for some $k\in\{2,3,4\}$. Since now $\frac{1}{d_3}\geq 1-\wt e_3-\frac{2}{d_2}> \frac{2}{3}-\frac{1}{2}$, we see that $d_3\leq 5$, so $T_3$ is also a tip. Then $\wt e_2=\frac{1}{d_2}$ and $\wt e_3=\frac{1}{d_3}$, so $\frac{1}{d_2}+ \frac{1}{d_3} \geq \frac{1}{2}$ and we conclude that $T_2=T_3=[4]$ and $T_1=[(k)]$ for some $k\in\{2,3\}$. It follows that $\wt e+\delta=2$, so $\ovk(W)=0$ and by \ref{3lem:$W$ properties} $\frac{1}{k+1}+\frac{1}{|G|}\geq \frac{1}{2}$. Then $|G|\leq 6$, so $G$ is abelian, because it is a small subgroup of $GL(2,\C)$. However, by \ref{1lem:eps properties}(iii) $\#\E=7+K\cdot E+\epsilon-k\geq 7$, a contradiction.

We are left with the case $T_1=[3,2]$, for which $\wt e_2+\frac{1}{d_2}+\wt e_3+\frac{1}{d_3}\geq \frac{6}{5}$. Now $t_2,t_3\neq 2$, so $\wt e_2,\wt e_3<\frac{1}{2}$. Suppose $t_2\geq 4$ or $t_3\geq 4$. Then $\wt e_2+\wt e_3<\frac{1}{2}+\frac{1}{3}$, so $\frac{1}{d_1}+\frac{1}{d_2}>\frac{1}{3}$ and we get $d_2=5$, hence $T_2=[5]$ or $T_2=[2,3]$. If $T_2=[5]$ then $\frac{1}{d_3}>\frac{4}{5}-\frac{1}{2}=\frac{3}{10}$. If $T_2=[2,3]$ then, since $t_3\geq 4$, $\wt e_3<\frac{1}{3}$ and $\frac{1}{d_3}> \frac{3}{5}-\frac{1}{3}=\frac{4}{15}$. In both cases we get $d_2\leq 3$, a contradiction. Thus $t_2=t_3=3$, so $\wt e_2+\wt e_3<1$ and we get $d_2\leq 9$. However, all admissible chains with discriminant $5\leq d\leq 9$ which end with a $(-3)$-curve satisfy $\wt e+\frac{1}{d}\leq\frac{3}{5}$ (cf. \ref{1lem:SmallChains}), the equality occurs only for $[2,3]$. Hence $\frac{1}{d_3}\geq\frac{3}{5}-\wt e_3>\frac{1}{10}$, so $d_3\leq 9$ too. This implies $T_2=T_3=[2,3]$, so $\wt e+\delta=2$, which gives $\ovk(W)=0$. By \ref{3lem:$W$ properties}(iii) $\frac{1}{|G|}\geq \frac{2}{5}$, a contradiction. \end{proof}

(3) $T_1=[2]$.

\begin{proof} Suppose $T_1=[3]$. We have $\wt e_2+\wt e_3+\frac{1}{d_2}+ \frac{1}{d_3}\geq\frac{4}{3}$, so since $\wt e_2+\wt e_3<\frac{2}{3}+\frac{1}{2}$, we get $\frac{1}{d_1}+\frac{1}{d_2}>\frac{1}{6}$, which gives $d_2\leq 11$.

\setcounter{case}{0}\bca Suppose $T_2\neq [3]$ or $(t_3',t_3)\neq (3,2)$. We prove that $d_3\leq 42$. For $d_2>6$ the inequality $\frac{1}{d_1}+\frac{1}{d_2}>\frac{1}{6}$ gives $d_3\leq 42$. We can therefore assume that $d_2\leq 6$. If $T_2=[3,2]$ then $\wt e_2+\frac{1}{d_2}=\frac{4}{5}$ and $t_3\neq 2$, so $\frac{1}{d_3}> \frac{4}{3}-\frac{4}{5}-\frac{1}{2}$ and $d_3\leq 29$. If $T_2=[4]$, $[5]$, $[6]$ or $[2,3]$ then $\wt e_2+\frac{1}{d_2}\leq\frac{3}{5}$ and since $\wt e_3< \frac{2}{3}$, we get $d_3\leq 14$. We are left with the case $T_2=[3]$, where we get $\wt e_3+ \frac{1}{d_3}\geq\frac{2}{3}$. If $t_3\geq 3$ then $\frac{1}{d_3}>\frac{2}{3}- \frac{1}{2}$, so $d_3\leq 5$. If $t_3=2$ and $t_2>3$ then $\frac{1}{d_3}>\frac{2}{3}-\frac{3}{5}$, so $d_3\leq 14$ and we are done.

Now note that whenever $d_3$ is bounded, by \ref{2rem:finite number of cases} there are finitely many possibilities for the weighted dual graphs of $D$ and $\E$. Using a computer program we checked that the conditions $d_2\leq 11$, $d_3\leq 42$, \ref{1lem:eps properties}(iii)-(iv), \ref{1lem:if eps<2}, \ref{1prop:possible $E$}, \ref{2lem:Zar and others for D fork} and \ref{1prop:basic S' properties}(iv) (which implies that $-d(D)/d(\E)$ is a square) are satisfied only in two cases:\benum[(i)]

\item $T_1=[3]$, $T_2=[3]$, $T_3=[3,(6)]$ and $\E=[2,3,4]$,

\item $T_1=[3]$, $T_2=[4]$, $T_3=[2,2,2]$ and $\E$ is a fork with a $(-2)$-curve as a branching component and maximal twigs $[2]$, $[2]$, $[2,2,3]$. \eenum

\noin In both cases $D$ contains $[3,1,2,2]$, a contradiction.\eca

\bca Suppose $T_2=[3]$ and $(t_3',t_3)=(3,2)$, write $T_3=T_0+[3,2]$. Using \ref{1lem:d() expansion} we check that the inequality $\wt e+\frac{1}{d_3}\geq\frac{2}{3}$ is equivalent to $d'(T_0^t)+3\geq d(T_0^t)$, so by \ref{3lem:eT+x/d=1 solutions} $T_3=[(k),3,2]$, $[3,(k),3,2]$, $[4,(k),3,2]$ or $[2,3,(k),3,2]$ for some $k\geq 0$. We conclude that $K\cdot T\leq 5$, hence \ref{2rem:finite number of cases} again reduces the problem to checking finitely many cases (here Noether formula implies $k\leq 9$, which gives $d_3\leq 102$). We checked that each of them leads to a contradiction with one of the conditions as in Case 1.
\qedhere \eca \end{proof}

It remains to prove that $\epsilon=2$. By (3) and \ref{3lem:$W$ properties}(v) we can assume $\ovk(W)=0$. For convenience we put formally $[3,(-1),3]=[4]$, then we have $d([3,(k-2),3])=4k$ for any $k\geq 1$. Suppose $\epsilon\leq 1$. By \ref{3lem:$W$ properties}(v) $2(K_{\ov S}+T+\E)\geq 0$, so by \ref{1lem:Eff-NegDef is Eff}(ii) $[2(K_{\ov S}+T^\#+\E^\#)]\sim U$ for some effective $U$. Then $K_{\ov S}+T^\#+\E^\#\equiv 0$ implies $U+\{2(K_{\ov S}+T^\#+\E^\#)\}\equiv0$, hence $2\Bk T_i$ and $2\Bk \E$ are $\Z$-divisors. Since $T_2$, $T_3$, $\E$ do not consist only of $(-2)$-curves, we obtain $2\Bk \E=\E$ and $2\Bk T_i=T_i$ for $i=2,3$. The latter equality holds only if $T_2$ and $T_3$ are of type $[3,(k),3]$ for some $k\geq -1$. Using \ref{2lem:Zar and others for D fork}(iv) we compute $\Bk^2\E=-\epsilon$, hence by \ref{1rem:-2 chains and forks} and \ref{1lem:ForkBk}(v) $\epsilon=1$ and $\E$ is a chain. Then we can write $\E=[3,(z-2),3]$ with $z\geq 1$. By \ref{1lem:eps properties}(iii) $x+y+z=11$, hence $1\leq x,y\leq 9$ and $$\frac{1}{x}+\frac{1}{y}+\frac{1}{11-x-y}\geq 2$$ by \ref{3lem:$W$ properties}(iii). This inequality is satisfied only for $(x,y)=(1,1)$ and $(x,y)=(1,9)$. However, in the first case $d(D)=0$, so $(x,y)=(1,9)$ and we get $T_2=[4]$, $T_3=[3,(7),3]$ and $\E=[4]$. By \ref{2lem:premin and short D elimination} there exists an affine ruling of $S$ extending to a $\PP^1$-ruling of $\ov S$. Since $B^2=-1$, $B$ is horizontal and the ruling has three singular fibers. This contradicts \ref{2lem:three fibers only for b=2}. \end{proof}

\bprop\label{3prop:k(W)=2} $\ovk(W)=2$.\eprop

\begin{proof}

Suppose $\ovk(W)\leq 1$. By \ref{3lem:$W$ properties} $\ovk(W)\leq 0$ and $b=1$. By \ref{3lem:if k(W)<=0 then T1=[2]} one of the maximal twigs of $D$ is $[2]$. We have also $\epsilon=2$, which gives $E=[3]$. Denote the coefficient of $E$ in $\Bk \E$ by $w_E$. We prove successive statements.

\med (1) \begin{minipage}[t]{0.93\textwidth}  If $w_E>\frac{1}{2}$ then $\E$ is a chain and $\Delta$ is connected. If $w_E=\frac{1}{2}$ then either $\E$ is a fork with maximal twigs $[3]$, $[2]$, $[2]$ or $\E=[2,3,2]$. \end{minipage}

\begin{proof} Suppose $\E$ is a fork. By \ref{1prop:if eps=0 or E^ fork}(iii) we know that $\Delta$ does not contain a fork and by \ref{2cor:Delta with three components} $E$ is not the branching component of $\E$, so $\E$ is of type (b1) (cf. \ref{1prop:possible $E$}) and the maximal twig of $\E$ containing $E$ is equal to $[(k),3]$ for some $k\geq 0$. Using \ref{1lem:ChainBk}(ii) and the definition of a bark of an admissible fork it is a straight computation to check that $w_E\leq\frac{1}{2}$ in each case and the equality occurs only for a fork with maximal twigs $[3]$, $[2]$, $[2]$. If $\E$ is a chain then $\E=[(m-1),3,(\wt m-1)]$ for some $m,\wt m\geq 1$ and $$w_E=\frac{m+\wt m}{m\wt m+m+\wt m}=1-\frac{1}{1+\frac{1}{m}+\frac{1}{\wt m}},$$ so $w_E\geq \frac{1}{2}$ if and only if $\frac{1}{m}+\frac{1}{\wt m}\geq 1$, hence (1) follows.\end{proof}

By \ref{2cor:S-Delta pre-minimally ruled} there exists a pre-minimal affine ruling of $S\setminus\Delta$, let $f:(\ov S^\dagger,D^\dagger+\Delta)\to \PP^1$ be its minimal completion. Since $\Sigma_{S_0}=0$, every singular fiber of $f$ has a unique $S_0$-component and this component is a $(-1)$-curve. We use the notation \ref{2no:premin ruling notation}. Since $b=1$ and $Z_1^2\leq -2$, $n=1$ and by \ref{2cor:two fibers equations} $(\wt h,h)=(2,3)$. Write $\Delta'=[(m-1)]$, $\wt \Delta=[(\wt m-1)]$ for some $m,\wt m\geq 1$. The maximal twig of $D^\dagger$ contained in the first branch of $F$, call it $T_2$, and the one contained in the second branch of $F$, call it $T_1$, are not touched by $\varphi_f$, hence they are maximal twigs of $D$.

Fibers of $\PP^1$-rulings cannot contain branching $(-1)$-curves, so since $b=1$, $\varphi_f$ touches the birational transform of $B$. Let $\ov S^\dagger\to\wt S\xra{\wt \rho} \ov S$ be the factorization of $\varphi_f$, such that the birational transform of $B$ is touched by $\wt \rho$ exactly once. Let $\wt \pi\:\wt S\to \wt U$ and $\pi\:\ov S\to U$ be the contractions of $T_1+C+\Delta'$ on respective surfaces.

$$\xymatrix{\ov S^\dagger\ar[r] & \wt S\ar[r]^{\wt \pi}\ar[d]^{\wt \rho} & \wt  U\ar[r]^{\eta}\ar[d]^{\rho}& \PP^1 \\   &\ov S\ar[r]^{\pi} &U &}$$

The centers of $\wt \rho$ and $\wt \pi$ are different, so there exists a birational morphism $\rho\:\wt U\to U$, such that $\rho \circ \wt \pi=\pi\circ \wt \rho$. Denote the birational transform of $B$ contained in $\wt U$ by $\wt B$. By definition $\wt B^2=0$. Consider the $\PP^1$-ruling $\eta:\wt U\to \PP^1$ induced by $\wt B$. Denote by $\wt T_3,\wt E\subseteq \wt U$ the reduced total inverse image of $T_3$ and the birational transform of $E$ respectively. Put $\wt D=T_2+\wt B+\wt T_3$. Let $D_2\subseteq T_2$ and $D_3\subseteq \wt T_3$ be the sections of $\eta$ contained in $\wt D$ and let $F'$ be a general fiber. Since $\Sigma_{S_0}=1$ for the ruling $\eta\circ\wt \pi$, there exists a unique singular fiber $F_1$ with $\sigma(F_1)=2$. Let $M_1$, $M_2$ be its $S_0$-components.

\med (2) \begin{minipage}[t]{0.95\textwidth} $M_1$ and $M_2$ are $(-1)$-curves. If $\eta$ has more than one singular fiber then $F_1=M_1+\wt \Delta+M_2$.  \end{minipage}

\begin{proof} Suppose there is another singular fiber $F_0$. Note that vertical $(-1)$-curves are $S_0$-components. We have $\sigma(F_0)=1$, so $F_0$ is a chain intersected in tips by $D_2,D_3$, otherwise there would be a loop in $\Supp D$. Then $F_0$ contains $T_3-D_2+T_2-D_2$, so $F_1$ does not contain $\wt D$-components. Since $M_i\cdot D=M_i\cdot (D_2+D_3)$, both $M_i$ intersect $D_2+D_3$, hence both have multiplicity one. It follows that $F_1=[1,(\wt m-1),1]$, so we are done. We can therefore assume that $F_1$ is the unique singular fiber of $\eta$. Suppose $F_1$ has only one $(-1)$-curve. Then $D_2$ and $D_3$ intersect tips of $F_1$ belonging to the first branch of $F_1$, so when we contract $F_1$ to a smooth fiber we touch $D_2+D_3$ at most once. This gives two disjoint sections of a $\PP^1$-ruling of a Hirzebruch surface, one negative and one non-positive, which is a contradiction.   \end{proof}

The morphism $\wt \pi$ contracts the fiber consisting of $T_1+C+\Delta'$, so since $h=3$, we can write $$\wt \pi=p_{2}\circ\sigma_2\circ p_1\circ \sigma_1,$$ where $p_1,p_2$ are sprouting blow-ups (with respect to the image of the fiber) and $\sigma_i$ are compositions of sequences of subdivisional blow-downs. Note that $p_1\circ\sigma_1$ is the contraction of $C+\Delta'$. Put $\sigma=\sigma_2\circ p_1\circ \sigma_1$ and let $R_i$ for $i=1,2$ be the exceptional divisors of $p_i$. We now analyze the contraction $\wt \pi$ and singular fibers of $\eta$ more closely.

\med (3) $\wt E\cdot(K_{\wt U}+\wt D)+E\cdot \sigma^*R_2=1$.

\begin{proof}  Let us use the common letter $E'$ for the birational transforms of $E$. Using \ref{1lem:blow-up and transforms} we check how the quantity $E'\cdot(K'+D')$, where $D'$ is the reduced total transform of $\wt D$ and $K'$ the canonical divisor on a respective intermediate surface between $\wt S$ and $\wt U$, changes under subsequent blow-downs. Since $\wt \rho$ is subdivisional with respect to $D$, at the beginning we have $$E'\cdot(K'+D')=E\cdot(K+D+C+\Delta')=1+E\cdot (C+\Delta').$$ Under $\sigma$ it decreases by $E'\cdot R_1=E\cdot \sigma_1^*R_1=E\cdot (C+\Delta')$ and under $p_2$ it decreases by $E'\cdot R_2=E\cdot \sigma^*R_2$. \end{proof}

(4) \begin{minipage}[t]{0.93\textwidth} There is a unique $(-1)$-curve $L$, such that $L\cdot \wt D>1$. It satisfies $K_{\wt U}+\wt D+L\equiv 0$. \end{minipage}

\begin{proof} We have $$K_{\wt U}\cdot(K_{\wt U}+\wt D)=K_U\cdot(K_U+\pi_*D)=K\cdot(K+D)+1=1,$$ so by Riemann-Roch's theorem $$h^0(-K_{\wt U}-\wt D)+h^0(2K_{\wt U}+\wt D)\geq K_{\wt U}\cdot(K_{\wt U}+\wt D).$$ If $2K_{\wt U}+\wt D\geq 0$ then $$0\leq \kp(K_{\wt U}+\wt D)= \kp(K_U+\pi_*D)= \kp(K+D+C+\Delta')=\kp(K+D),$$ where the last equality follows from \ref{1lem:Eff-NegDef is Eff}(i), and this contradicts $\kp(K+D)=\ovk(S)=-\8$. We get $-K_{\wt U}-\wt D\geq 0$. Write $-K_{\wt U}-\wt D=\sum C_i$ for some irreducible $C_i$'s, such that $C_i^2<0$ (cf. \ref{1lem:linear systems}(ii)). For a fiber $F'$ of $\eta$ we have $F'\cdot(K_{\wt U}+\wt D)=0$, so $C_i$'s are vertical.

Each $S_0$-component $L$ of a singular fiber intersects $\wt D$ and by (2) it is a $(-1)$-curve. Suppose each satisfies $L\cdot \wt D=1$. Then $F_1$ is the only singular fiber of $\eta$. Indeed, if $F'\neq F_1$ is a singular fiber then $\sigma(F')=1$ and since $\Supp \wt D$ does not contain a loop, $F'$ is a chain, so its exceptional component does not satisfy our assumption. $F_1\cap \wt D$ has two connected components (which may be points), let $R\subseteq M_1+\wt \Delta+M_2$ be a chain connecting them. By assumption $R\neq M_1,M_2$, so $R$ contains both $M_i$. It follows that $R$ contains a divisor with zero discriminant, which is possible only if $F_1=[1,(\wt m-1),1]$, hence $T_2=D_2$ and $T_3=D_3$. If we now look at the pre-minimal ruling of $S\setminus\Delta$ then we see that $\wt Z_l$ and $Z_l$ are irreducible, so $\wt G$ and $G$ are $(-2)$-curves, which implies that $D$ contains a component with non-negative self-intersection, a contradiction. Thus there is an exceptional $S_0$-component $L$, such that $L\cdot \wt D>1$.

Note that if for some $i\in\{2,3\}$ the section $D_i$ intersects $L$ then $D_i$ is a maximal twig of $\wt D$, because $D_i\cdot F=1$. It follows that $L\cdot \wt D=2$. Since $(-K_{\wt U}-\wt D)\cdot L=1-\wt D\cdot  L<0$, $L$ appears among $C_i$'s. However, $-K_{\wt U}-\wt D-L$ is vertical and satisfies $$(-K_{\wt U}-\wt D-L)^2=K_{\wt U}\cdot(K_{\wt U}+\wt D)-1=0$$ so $-K_{\wt U}-\wt D-L\equiv\alpha F$ for some $\alpha\geq 0$. Multiplying by $D_i$ for $i=2,3$ we get $\beta_{\wt D}(D_i)+L\cdot D_i=2-\alpha$. For $\alpha>0$ we would obtain $\beta_{\wt D}(D_2)=\beta_{\wt D}(D_3)=1$ and $L\cdot D_2=L\cdot D_3=0$, which is impossible, as $L\cdot \wt D>0$. Thus $K_{\wt U}+\wt D+L\equiv 0$. If $L'$ is another $(-1)$-curve, such that $L'\cdot \wt D>1$, then $-L'\cdot L=L'\cdot (K_{\wt U}+\wt D)>0$, hence $L'=L$. \end{proof}

(5) $2\leq E\cdot \sigma^*R_2=1+E\cdot L\leq 3$.

\begin{proof} Intersecting $K_{\wt U}+\wt D+L\equiv 0$ with components of $\wt D+\wt \Delta$ we see that $L\cdot \wt \Delta=0$ and $L$ intersects $\wt D$ only in tips, each tip once. It follows that $\rho$ and $\pi$ do not touch $L$. Intersecting $$K+T+\E\equiv \lambda\cal P+\Bk T+\Bk \E$$ with $L$ we get $$E\cdot L(1-w_E)\leq (\Bk T_2+\Bk T_3)\cdot L-1.$$ We have $(\Bk T_1+\Bk T_3)\cdot L<2$, otherwise $T_2$ and $T_3$ would be $(-2)$-chains, which is impossible by \ref{3prop:no small 0-chains in D}. Thus $E\cdot L<\frac{1}{1-w_E}$. By (3) we get $$E\cdot \sigma^*R_2=1-\wt E\cdot(K_{\wt U}+\wt D)=1+E\cdot L<1+\frac{1}{1-w_E}.$$ By (2) either $w_E\leq \frac{1}{2}$ or $\E=[3,(n-1)]$ for some $n\geq 1$ and then $\frac{1}{1-w_E}=2+\frac{1}{n}\leq 3$. In any case $E\cdot \sigma^*R_2\leq 3$.

Consider the ruling $\eta\circ\wt \pi:\wt S\to \PP^1$. Let $\mu_C$ and $\mu_{\Delta}$ be the coefficients in $\sigma^*R_2$ of $C$ and respectively of a component of $\Delta'$ intersecting $E$ (put $\mu_{\Delta}=0$ for $\Delta'=0$). Clearly, $\wt \rho$ does not touch $T_1+C+\Delta'+E$. We have $E\cdot \sigma^*R_2=\mu_CC\cdot E+\mu_\Delta$ and $\mu_\Delta<\mu_C$. Note that $E\cdot \sigma^*R_2\geq 2$, otherwise $E\cdot (C+\Delta')\leq 1$, a contradiction with \ref{1prop:basic S' properties}(ii). \end{proof}

(6) $T_1=[(k),3]$ for some $k\geq 1$. $\E=[3,2]$.

\begin{proof} Suppose first that $\#T_1=1$. Then $E\cdot \sigma^*R_2=E\cdot F'$ for a generic fiber $F'$ of $\eta\circ\wt \pi$. By (5) we have $$2\leq E\cdot L+1=E\cdot  F'=\mu_C C\cdot E+\mu_\Delta\leq 3.$$ Suppose $L\nsubseteq F_1$ (cf. (2)). The fiber containing $L$ has $\sigma=1$, so $\mu(L)\geq 2$ and since $\mu(L)E\cdot L\leq E\cdot F'\leq 3$, we get $E\cdot F'=E\cdot L+1=2$. Then $F_1=M_1+\wt \Delta+M_2$ by (2), because $L$ is contained in some singular fiber. Since both $M_i$ intersect $\wt D$, we have $\wt D\cdot M_1=\wt D\cdot M_2=1$. By \ref{1prop:basic S' properties}(ii) $\E\cdot M_1,\E\cdot M_2\geq 2$, so $\wt \Delta\neq 0$ and then $E\cdot \wt \Delta=E\cdot (F'-M_1-M_2)\leq 0$, a contradiction. Therefore $L\subseteq F_1$, say $L=M_1$. By (4) $\wt D\cdot M_2\leq 1$, so $\E\cdot M_2\geq 2$ by \ref{1prop:basic S' properties}(ii). We have $E\cdot M_2\leq E\cdot (F_1-L)=1$, so $0\neq \wt \Delta\subseteq F_1$ and $E\cdot M_2\leq E\cdot (F_1-L-\wt \Delta)\leq 0$. Then $\E\cdot M_2=\wt \Delta\cdot M_2\leq 1$, a contradiction. Thus $\#T_1>1$.

Suppose $\mu_\Delta=0$. Then $\Delta'=0$, so $C\cdot E\geq 2$. Since $\mu_C C\cdot E+\mu_\Delta\leq 3$, we get $\mu_C=1$, so $T_1=[(k)]$ for some $k\geq 0$. Since $\#T_1>1$, $D$ contains $[2,1,2]$ by \ref{3lem:if k(W)<=0 then T1=[2]}, a contradiction with \ref{3prop:no small 0-chains in D}. Thus $\mu_\Delta>0$. We get $\mu_C>1$ and then $\mu_C=2$, $\mu_\Delta=1$ and $C\cdot E=1$. As $\#T_1>1$, it follows that $T_1$ is $[(k),3]$ or $[3,(k)]$ for some $k\geq 1$. However, in the latter case the equality $h=3$ does not hold. Thus $T_1=[(k),3]$ for some $k\geq 1$. We conclude that $\Delta'=[2]$ and $E\cdot  \sigma^*R_2=3$, so $E\cdot L=2$. Since $E\cdot L<\frac{1}{1-w_E}$ (cf. (5)), we get $\wt \Delta=0$ by (1). \end{proof}

(7) $T_2=[2]$.

\begin{proof} Recall that $T_2$ is the maximal twig of $D$ contained in the first branch of $F$ (a fiber of $f$). Suppose $T_2\neq [2]$. By (6) and \ref{3lem:if k(W)<=0 then T1=[2]} $T_3=[2]$, so since $\#T_3=1$, $f$ is not almost minimal. Thus by \ref{2prop:pre-minimal ruling exists} the morphism $\varphi_f:\ov S^\dagger\to \ov S$ minimalizing $D^\dagger$ contracts precisely $H^\dagger+\wt Z_1$ and touches $Z_1$ at least four times. However, since $\wt \Delta=\emptyset$, $\wt G+\wt Z_u+\wt Z_1$ consists of $(-2)$-curves, hence $\varphi_f$ touches $Z_1$ at most once, a contradiction. \end{proof}

From (7) we see that $F$ is produced by the following sequence of characteristic pairs (cf. \ref{1def:char pairs} and \ref{2no:pairs notation}): $\binom{4k+4}{2k+2}$, $\binom{2k+2}{2}$, $\binom{2}{1}$, so the pairs $\binom{\un c_i}{\un p_i}$ are $\binom{2k+2}{k+1}$, $\binom{k+1}{1}$. By (6) $C\cdot E=1$ and $\kappa=2C\cdot E+1=3$. The second fiber $\wt F$ of $f$ is produced by the sequence $\binom{c}{p}$, $\binom{1}{1}$ for some $c,p\geq 1$. We have $$\wt \kappa c=d=\kappa \un c_1=6k+6.$$ By \eqref{eq:pairs-general I} $3d+1=\kappa (2k+2+k+1+1)+\wt \kappa(c+p)$, hence $\wt \kappa p=3k+1$. Then $$\wt \kappa=gcd(\wt \kappa c,\wt \kappa p)=gcd(6k+6,3k+1)=gcd(4,3k+1),$$ so $\wt \kappa\in\{2,4\}$ ($\wt C$ would be simple for $\wt \kappa=1$). On the other hand \eqref{eq:pairs-general II} gives $$d^2+3=\wt \kappa^2c p+\wt \kappa^2+9(2(k+1)^2+k+1)+3C\cdot E+C\cdot E+1,$$ hence $\wt \kappa^2=3k+1$. For $\wt \kappa=2$ we get $k=1$, so $(c,p)=(6,2)$, which contradicts the relative primeness of $c$ and $p$. Thus $\wt \kappa=4$ and we get $k=5$ and $(c,p)=(9,4)$. Then $\wt G+\wt Z_u=[3,2,2,2]$ and $\wt Z_l=[2,5]$, so $T_3=[2,4]$. Then $\wt e+\delta=\frac{3}{7}+1+\frac{7}{13}<1$, a contradiction with \ref{3lem:$W$ properties}(iv). \end{proof}

\bcor\label{3cor:final bounds for E^} $\E$ is one of $[2,3]$, $[3]$, $[4]$, $[5]$ and $\epsilon\in \{1,2\}$. \ecor

\begin{proof} By \ref{3prop:k(W)=2} $\ovk(W)=2$, so by \ref{3lem:$W$ properties}(iii) and \ref{2lem:Zar and others for D fork}(ii) we have $\epsilon\neq 0$ and $1>\delta>1-\frac{1}{|G|}$. Suppose $|G|\geq 7$ and assume $d_1\leq d_2\leq d_3$.  For $d_1\geq 3$ we get $d_2=3$ and $d_3\leq 5$. For $d_1=2$ we have $d_2\geq 3$ and the inequality gives $d_2\leq 5$ and $\frac{1}{d_3}> \frac{6}{7}-\frac{1}{2}-\frac{1}{3} =\frac{1}{42}$, so $d_3\leq 41$. By \ref{2rem:finite number of cases} there are only finitely many possibilities for the weighted dual graphs of $\E$ and $D$. Using a computer program we checked that with the above bounds conditions \ref{1lem:eps properties}(iii), \ref{1prop:possible $E$}, \ref{2lem:Zar and others for D fork} and \ref{1prop:basic S' properties}(iv) can be satisfied only for $\E=[4]$, which contradicts our assumption. We conclude that $|G|\leq 6$, so $\E$ is one of: $[2,3]$, $[3]$, $[4]$, $[5]$, $[6]$. However, $[6]$ is ruled out by \ref{1cor:possible KE and eps}. \end{proof}

\med\section{Special cases}\label{sec:special cases}

By section \ref{sec:surface W} we know that $\ovk(W)=2$ and $(\epsilon,\E)\in \{(2,[2,3]), (2,[3]), (1,[4]), (1,[5])\}$. We will rule out these cases now. Let $f\:(\ov S^\dagger,D^\dagger)\to \PP^1$ be a minimal completion of a pre-minimal affine ruling of $S\setminus\Delta$ (see Fig. \ref{fig:pre-minimal ruling}). We use the notation of \ref{2no:premin ruling notation}. Let $(x,y,z)$ with $x\leq y\leq z$ be the ordering of $(d_1,d_2,d_3)$, where as before $d_i=d(T_i)$ are discriminants of maximal twigs of $D$. By \ref{3lem:$W$ properties} we have $1>\delta>1-\frac{1}{|G|}\geq \frac{2}{3}$, where $|G|=d(\E)$, so $x\leq 4$ and $y\leq 11$.

\blem\label{4lem:x-y-|G| cases} One of the following cases occurs:\benum[(i)]

\item $(x,y)=(3,3)$ and $\E=[3]$,

\item $(x,y)=(2,3)$ and $\E\in \{[2,3], [3], [4], [5]\}$,

\item $(x,y)=(2,4)$ and $\E$ is either $[3]$ or $[4]$,

\item $(x,y)\in\{(2,5),(2,6)\}$ and $\E=[3]$.\eenum

\noin In particular, the two maximal twigs of $D$ corresponding to $x$ and $y$ belong to \\$\cal L=\{[2], [2,2], [2,2,2], [2,2,2,2], [2,2,2,2,2], [3], [4], [5], [6], [2,3], [3,2]\}$.\elem

\begin{proof} Suppose $z\leq 41$. Given an upper bound for $z$ there are finitely many possible weighted dual graphs of $D$. We used a computer program, which showed that for $x\leq 4$, $y\leq 11$, $z\leq 41$ conditions  \ref{1prop:basic S' properties}(iv), \ref{1lem:eps properties}(iii)-(iv), \ref{1lem:if eps<2}, \ref{2lem:Zar and others for D fork}  and \ref{3lem:$W$ properties}(iii) are satisfied only in three cases:\benum[(i)]

\item $b=1$, $T_1=[2]$, $T_2=[4]$, $T_3=[(8),4]$ and $\E=[4]$,

\item $b=2$, $T_1=[2]$, $T_2=[2,2]$, $T_3=[4,(6)]$ and $\E=[4]$,

\item $b=2$, $T_1=[2]$, $T_2=[2,2,2]$, $T_3=[3,3,(4)]$ and $\E=[4]$. \eenum

\noin These are included above, so we are done. Now suppose $z\geq 42$. For $x\geq 4$ we get $\frac{1}{z}>1-\frac{1}{|G|}-\frac{1}{2}\geq \frac{1}{6}$, which is impossible. For $x=3$ we have $\frac{1}{y}+\frac{1}{|G|}>\frac{2}{3}-\frac{1}{42}$, which gives $|G|=y=3$. Since $\delta<1$, for $x=2$ we have $y\geq 3$ and $\frac{1}{y}+\frac{1}{|G|}>\frac{1}{2}-\frac{1}{42}$, hence $y\leq 6$ and the bounds on $\E$ follow.\end{proof}

\bcor\label{4cor:there are two fibers} The ruling $f$ has two singular fibers and $\wt h=2$.\ecor

\begin{proof} By \ref{2cor:two fibers equations} $f$ has more than one singular fiber and it has at most three because $D$ is a fork. Each contains a unique $S_0$-component. Suppose it has three. Then $D^\dagger=D$ and since $x\leq 3$, for one of the singular fibers, say $F_1$, $F_1\cap D$ has at most two components, hence $F_1$ is a chain and $\Delta\cap F_1\neq \emptyset$. Then $\E=[2,3]$ and $\Delta\subseteq F_1=[2,1,2]$. It follows that the maximal twigs contained in other singular fibers of $f$ have more than two components, a contradiction with \ref{4lem:x-y-|G| cases}. Assume $\wt h\leq h$. Since $D$ is a fork, $\wt h \leq 2$. By \ref{2cor:two fibers equations} $\wt h=2$. \end{proof}

Let $T_1$, $T_2$ be the maximal twigs of $D$ contained respectively in the second and in the first branch of $F$. (The role of $T_i$'s is not symmetric because of this, that is exactly why we do not assume $d_1\leq d_2\leq d_3$, but use $x,y,z$ instead.) Clearly, they are also maximal twigs of $D^\dagger$ and $\varphi_f$ contracts the chain $H^\dagger+\wt Z_1+\wt Z_u$ to $T_3$.

We rewrite the equations of \ref{2prop:pairs} for two fibers. Put $\alpha=n+\epsilon+K\cdot E-4$, then $h=3+\alpha$ and $0\leq \alpha\leq n$. Put $\binom{\wt{\un c}_1}{\wt {\un p}_1}=\binom{\wt c}{\wt p}$, $\binom{\un c_1}{\un p_1}=\binom{c}{p}$ and $\binom{\un c_{h-1}}{\un p_{h-1}}=\binom{c'}{p'}$. Since $T_1$ is a chain, we have $\binom{\un c_2}{\un p_2}=\binom{\un c_3}{\un p_3}=\ldots=\binom{\un c_{h-2}}{\un p_{h-2}}=\binom{c'}{c'}$. Recall that $\rho=\kappa C\cdot E+c_h'C\cdot E+c_h'$. We have $\rho=\kappa^2$ for $\Delta'=0$ and $\rho=\frac{1}{2}(\kappa^2+1)$ for $\Delta'=[2]$, analogously for $\wt \rho$. In any case $\rho\leq \kappa^2$ and $\wt \rho^2\leq \wt \kappa^2$ (in fact these bounds hold in general, which can be shown by a straightforward computation). Recall that $\kappa,\wt \kappa\geq 2$ by \ref{1prop:basic S' properties}(ii). We have $d=c\kappa=\wt c\wt \kappa$, so we can write \eqref{eq:pairs-general I} as:

\begin{align}dn+\gamma-2&=\kappa (p+\alpha c'+p')+\wt \kappa \wt p .\label{eq:I}\end{align}

\noin Multiplying the above equation by $d$ and subtracting \eqref{eq:pairs-general II} we obtain:

\begin{align}d(\gamma-2)-\gamma&=\kappa^2(c-c')(\alpha c'+p')-\rho -\wt \rho.\label{eq:III}\end{align}

\bsrem Knowing the dual graph of $Z_l$ it is easy to determine ${c}/{c'}$ and ${p}/{c'}$. One has ${c}/{c'}=d(G+Z_u)=d(Z_l)$ and ${p}/{c'}=d(Z_u)=d(Z_l)-d(Z_l-Z_{ll})$ (cf. Appendix of \cite{KR-C*onC3}).\esrem

\brem\label{4rem:F determines wtF+H} For a fixed weighted dual graph of $F$ there are finitely many possible weighted dual graphs of $\wt F+H$. \erem

\begin{proof} If the (weighted) dual graph of $F$ is known then we know $c,p,c',p'$. The equation \eqref{eq:I} gives $$n(c-c')+\frac{\gamma-2}{\kappa}=p+(\epsilon+K\cdot E-4)c'+p'+\frac{\wt \kappa\wt p}{\kappa},$$ so $n(c-c')<p+p'+c\leq 2c$, hence $n<2+\frac{2c'}{c-c'}\leq 4$. Since now $\alpha$ is bounded, it is enough to bound $\kappa$, because then $d, \rho$, and hence $\wt c,\wt p,\wt \kappa, \wt \rho$ are bounded. We have $\wt c\wt \kappa=c\kappa$, so $\wt \kappa|c\cdot gcd(\kappa,\wt \kappa)$. By \eqref{eq:I} $gcd(\kappa,\wt \kappa)|\gamma-2$ and since $\gamma-2\in\{1,2,3\}$, we get $\wt \kappa|c(\gamma-2)$ and then $\wt \kappa\leq 3c$. Therefore $\wt \kappa$ and $\wt \rho$ are bounded. The coefficient of $\kappa$ in \eqref{eq:III} does not vanish, so \eqref{eq:III} is a nontrivial polynomial equation for $\kappa$ of degree at most two, so we are done.\end{proof}

\blem\label{4lem:T1 or T2 not on list} $d_1\leq 6$ if and only if $d_2>6$. \elem

\begin{proof} By \ref{4lem:x-y-|G| cases} $d_1\leq 6$ or $d_2\leq 6$. Suppose $d_1\leq 6$ and $d_2\leq 6$. Clearly, having the dual graph of $T_1$, there are only finitely many possibilities for the dual graphs of $T_1+C+\Delta'$, in each case $Z_1^2$ is determined. On the other hand, $T_2=Z_l$ and $(G+Z_u)^t$ are \emph{adjoint chains} (cf. \cite[4.7]{Fujita}), i.e. $e(G+Z_u)=1-e(\wt Z_l)$, so the dual graph of $G+Z_u$ is determined by $T_2$. Then by \ref{4rem:F determines wtF+H} there is finitely many possibilities for the dual graphs of $\wt F+H$. We use a computer program which for given $F$ (in terms of $(c,p,c',p')$) computes possible $(\gamma, n, \kappa, \rho, \wt\kappa, \wt c, \wt p, \wt \rho)$ using the algorithm sketched in \ref{4rem:F determines wtF+H} and checks if \eqref{eq:I} and \eqref{eq:III} can be satisfied. In each case (there may be many solutions) the maximal twig $T_3$ is determined and the program returns only these, for which conditions $\delta+\frac{1}{|G|}>1$, \ref{1lem:eps properties}(iii)-(iv), \ref{2lem:Zar and others for D fork}, $\sqrt{-d(D)/d(\E)}\in \Z$ and \ref{1lem:if eps<2} hold, these are:\benum[(i)]

\item $(n,\gamma,\kappa,\wt\kappa)=(1,4,4,2),\binom{c}{p}=\binom{4}{1},\binom{c'}{p'}=\binom{1}{1}, \binom{\wt c}{\wt p}=\binom{8}{5}$; $b=2$, $T_1=[2]$, $T_2=[(3)]$, $T_3=[3,3,(4)]$,

\item $(n,\gamma,\kappa,\wt\kappa)=(1,4,4,2),\binom{c}{p}=\binom{4}{3},\binom{c'}{p'}=\binom{1}{1}, \binom{\wt c}{\wt p}=\binom{8}{1}$; $b=1$, $T_1=[2]$, $T_2=[4]$, $T_3=[(8),4]$,

\item $(n,\gamma,\kappa,\wt\kappa)=(2,4,4,2),\binom{c}{p}=\binom{2}{1},\binom{c'}{p'}=\binom{1}{1}, \binom{\wt c}{\wt p}=\binom{4}{3}$; $b=2$, $T_1=[2,2]$, $T_2=[2]$, $T_3=[4,(6)]$.\eenum

In cases (i) and (ii) we have $-d(D)/d(\E)=4$ and $gcd(c,\wt c)=4$, in case (iii) $-d(D)/d(\E)=1$ and $gcd(c,\wt c)=2$. By \ref{2cor:two fibers equations}(iii) this is a contradiction. \end{proof}

We are ready to finish the proof of our main result.

\begin{proof}[Proof of the Main Theorem \ref{thm:main result}] As before, let $S'$ be a singular $\Q$-homology plane and let $S_0$ be its smooth locus. Suppose $\ovk(S')=-\8$ and $\ovk(S_0)=2$. With the notation as above by \ref{4lem:T1 or T2 not on list} and \ref{4lem:x-y-|G| cases} $T_3\in\cal L$. We first prove that $f$ is almost minimal. Suppose not. Then by \ref{2prop:pre-minimal ruling exists} $\wt \Delta=0$ and $\varphi_f$ contracts $H^\dagger+\wt Z_1$, where $H^\dagger=\wt Z_u+\wt G+H+G+Z_u$. Furthermore, $\varphi_f$ touches $\wt Z_1$ once and $Z_1$ $x$ times, where $x=1-\wt Z_{lu}^2\geq 4$. It follows that $n=1$, $\wt Z_1^2=-2$ and $Z_1^2=\wt Z_{lu}^2-b-1$. For a given weighted dual graph of $T_3$ the dual graph of $\wt G+\wt Z_u$ is determined uniquely. Indeed, $\wt G+\wt Z_u$ and $\wt Z_l^t$ are adjoint chains, so $e(\wt G+\wt Z_u)=1-e(\wt Z_l)$. Similarly, $e(G+Z_u)=1-e(Z_l)$. By the properties of $\varphi_f$ the chain $\wt C+\wt Z_1+H^\dagger$ has zero discriminant, so the snc-minimalization of $\wt G+\wt Z_u+\wt C$ is adjoint to $(G+Z_u)^t$, and hence has the same weighted dual graph as $Z_l$. Therefore $\wt Z_l$ determines the weighted dual graph of $H^\dagger+Z_1+Z_l$. Note that since $Z_1$ is touched more than once, $\wt Z_1+\wt Z_u$ cannot consist of $(-2)$-curves, so $\#T_3>1$. We now rule out the remaining cases.

\bca $T_3=[3,2]$.

We have $\wt Z_l=[3,3]$, so $\wt G+\wt Z_u=[2,3,2]$ and hence $d(T_2)=d([2,3,2,2,1])=d([2,2])=3$. Then $(x,y)=(3,5)$ by \ref{4lem:T1 or T2 not on list} and this contradicts \ref{4lem:x-y-|G| cases}. \eca

\bca $T_1=[2,3]$.

We have $\wt Z_l=[2,4]$, so $\wt G+\wt Z_u=[3,2,2]$ and hence $T_2$ is a minimalization of $[3,2,2,2,1]$, which is $[2]$. Then $(x,y)=(2,5)$, so $\E=[3]$ by \ref{4lem:x-y-|G| cases}. We have $\binom{\wt c}{\wt p}=\binom{7}{3}$ and $\binom{c}{p}=\binom{2c'}{c'}$, so $\kappa|d=7\wt \kappa$ and $gcd(\kappa,\wt \kappa)|\gamma-3$, hence $\kappa=7$ and $\wt \kappa=2c'$. However, \eqref{eq:I} gives $7p'=c'+1$ and then \eqref{eq:III} implies that $3(c')^2-7c'-46=0$, a contradiction with $c'\in \N$. \eca

\bca $T_1=[(k)]$ for some $k\in \{2,3,4,5\}$.

We have $\wt Z_l=[(k-1),3]$, so $\wt G+\wt Z_u=[k+1,2]$ and hence $T_2$ is a minimalization of $[k+1,2,2,1]$, which is $[k]$. Then by \ref{4lem:T1 or T2 not on list} $T_1\not\in \cal L$, so $(x,y)=(k,k+1)$ and we get $k=2$ by \ref{4lem:x-y-|G| cases}. We have $\binom{\wt c}{\wt p}=\binom{5}{2}$ and $\binom{c}{p}=\binom{2c'}{c'}$. Then $5\wt\kappa=d=2c'\kappa$, so by \eqref{eq:I} $\kappa c'(\alpha-1)=\gamma-2-\kappa p'-2\wt\kappa$. The left hand side is negative, so $\alpha=0$, i.e. $K\cdot E+\epsilon=3$. Suppose $\gamma=3$. By \eqref{eq:I} $gcd(\kappa,\wt \kappa)=1$, so $\kappa=5$. We get $c'=5p'-1$ and then \eqref{eq:III} implies $(c')^2-5c'+3-\rho=0$. For $\kappa=5$ we get $\rho=25$ or $\rho=13$, a contradiction with $c'\in \Z$. Thus $\gamma=4$ and now $gcd(\kappa,\wt \kappa)|2$, so $\kappa\in\{2,5,10\}$. We check that \eqref{eq:I} and \eqref{eq:III} lead to a contradiction for $\kappa\neq 2$ and for $\kappa=2$ give $\binom{c'}{p'}=\binom{25}{6}$. Then $T_1=[(3),7,(6)]$ and $b=2$, hence $d(D)=-25$, a contradiction with \ref{2cor:two fibers equations}(iii). \eca

Thus $f$ is almost minimal. Suppose $n>1$. Then $D^\dagger=D$ and $\wt h\geq 2$, so $\#T_3\geq 5$ and in fact $T_3=[(5)]$ because $T_3\in \cal L$. We get $\wt G+\wt Z_u=[2]$ and $G+Z_u=[2]$, so $\binom{c}{p}=\binom{2c'}{c'}$ and $\binom{\wt c}{\wt p}=\binom{2}{1}$, hence $\wt \kappa=\frac{d}{\wt c}=c'\kappa$. By \eqref{eq:I} we get $1<\kappa|\gamma-2$, so $\gamma\neq 3$ and hence $\Delta=0$. Then by \eqref{eq:III} $\kappa|\gamma$, so $\kappa=2$ and $\E=[4]$. We get $\alpha=1$ and then \eqref{eq:I} gives $p'=c'+1$, which contradicts $p'\leq c'$.

Since $f$ is almost minimal, $\varphi$ does not contract $\wt Z_1$, so $\#T_3\geq 2$. Moreover, if $\#T_3=2$ then $\#\wt Z_l=1$, so $\wt G+\wt Z_u$ consists of $(-2)$-curves and since $\varphi_f$ has to contract $G$, we see that $\wt Z_1$ is touched at least twice by $\varphi_f$. The latter shows that if $\#T_3=2$ then $\wt Z_1^2\leq -4$, which contradicts $\#\Delta\leq 1$. Therefore $T_3=[(k)]$ for some $k=3,4,5$.

By \ref{4lem:x-y-|G| cases} $\E=[4]$ or $\E=[3]$. In particular, $\alpha=0$ and $\Delta=0$. The latter yields $\wt Z_1^2=-2$. Now $\wt Z_l$ consists of $(-2)$-curves, so $\wt Z_u=0$. Let's write $\wt Z_l=[(s)]$ and $\wt G=[s+1]$ for some $s\geq 1$. Since $\varphi_f$ does not contract $\wt Z_1$, it cannot contract $\wt G$. This gives $s\geq 2$, as $n=1$. Suppose $G\neq [2]$. Then $\#T_3\leq 5$ implies $s=2$, $Z_u=0$ and $G=[3]$, so $d_2=3$. By \ref{4lem:T1 or T2 not on list} we get $(x,y)=(3,6)$, a contradiction with \ref{4lem:x-y-|G| cases}. Thus $G=[2]$, so $\varphi_f$ touches $\wt G$ at least twice, which gives $s\geq 3$. Now $k\leq 5$ implies $s=3$ and $Z_u=0$. By \ref{4lem:x-y-|G| cases} $\E=[3]$. We have $\binom{\wt c}{\wt p}=\binom{4}{1}$ and $\binom{c}{p}=\binom{2c'}{c'}$. Then $4\wt\kappa=d=2c'\kappa$ and $gcd(\kappa,\wt \kappa)=1$, so $\kappa=2$. Now \eqref{eq:I} gives $c'=2p'-1$, so by \eqref{eq:III} $(c')^2-2c'=1$, a contradiction. \end{proof}

\med
\bibliographystyle{amsalpha}
\bibliography{bibl}

\providecommand{\bysame}{\leavevmode\hbox to3em{\hrulefill}\thinspace}
\providecommand{\MR}{\relax\ifhmode\unskip\space\fi MR }
\providecommand{\MRhref}[2]{%
  \href{http://www.ams.org/mathscinet-getitem?mr=#1}{#2}
}
\providecommand{\href}[2]{#2}
\begin{thebibliography}{BHPVdV04}

\bibitem[BHPVdV04]{BHPV}
Wolf~P. Barth, Klaus Hulek, Chris A.~M. Peters, and Antonius Van~de Ven,
  \emph{Compact complex surfaces}, second ed., Ergebnisse der Mathematik und
  ihrer Grenzgebiete. 3. Folge. A Series of Modern Surveys in Mathematics
  [Results in Mathematics and Related Areas. 3rd Series. A Series of Modern
  Surveys in Mathematics], vol.~4, Springer-Verlag, Berlin, 2004.

\bibitem[Bri68]{Brieskorn}
Egbert Brieskorn, \emph{Rationale {S}ingularit\"aten komplexer {F}l\"achen},
  Invent. Math. \textbf{4} (1967/1968), 336--358.

\bibitem[CS10]{fake_proj_spaces}
Donald~I. Cartwright and Tim Steger, \emph{Enumeration of the 50 fake
  projective planes}, C. R. Math. Acad. Sci. Paris \textbf{348} (2010),
  no.~1-2, 11--13.

\bibitem[Fuj79]{Fuj-Zar}
Takao Fujita, \emph{On {Z}ariski problem}, Proc. Japan Acad. Ser. A Math. Sci.
  \textbf{55} (1979), no.~3, 106--110.

\bibitem[Fuj82]{Fujita}
\bysame, \emph{On the topology of noncomplete algebraic surfaces}, J. Fac. Sci.
  Univ. Tokyo Sect. IA Math. \textbf{29} (1982), no.~3, 503--566.

\bibitem[GM92]{GM-Affine-lines}
R.~V. Gurjar and M.~Miyanishi, \emph{Affine lines on logarithmic {${\bf
  Q}$}-homology planes}, Math. Ann. \textbf{294} (1992), no.~3, 463--482.

\bibitem[Hir53]{Hirzebruch-cyclic_singularities}
Friedrich Hirzebruch, \emph{\"{U}ber vierdimensionale {R}iemannsche {F}l\"achen
  mehrdeutiger analytischer {F}unktionen von zwei komplexen
  {V}er\"anderlichen}, Math. Ann. \textbf{126} (1953), 1--22.

\bibitem[Iit82]{Iitaka}
Shigeru Iitaka, \emph{Algebraic geometry}, Graduate Texts in Mathematics,
  vol.~76, Springer-Verlag, New York, 1982, An introduction to birational
  geometry of algebraic varieties, North-Holland Mathematical Library, 24.

\bibitem[Kob90]{Kob}
Ryoichi Kobayashi, \emph{Uniformization of complex surfaces}, K\"ahler metric
  and moduli spaces, Adv. Stud. Pure Math., vol.~18, Academic Press, Boston,
  MA, 1990, pp.~313--394.

\bibitem[Kor93]{Koras-A2}
Mariusz Koras, \emph{A characterization of {$\bold A\sp 2/\bold Z\sb a$}},
  Compositio Math. \textbf{87} (1993), no.~3, 241--267.

\bibitem[KR99]{KR-C*onC3}
Mariusz Koras and Peter Russell, \emph{{${\bf C}\sp *$}-actions on {${\bf C}\sp
  3$}: the smooth locus of the quotient is not of hyperbolic type}, J.
  Algebraic Geom. \textbf{8} (1999), no.~4, 603--694.

\bibitem[KR07]{KR-ContrSurf}
\bysame, \emph{Contractible affine surfaces with quotient singularities},
  Transform. Groups \textbf{12} (2007), no.~2, 293--340.

\bibitem[Lan03]{Langer}
Adrian Langer, \emph{Logarithmic orbifold {E}uler numbers of surfaces with
  applications}, Proc. London Math. Soc. (3) \textbf{86} (2003), no.~2,
  358--396.

\bibitem[Miy01]{Miyan-OpenSurf}
Masayoshi Miyanishi, \emph{Open algebraic surfaces}, CRM Monograph Series,
  vol.~12, American Mathematical Society, Providence, RI, 2001.

\bibitem[MS91a]{MiSu-hPlanes}
M.~Miyanishi and T.~Sugie, \emph{Homology planes with quotient singularities},
  J. Math. Kyoto Univ. \textbf{31} (1991), no.~3, 755--788.

\bibitem[MS91b]{MiSu-Qhp_with_C**}
Masayoshi Miyanishi and Tohru Sugie, \emph{{${\bf Q}$}-homology planes with
  {${\bf C}\sp {**}$}-fibrations}, Osaka J. Math. \textbf{28} (1991), no.~1,
  1--26.

\bibitem[MT84]{MiTs-PlatFibr}
Masayoshi Miyanishi and Shuichiro Tsunoda, \emph{Noncomplete algebraic surfaces
  with logarithmic {K}odaira dimension {$-\infty$} and with nonconnected
  boundaries at infinity}, Japan. J. Math. (N.S.) \textbf{10} (1984), no.~2,
  195--242.

\bibitem[Mum61]{Mumford}
David Mumford, \emph{The topology of normal singularities of an algebraic
  surface and a criterion for simplicity}, Inst. Hautes \'Etudes Sci. Publ.
  Math. (1961), no.~9, 5--22.

\bibitem[Pal08]{Palka-classification}
Karol Palka, \emph{Classification of singular {$\bold Q$}-acyclic surfaces with
  smooth locus of non-general type},
  \href{http://arxiv.org/abs/0806.3110}{arXiv:0806.3110} (2008).

\bibitem[Pal09]{Palka-k(S_0)=0}
\bysame, \emph{Exceptional singular {$\bold Q$}-homology planes},
  \href{http://arxiv.org/abs/0909.0772}{arXiv:0909.0772} (2009).

\bibitem[Pal10]{Palka-review}
\bysame, \emph{Recent progress in the geometry of {$\bold Q$}-acyclic
  surfaces}, \href{http://arxiv.org/abs/1003.2395}{arXiv:1003.2395} (2010).

\bibitem[PS97]{PS-rationality}
C.~R. Pradeep and Anant~R. Shastri, \emph{On rationality of logarithmic {$\bold
  Q$}-homology planes. {I}}, Osaka J. Math. \textbf{34} (1997), no.~2,
  429--456.

\bibitem[Rus80]{Russell2}
Peter Russell, \emph{Hamburger-{N}oether expansions and approximate roots of
  polynomials}, Manuscripta Math. \textbf{31} (1980), no.~1-3, 25--95.

\end{thebibliography}
\end{document}